\documentclass[a4paper,11pt]{article}
 
\usepackage{amsthm}
\usepackage[abbrvbib,preprint]{jmlr2e}
\RequirePackage{amsmath,amssymb,amsfonts,algorithm,bm,bbm,xcolor,graphicx,fancyhdr,lastpage,mathtools,comment} 
\usepackage[noend]{algpseudocode}
\usepackage[shortlabels]{enumitem}
\usepackage{cmap,lmodern} \usepackage[T1]{fontenc} 
\definecolor{navy}{RGB}{20,0,105}
\usepackage[capitalise]{cleveref}		
\usepackage[titletoc,title]{appendix}

\newcommand{\dif}{\mathop{}\!\mathrm{d}} 
  \newcommand{\du}{\dif u}     \newcommand{\dx}{\dif x}  \newcommand{\dy}{\dif y} \newcommand{\dz}{\dif z}      \newcommand{\dmu}{\dif \mu}

  \newcommand{\QQ}{\mathbb{Q}} \newcommand{\NN}{\mathbb{N}} \newcommand{\ZZ}{\mathbb{Z}} \newcommand{\RR}{\mathbb{R}}   \newcommand{\II}{\mathbbm{1}}  
\newcommand{\DD}{\mathbb{D}}
\newcommand{\Aa}{\mathcal{A}} \newcommand{\Bb}{\mathcal{B}} \newcommand{\Cc}{\mathcal{C}} \newcommand{\Dd}{\mathcal{D}}    \newcommand{\Hh}{\mathcal{H}} \newcommand{\Ii}{\mathcal{I}}           \newcommand{\Tt}{\mathcal{T}} \newcommand{\Uu}{\mathcal{U}} \newcommand{\Vv}{\mathcal{V}}  \newcommand{\Xx}{\mathcal{X}} \newcommand{\Yy}{\mathcal{Y}} 
\DeclarePairedDelimiter{\norm}{\lVert}{\rVert}
\DeclarePairedDelimiter{\abs}{\lvert}{\rvert}
\DeclarePairedDelimiter{\braces}{ \{ }{ \} }
\DeclarePairedDelimiter{\brackets}{(}{)}
\DeclarePairedDelimiter{\sqbrackets}{[}{]}

\DeclarePairedDelimiter{\ip}{\langle}{\rangle}
\DeclarePairedDelimiter{\ceil}{\lceil}{\rceil}
\DeclarePairedDelimiter{\floor}{\lfloor}{\rfloor}

\DeclareMathOperator*{\argmax}{argmax}
\DeclareMathOperator{\Var}{Var}

\DeclareMathOperator{\sign}{sign}

\DeclareMathOperator{\KL}{KL} 
\DeclareMathOperator{\FDP}{FDP}
\DeclareMathOperator{\FDR}{FDR}
\DeclareMathOperator{\postFDR}{postFDR}

\DeclareMathOperator{\mFDR}{mFDR}
\DeclareMathOperator{\mTDR}{mTDR}
\DeclareMathOperator{\TDR}{TDR}
\DeclareMathOperator{\diag}{diag}
\newcommand{\vphi}{\varphi}
\newcommand{\eps}{\varepsilon}

\DeclareMathOperator{\Markov}{Markov}
\DeclareMathOperator{\sep}{sep}

\newcommand{\lambdamax}{\lambda_{\textnormal{max}}}
\newcommand{\transpose}{\intercal}

\newcommand{\paragraphi}[1]{\emph{#1.}}

\theoremstyle{definition} 
\newtheorem{assumption}{Assumption} 
\newtheorem{condition}[assumption]{Condition}
\newtheorem*{assumption*}{Assumption*}

\newtheorem{manualassumptioninner}{Assumption}
\newenvironment{manualassumption}[1]{%
	\renewcommand\themanualassumptioninner{#1}%
	\manualassumptioninner
}{\endmanualassumptioninner}
\theoremstyle{remark} 
\newtheorem*{remark*}{Remark} 
\newtheorem*{remarks}{Remarks} \newtheorem*{examples}{Examples}
\Crefname{appsec}{Appendix}{Appendices}
\Crefname{assumption}{Assumption}{Assumptions} \Crefname{assumption*}{Assumption}{Assumptions}
\Crefname{equation}{}{}
\crefrangelabelformat{equation}{eqs.~(#3#1#4) to~(#5#2#6)}
\Crefname{enumi}{}{}
\newlist{lemenum}{enumerate}{1} 
\setlist[lemenum]{label=\alph*., ref=\arabic{theorem}\alph*}
\crefalias{lemenumi}{lemma} 
\newlist{thmenum}{enumerate}{1} 
\setlist[thmenum]{label=\Alph*., ref=\arabic{theorem}\Alph*}
\crefalias{thmenumi}{theorem} 
\Crefname{manualassumption}{Assumption}{Assumptions}
\Crefname{theorem}{Theorem}{Theorems}
\Crefname{example}{Example}{Examples}
\Crefname{lemma}{Lemma}{Lemmas}
\Crefname{lemma*}{Lemma}{Lemmas}
\Crefname{proposition}{Proposition}{Propositions}
\Crefname{condition}{Condition}{Conditions}
\renewcommand{\namecref}{\lcnamecref}

\algnewcommand{\Initialize}[1]{%
	\State \textbf{initialize}
	\parbox[t]{.8\linewidth}{\raggedright #1}
}
\algnewcommand{\Input}[1]{%
	\State \textbf{input}
	\parbox[t]{.8\linewidth}{\raggedright #1}
}
\algnewcommand{\Estimate}[1]{%
	\State \textbf{estimate}
	\parbox[t]{.8\linewidth}{\raggedright #1}
}

\algnewcommand{\Choose}[1]{%
	\State \textbf{choose}
	\parbox[t]{.8\linewidth}{\raggedright #1}
}

\algnewcommand{\Then}[1]{%
	\State \textbf{then}
	\parbox[t]{.8\linewidth}{\raggedright #1}
}

\algnewcommand{\Output}[1]{%
	\State \textbf{output}
	\parbox[t]{.8\linewidth}{\raggedright #1}
}

\newcommand{\given}{\,|\,}
\newcommand{\te}{\theta}
\newcommand{\si}{\sigma}
\newcommand{\al}{\alpha}

\newcommand{\EM}{\ensuremath}
\newcommand{\cC}{\EM{\mathcal{C}}}
\newcommand{\cF}{\EM{\mathcal{F}}}

\newcommand{\bef}{{\bf f}}
\newcommand{\beg}{{\bf g}}
\newcommand{\bev}{{\bf v}}
\newcommand{\bet}{{\bf T}}

\def\1{1\!{\rm l}}

\definecolor{blendedblue}{rgb}{0.2,0.2,0.7}

\jmlrheading{}{}{}{}{}{}{Kweku Abraham, Isma{\"e}l Castillo and Elisabeth Gassiat}
\ShortHeadings{Multiple Testing in Nonparametric HMMs}{Abraham, Castillo and Gassiat}
\firstpageno{1}

\title{Multiple Testing in Nonparametric Hidden Markov Models: An Empirical Bayes Approach}
\author{\name Kweku Abraham \email kweku.abraham@universite-paris-saclay.fr \\
	\addr Universit{\'e} Paris-Saclay, CNRS, Laboratoire de Math{\'e}matiques d'Orsay, \\ 91405 Orsay, France
	\AND 
\name Isma{\"e}l Castillo \email ismael.castillo@sorbonne-universite.fr\\
\addr Sorbonne Universit{\'e}, Laboratoire de Probabilit\'es, Statistique et Mod\'elisation,\\
4 Place Jussieu, 75005 Paris, France
	\AND
	\name Elisabeth Gassiat \email elisabeth.gassiat@universite-paris-saclay.fr \\
	\addr Universit{\'e} Paris-Saclay, CNRS, Laboratoire de Math{\'e}matiques d'Orsay, \\ 91405 Orsay, France
}

\editor{}

\begin{document}
	\maketitle
\begin{abstract}%
Given a nonparametric Hidden Markov Model (HMM) with two states, the question of constructing efficient multiple testing procedures is considered, treating one of the states as an unknown null hypothesis. A procedure is introduced, based on nonparametric empirical Bayes ideas, that controls the False Discovery Rate (FDR) at a user--specified level. Guarantees on power are also provided, in the form of a control of the true positive rate. One of the key steps in the construction requires supremum--norm convergence of preliminary estimators of the emission densities of the HMM. We provide the existence of such estimators, with convergence at the optimal minimax rate, for the case of a HMM with $J\ge 2$ states, which is of independent interest.

\end{abstract}
\begin{keywords}
efficient multiple testing, hidden Markov models, false discovery rate, true discovery rate, supremum norm estimation, minimax rate
\end{keywords}
	
\section{Introduction}\label[section]{sec:Introduction}

\subsection{Aim of the Paper}\label[section]{sec:Aim}
We consider the problem of multiple testing in a hidden Markov model (HMM) setting. Given data $(X_i : i\leq N)$ whose distribution is governed by an unknown categorical variable $\theta=(\theta_i: i\leq N)$ drawn from a Markov chain, one seeks to test the null hypotheses $H_{0,i}: \theta_i=0$ against the alternatives $H_{1,i}: \theta_i\neq 0$ simultaneously for $i=1,\dots N$. In seeking procedures with optimal properties with respect to multiple testing--measures of risk, for example with controlled False Discovery Rate (FDR) and maximal `power' as measured by the True Discovery Rate (TDR), it is natural to consider thresholding based on the probabilities of the $\theta_i$'s being zero conditional on the observations $X_i$, $i=1,\ldots,N$ (see \cref{sec:TestingProcedureDefinition}).
The conditional probabilities are simply posterior probabilities in the Bayesian world, and smoothing probabilities in the latent variables vocabulary. They will (mainly) be called  $\ell$-value's in this work.
 
A first such procedure is one that rejects all coordinates whose $\ell$--value is below a user-specified level $t$, see e.g.\ \cite{ETST2001, Efron2007}. From the Bayesian point of view, this is directly related to the Bayes factor for testing the individual coordinate. The procedure we consider in this paper is still one based on $\ell$--value thresholding, but with a data--dependent threshold chosen in such a way that the expected false discovery rate conditional on the data is equal or very close to $t$ \citep{muelleretal04, SC09}; this typically yields an FDR close to $t$, and, being less `conservative' than the previous procedure, enjoys certain optimality properties. There are other alternatives, such as so-called $q$--value procedures \citep{Sto03}, that are based on controlling `marginal' versions of the FDR.

Of course, it is rare that the model parameters are known in practice, so that instead of the `oracle' procedures described above (so called because the $\ell$--values depend on the parameters), calculations are based on first estimating these in the chosen modeling: the `empirical' Bayes method. We consider here a \emph{nonparametric} HMM setting, with unknown parameters corresponding to emission densities and to characteristics of the underlying Markov chain. The key question addressed in this paper is this: 
\emph{what is the impact of the estimation step on the FDR?} 
More precisely, one would like to understand whether the discussed thresholding procedures still (asymptotically) maintain multiple testing optimality properties when the parameters are estimated. 
We also note that, being in a nonparametric setting, the loss function chosen to measure the quality of estimation may have more influence over the plug--in operation compared to the parametric situation.
Our main contributions can be summarised as follows.
\begin{itemize}
	\item	Our first main results, \cref{thm:BFDRcontrol,thm:BFNRoptimality},
	show theoretically that  in the nonparametric HMM setting  an empirical Bayesian procedure attains the target FDR level and enjoys TDR optimality. 
	The proofs of these two theorems are partly based on a result in \cite{dCGlC17}, which shows how control of plug--in estimators propagates to give control of $\ell$--value errors.
	A key step is to have good \emph{supremum--norm} estimators, in contrast to the $L^2$--norm estimators previously found in the literature. 
	\item
	Our second main results, which are both key to obtaining the first and also of independent interest, concern supremum--norm estimation of emission densities in nonparametric HMMs. We provide estimators, and prove in \cref{thm:DiscreteEstimation,thm:LinftyEstimation}   that the supremum--norm risk of these estimators achieves the parametric convergence rate $N^{-1/2}$ for discrete observations (where the set of possible values is countable), and the convergence rate  $(N/\log N)^{-s/(2s+1)}$, familiar from the classical i.i.d.\ density estimation setting and also proved to be optimal in the HMM context (see Proposition \ref{prop:MinimaxLowerBound}), for H\"older densities with regularity $s$.
\end{itemize}

Our key question connects with the frequentist analysis of the behaviour of empirical Bayes procedures, a topic currently under rapid development: we briefly review some connections at the end of \cref{sec:setting}. It has previously been considered in an i.i.d.\ setting in \cite{SC07}, in a graph setting (with a $q$--value based procedure) in \cite{RRV19}, and, most pertinently, in a \emph{parametric} HMM setting in \cite{SC09}, wherein it is argued that first estimating the model parameters leaves the FDR asymptotically unchanged. 
Parametric modeling of HMMs is known, however, to lead to poor results in many applications, as shown for instance in \cite{YPRH11}. 
We draw attention also to the extensive  simulations conducted and discussed in \cite{WSZ19} for real valued observations, and in \cite{SW20} for count data. These latter two works demonstrate empirically that the FDR and TDR are badly impacted by parametric modeling in case of misspecification, while nonparametric empirical Bayes methods as considered here closely match the optimal behaviour of oracle $\ell$--value procedures. 

A further advantage of modelling the HMM densities nonparametrically is that it ensures our results allow for fairly arbitrary distributions under the null hypothesis. 
In contrast, many common multiple testing procedures -- including the original Benjamini--Hochberg procedure -- assume that the null distribution is known. One can of course adjust such procedures to use an estimated null hypothesis, but there are so far only a few settings in which it has been proved that this plug-in step has no negative effect on the desired properties of the procedures. We refer to the recent work by \cite{RV19} for more discussion concerning this issue.

Finally, we note that as well as enabling the plug-in results which yield control of the FDR, estimating the emission densities in terms of the supremum norm is useful in its own right. Indeed, practically speaking, results of this type justify that plots of density estimators will be visual close to the original density. Such estimators can also be helpful for identifying change points, estimating level sets, and constructing confidence bands for uncertainty quantification. 

\subsection{Context}\label[section]{sec:Context}
Let us now place these results in the broader multiple testing and HMM contexts.

\paragraphi{Multiple testing} The problem of identifying relevant variables among a large number of possible candidates is ubiquitous with high dimensional data: indeed, multiple testing methods are very popular in the analysis of genomic data, in astrostatistics, and in imaging, to name just a few practical applications. Since the seminal work of \cite{BH1995}, controlling the FDR has been the goal of much of the extensive literature on the subject.

Early works tended to assume i.i.d.\ data. \cite{Efron2007Correlations} noted that ignoring dependence and using methods designed for FDR control with independent data could result in either too conservative or too liberal procedures, showing that dependence must carefully being taken into account.  A number of works, including those of \cite{BY01}, \cite{Farcomeni07}, \cite{FinnerDickhausRoters07} and \cite{Wu08}, have shown that under certain assumptions on the dependence structure, some multiple testing procedures designed for independent case (such as the step-up Benjamini--Hochberg procedure) still control the FDR below a given target level. Such procedures, although having guaranteed FDR even under dependence, may suffer from being too conservative. 

The control of power in dependent data settings is less developed. Some works in this direction include those of \cite{XCML11} and of \cite{HellerRosset20} which consider the `general two group model', wherein the $\theta_i$'s are independent and identically distributed, but for each $i$  the distribution of $X_{i}$ given $\theta$ may depend on the whole vector $\theta$ and not only on $\theta_i$. In some settings, such as with genetic data, allowing for the $\theta_i$'s themselves to be dependent can however be more natural, and the HMM model for $X$ considered here allows for a natural local structure of $\theta$ -- while still remaining tractable -- by modelling it as a Markov chain.

\paragraphi{Hidden Markov models}
HMMs have been widely used for applications as varied as speech modelling, computational finance and gene prediction since works of Baum, Petrie and coauthors introduced practical algorithms and proved parametric estimation rates in a discrete data setting \citep{Petrie67,BP66,BPSW70}. Later works, including those of \cite{BRR98} and of \cite{DM01}, extended these proofs to allow parametric modelling of the emission distributions. 

Recently, \cite{GCR16} opened the possibility that consistency holds also when the emission densities are modelled nonparametrically by proving identifiability under mild conditions. \cite{AHK12} introduced in the parametric case a spectral method which was then generalised in \cite{dCGL16} and \cite{Lehericy18} to indeed give consistency at a usual rate in the nonparametric setting. These nonparametric works however focus on $L^2$--estimation, and do not immediately generalise to give rate-optimal supremum norm estimation: indeed, attempting to apply a typical wavelet method of estimating individual coefficients at a parametric rate and aggregating, one runs into an alignment issue arising from the fact that the emission densities are identifiable only up to a permutation. An insight of the current work is that returning to the spectral method and using a kernel based estimator allows this issue to be bypassed.

\subsection{Setting}\label[section]{sec:FrameworkAndAssumptions}
\label[section]{sec:setting}

Consider a hidden Markov model (HMM), in which the observations $X=(X_n)_{n\leq N}$ satisfy
\begin{equation}
\label[equation]{eqn:def:model}
\begin{split}
X_n\mid \theta&\sim f_{\theta_{n}}, \quad 1\leq n \leq N,\\
\theta=(\theta_n)_{n\leq N} &\sim \Markov(\pi,Q),\end{split}
\end{equation} 
and, conditional on $\theta$, the entries of $X$ are independent. The vector $\theta$ of `hidden states' takes values in $\braces{0,1}^N$ (we will later also consider the case where $\theta$ takes values in $\braces{1,\dots,J}^N$ for some $J\geq 2$) and $\Markov(\pi,Q)$ denotes a Markov chain of initial distribution $\pi=(\pi_0,\pi_1)$, 
 and $2\times 2$ transition matrix $Q$.  
The `emission densities' $f_0,f_1$ are probability densities with respect to some common dominating measure $\mu$ on a measurable space $\Xx$. For simplicity we will assume that $\mu$ is either Lebesgue measure on $\RR$ or counting measure on $\ZZ\subset \RR$; our results adapt straightforwardly to the $d$--dimensional setting, and in principle versions should hold for more general measure spaces (see the discussion in \cref{sec:Extensions}). We use $H=\braces{Q,\pi,f_0,f_1}$ to denote a generic set of parameters for the HMM. We denote by $\Pi_H$ the law of $(X,\theta)$ in \cref{eqn:def:model}, and by extension also the marginal laws of $X$ and $\theta$. 
 We write $E_H$ to denote the expectation operator associated to $\Pi_H$. 

The goal of multiple testing is to provide a procedure $\vphi=\vphi(X)$ which identifies well for which $i$ we have signal ($\theta_{i}\not=0$).  We will measure the performance of $\vphi$ through the false discovery rate (FDR) and the true discovery rate (TDR). Defining the false discovery proportion (FDP) at $\theta$ as
\begin{equation} \label[equation]{eqn:def:FDP}
\FDP_\theta(\vphi):= \frac{\sum_{i=1}^N \II\braces{\theta_{i}=0, \vphi_i=1}}{1\vee \brackets[\big]{\sum_{i=1}^N \vphi_i}},
\end{equation}
the FDR at $\theta$ is given by
\begin{equation}\label[equation]{eqn:def:FDR}
\FDR_\theta(\vphi):= E[\FDP_\theta(\vphi(X)) \mid \theta].
\end{equation}
We consider the average false discovery rate for $\theta$ generated according to the `prior' law $\Pi_{H}$:
\begin{equation}\label[equation]{eqn:def:BFDR}
\FDR_{H}(\vphi):=E_{\theta\sim \Pi_{H}} \FDR_\theta(\vphi)\equiv E_{(X,\theta) \sim \Pi_{H}} \FDP_\theta(\vphi),
\end{equation} 
and we define the `posterior FDR' as the FDR obtained by drawing $\theta$ from its posterior 
\begin{equation}\label[equation]{eqn:def:postFDR1} \postFDR_{H}(\vphi) =\postFDR_H(\vphi;X) := E_H [\FDP_\theta(\vphi) \mid X].
\end{equation}
The true discovery rate is defined as the expected proportion of signals which are detected by a procedure:
\begin{equation}\label{eqn:def:TDR} \TDR_H(\vphi) = E_H \sqbrackets[\Big]{\frac{\sum_{i=1}^N \II\braces{\theta_i=1,\vphi_i=1}}{1\vee \brackets{\sum_{i=1}^N \theta_i}}}. \end{equation}

\paragraphi{Bayesian formulation and latent variable formulation} Let $P_0$ denote the 
``true'' distribution of the data $X$ arising from model \eqref{eqn:def:model}. 
If, in \eqref{eqn:def:model}, the distribution of $\te$ is interpreted as ``prior'' distribution (it is of course an ``oracle prior'', as $\pi,Q$  are components of the unknown ``true'' parameter $H=(\pi,Q,f_0,f_1)$), the distribution of $X=(X_n)_{n\le N}$ in the (oracle) Bayesian setting is simply the true distribution $P_0$. Of course, one may also avoid the Bayesian vocabulary and simply view model \eqref{eqn:def:model} as a latent variable model: under such point of view, $\ell$--values are known as smoothing probabilities and $\te\given X$ is simply a conditional distribution. We find it convenient to nevertheless use Bayesian terminology. Partially this is in accordance with classical decision theory, wherein Bayesian terminology is commonly used for describing optimal classifiers (indeed, as \cite{Sto03} observed, ``classical classification theory seems to be a bridge between Bayesian modeling and hypothesis testing''). It is also helpful preparation for considering a setting where $\theta$ is fixed and non-random, as dicussed next. 

\paragraphi{Connection with frequentist analysis of Bayesian procedures} Recent years have seen notable progress on providing frequentist validations of the use of posterior distributions for inference, with most results concerning the estimation task, and more recently also uncertainty quantification and confidence sets \citep{GVdV_book}. One can consider using the HMM model \eqref{eqn:def:model} not because one believes $\theta$ is genuinely random with a Markov structure, but rather as a way to model some block structure of a fixed true $\theta$, wherein neighbour coordinates of $X$ have a higher chance of coming from the same distribution. The first results in this spirit in a multiple testing setting were obtained recently for sparse sequences (without block structure) in \cite{CR20}. We plan to investigate in further work the Bayesian procedure studied in this paper for structured sequences of fixed $\theta$ where the HMM modeling will then be a Bayesian prior.   

We also note that the results we obtain below still constitute a (partial) frequentist Bayes validation, in the following sense. Consider a standard Bayesian approach where $\theta$ is viewed as parameter and given a HMM prior, but not the other parameters $(f_0,f_1,\pi, Q)$, which are estimated separately. Then Theorems \ref{thm:BFDRcontrol} and \ref{thm:BFNRoptimality} below prove that if the true (frequentist) data generating distribution is some nonparametric HMM, then the empirical Bayes procedure derived from the posterior on $\theta$ behaves consistently from the multiple testing viewpoint: its FDR is controlled with optimality guarantees on the TDR. It is a less strong frequentist analysis than under an arbitrary fixed $\theta_0$, but it validates the frequentist use of the procedure assuming that the data comes from some (fairly arbitrary) non-parametric HMM: this still allows one to capture many typical signals with varied latent densities. 

\subsection{Outline of the Paper}\label[section]{sec:outline}

In Section \ref{sec:EmpiricalBayesProcedure} we introduce 
our multiple testing procedure and establish its asymptotic performance in 
\cref{thm:BFDRcontrol,thm:BFNRoptimality}.
Section \ref{sec:SupNormEstimation} is devoted to the estimation of the emission densities, with asymptotic supremum norm control established in  \cref{thm:DiscreteEstimation,thm:LinftyEstimation}. 
We also give in \cref{prop:MinimaxLowerBound} a lower bound for the estimation of  H\"older emission densities with regularity $s$ in the HMM context.
Finally, \cref{prop:ExistenceOfEstimators} gives examples of how to overcome the `label switching' issue, present in the HMM setting as for mixture models, in order to know which estimator corresponds to the null state and which to the alternative. This allows us to avoid the assumption, common to many multiple testing methods, that the distribution of the data under the null is known.

In Section \ref{sec:DiscussionOfMainResults}, we provide a detailed discussion of our assumptions and comparisons of our results with the literature. We also explain how the rates of convergence of our emission densities estimators can be understood as minimax rates of convergence in supremum norm. 

Proofs of the main theorems are given in \cref{sec:Proofs}. Intermediate results useful for these proofs are given in \cref{sec:LemmasForMultipleTesting,sec:lemmasforLinftyEstimation}. \Cref{sec:ProofOfLowerBound} gives a proof of a minimax lower bound.
For the reader's convenience, the notation introduced throughout the paper is gathered in \cref{sec:Notation}.

\section{The Empirical Bayesian Procedure}\label[section]{sec:EmpiricalBayesProcedure}
\subsection{Definition}\label[section]{sec:TestingProcedureDefinition}
We analyse an empirical Bayesian approach to the multiple testing problem, based on thresholding by the posterior (smoothing) probabilities, here called the  `$\ell$--values' and also known in the literature as the `local indices of significance' \citep{ETST2001, Efron2007,SC09}:
\begin{equation}
\label[equation]{eqn:def:Lvals}
\ell_i(X)\equiv \ell_{i,H}(X)=\Pi_{H}(\theta_i=0 \mid X).
\end{equation}

In the `oracle' setting (where the parameter $H$ is known), it is well known that the optimal (weighted) classification procedure is an $\ell$--value thresholding procedure; that is, it is $\vphi_{\lambda,H}$ for some $\lambda$, where 
\begin{equation} \label[equation]{eqn:def:PhiLambda} 
	\vphi_{\lambda,H}(X)=(\II\braces{\ell_{i,H}(X)<\lambda})_{i\leq N}.
\end{equation}
It has been shown in \cite{SC09} that this class of procedures (possibly with data-driven thresholds) is also optimal in a multiple testing sense, in that a procedure making false discoveries at a pre-specified rate and maximising a suitable notion of the multiple testing power is necessarily an $\ell$--value thresholding procedure.

The FDR is the expectation of the posterior FDR, so that using the latter (which is observable) to choose the threshold is a natural approach.
When the parameter $H$ is unobserved, we use an estimator $\hat{H}=(\hat{Q},\hat{\pi},\hat{f}_0,\hat{f}_1)$ instead (to be constructed later), and so we
are led to the procedure $\vphi_{\hat{\lambda},\hat{H}}$, where
\begin{equation}
\label[equation]{eqn:def:hatlambda}
\hat{\lambda}=\hat{\lambda}(\hat{H},t):=\sup\braces{\lambda : \postFDR_{\hat{H}}(\vphi_{\lambda,\hat{H}})\leq t}.\end{equation} 
\smallskip
We also note an alternative characterisation of the threshold $\hat{\lambda}$. In view of the definitions \cref{eqn:def:postFDR1,eqn:def:Lvals}, we have the following expression for the posterior FDR:
	\begin{equation}\label[equation]{eqn:postFDR} \postFDR_{H}(\vphi) =\frac{\sum_{i=1}^N \ell_{i,H} \vphi_i}{1\vee(\sum_{i=1}^n \vphi_i)}.\end{equation} 
	That is, the posterior FDR of a procedure $\vphi$ is the average of the selected $\ell$--values. Consequently, the procedure $\vphi_{\hat{\lambda},\hat{H}}$ must threshold at one of the ``empirical $\ell$--values'' (i.e.\ at some $\hat{\ell}_i=\ell_{i,\hat{H}}$), as  $\postFDR_{\hat{H}}(\vphi_{\lambda,\hat{H}})$ only changes when $\lambda$ crosses such a threshold. The threshold $\hat{\lambda}$ can therefore equivalently be expressed, as in \cite{SC09}, as $\hat{\lambda}=\hat{\ell}_{(\hat{K}+1)}$, with $\hat{\ell}_{(i)}$ denoting the $i$th order statistic\footnote{We define the order statistics so that repeats are allowed: the order statistics are defined by the fact that $\braces{\ell_i, i\leq N}=\braces{\ell_{(j)}, j\leq N}$ as a multiset ($\forall x\in\RR$, $\#\braces{i : \ell_i=x}=\#\braces{i : \ell_{(i)}=x}$) and $\ell_{(1)}\leq \ell_{(2)}\leq \dots \leq \ell_{(N)}$.} of $\braces{\ell_{i,\hat{H}}: 1\leq i \leq N}$,  where $\hat{K}$ is defined by
	\begin{equation}
		\label[equation]{eqn:def:HatK}
		\frac{1}{\hat{K}} \sum_{i=1}^{\hat{K}} \hat{\ell}_{(i)} \leq t < \frac{1}{\hat{K}+1} \sum_{i=1}^{\hat{K}+1} \hat{\ell}_{(i)}.
	\end{equation}
	(By convention the left inequality automatically holds in the case $\hat{K}=0$, and we define $\hat{\ell}_{(N+1)}:=\infty$ so that the right inequality automatically holds in the case $\hat{K}=N$.) Note that $\hat{K}$ is well defined and unique, by monotonicity of the average of nondecreasing numbers. This monotonicity also makes clear the following dichotomy: 
	\begin{equation}\label[equation]{eqn:HatLambdaCharacterisesPostFDRSign}\postFDR_{\hat{H}}(\vphi_{\lambda,\hat{H}})\leq t \iff \lambda \leq \hat{\lambda}.\end{equation}

	If there are no ties, the procedure $\vphi_{\hat{\lambda},\hat{H}}$ necessarily rejects $\hat{K}$ of the null hypotheses. In the case of ties, it may reject fewer, and to avoid potential conservativity, we therefore consider a slightly adjusted procedure $\hat{\vphi}$. 
	\begin{definition}\label[definition]{def:HatPhi}
	Define $\hat{\vphi}=\hat{\vphi}^{(t)}$ to be a procedure rejecting exactly $\hat{K}$ of the hypotheses with the smallest $\hat{\ell}_i$ values, choosing arbitrarily in case of ties, where $\hat{K}$ is defined by \cref{eqn:def:HatK}. We write $\hat{S}_0$ for the rejection set \[\hat{S}_0=\braces{i\leq N : \hat{\vphi}_i=1},\] and we note that by construction we have $\abs{\hat{S}_0}=\hat{K}$ and \[\braces{i : \hat{\ell}_i(X)<\hat{\lambda}}\subseteq \hat{S}_0\subseteq \braces{i : \hat{\ell}_i(X)\leq \hat{\lambda}}.\]
	\end{definition}
We make the following assumptions on the parameters. The assumptions are not particularly restrictive, and are discussed in detail in \cref{sec:applicability}.
\begin{assumption}
	\label[assumption]{ass:fpi} \begin{enumerate} \item \label{ass:ExistsCNu} There exists a constant $\nu>0$ such that \begin{equation*}\label[equation]{eqn:ExistCNu}\max_{j=0,1} E_{X\sim f_j}(\abs{X}^\nu)<\infty.\end{equation*} 
		\item \label{ass:ExistsX*}
		There exists $x^*\in \RR\cup\braces{\pm \infty}$ such that either \begin{align*}f_1(x)/f_0(x)&\to \infty, \quad  \text{as }x\uparrow x^*, \quad \text{ or} \\
	f_1(x)/f_0(x)&\to \infty, \quad  \text{as }x\downarrow x^*,
	 \end{align*}
		 where we take the conventions that $1/0=\infty$, $0/0=0$. [In the case where $\mu$ is counting measure on $\ZZ$ and $x^*\not\in\braces{\pm\infty}$, the limits are interpreted to mean that $f_1(x^*)>0$ and $f_0(x^*)=0$.] 
	\end{enumerate}
\end{assumption}

\begin{assumption}\label[assumption]{ass:Qpi}
	\begin{enumerate}
				\item  \label{ass:Qfullrank}The matrix $Q$ has full rank (i.e.\ its two rows are distinct), and 
				\begin{equation*}\delta:=\min_{i,j} Q_{i,j}>0.\end{equation*} 
		\item  The Markov chain is stationary: the initial distribution $\pi=(\pi_0,\pi_1)$ is the invariant distribution for $Q$. 
	\end{enumerate}
\end{assumption}
Throughout we will write
\begin{equation}\label[equation]{eqn:def:fpi}
	f_\pi(x)=\pi_0f_0(x)+\pi_1 f_1(x)
\end{equation}
for the marginal distribution of each $X_i$, $i\leq N$, under \cref{ass:Qpi}; note that necessarily $\min(\pi_0,\pi_1)\geq \delta$ under the assumption. We note the following illustrative examples of pairs of densities with respect to the Lebesgue measure $\mu=\dx$ which satisfy both parts of \cref{ass:fpi}. 
\begin{examples}
	\begin{enumerate}[i.]
		\item $f_j(x)=\phi(x-\mu_j)$, where $\phi$ is the density of a standard normal random variable and $\mu_1\neq \mu_2$. 
\item $f_0$ is the density of any normal random variable, and $f_1$ is the density 
of any Cauchy random variable, or any other distribution with polynomial tails.
\item $f_0,f_1$ are compactly supported densities, and the support of $f_1$ is not a subset of the support of $f_0$.
\item $f_0,f_1$ are the densities of Beta random variables, $f_j(x)=c_j x^{\alpha_j-1}(1-x)^{\beta_j-1}\II\braces{x\in [0,1]}$ for a normalising constant $c_j$, and $\alpha_0>\alpha_1$ or $\beta_0>\beta_1$ (or both).
	\end{enumerate}
\end{examples}

\subsection{Theoretical guarantees}\label[section]{sec:TestingGuarantees}
Our main result shows that for suitably chosen $\hat{H}=(\hat{Q},\hat{\pi},\hat{f}_0,\hat{f}_1)$, the procedure $\hat{\vphi}$ achieves an FDR upper bounded by the level $t$ chosen by the user, at least asymptotically. The existence of estimators with suitable consistency properties is shown in the next \namecref{sec:SupNormEstimation} under mild further assumptions. Here $\norm{\cdot}$ denotes the usual Euclidean norm for vectors (and later also the corresponding operator norm for matrices), $\norm{\cdot}_F$ denotes the Frobenius matrix norm $\norm{A}_F^2=\sum_{ij} A_{ij}^2$, and $\norm{\cdot}_\infty$ denotes the $L^\infty$ (supremum) norm on functions taking values in $\RR$.
\begin{theorem} \label[theorem]{thm:BFDRcontrol}
Grant \cref{ass:fpi,ass:Qpi}.  Suppose that for some $u>1+\nu^{-1}$ and some sequence $\eps_N$ such that $\eps_N (\log N)^u\to 0$, the estimators $\hat{Q},\hat{\pi}$ and $\hat{f}_j,~j=0,1$ 
satisfy \begin{equation} \label[equation]{eqn:ConsistencyAssumption}\Pi_H(\max\braces{\norm{\hat{Q}-Q}_F,\norm{\hat{\pi}-\pi},\norm{\hat{f}_0-f_0}_\infty,\norm{\hat{f}_1-f_1}_\infty}> \eps_N)\to 0, \quad \text{as }N\to \infty.\end{equation} 
Then for $\hat{\vphi}$ the multiple testing procedure of \cref{def:HatPhi} we have \[ \FDR_{H} (\hat{\vphi})\to \min(t,\pi_0).\] 
\end{theorem} 
As alluded to, the construction of $\hat{\vphi}$ suggests it should have close to optimal power, and the following result shows that this is indeed true under an extra condition on the distribution of $(f_1/f_0)(X_1)$. The extra condition is only used to prove a property of the limiting $\ell$--values, so that a version of \cref{thm:BFNRoptimality} may also hold in the discrete setting -- see the discussion in \cref{sec:Extensions}. 
As is common in the literature (again see \cref{sec:Extensions}), the precise notion of power is given by the marginal true discovery rate (mTDR), the average proportion of true signals which a testing procedure discovers: 
\begin{equation} \label[equation]{eqn:def:mTDR} \mTDR_H (\vphi) = \frac{E_H\#\braces{i : \theta_i=1,\vphi_i=1}}{E_H\#\braces{i: \theta_i=1}}. \end{equation}
The marginal FDR is defined correspondingly:
\begin{equation}\label[equation]{eqn:def:mFDR} \mFDR_H(\vphi) = \frac{E_H\#\braces{i : \theta_i=0,\vphi_i=1}}{E_H\#\braces{i: \vphi_i=1}}, \end{equation}
with the convention that $0/0=0$. 
These `marginal' quantities are, by concentration results, close to the original quantities $\TDR_H(\vphi)$, $\FDR_H(\vphi)$ for many procedures, including $\hat{\vphi}$ (as is implied by ideas in the proof of the following result).
\begin{theorem}\label[theorem]{thm:BFNRoptimality}
In the setting of \cref{thm:BFDRcontrol}, additionally grant that the distribution function of the random variable $(f_1/f_0)(X_1)$ is continuous and strictly increasing. 
Then the procedure $\hat{\vphi}$ of \cref{thm:BFDRcontrol} satisfies the following as $N\to \infty$:
\begin{align*} \mTDR_H(\hat{\vphi}) &= \sup\braces{ \mTDR_H(\psi) : \mFDR_H(\psi)\leq \mFDR_H(\hat{\vphi})}+o(1) \\
	&= \sup\braces{ \mTDR_H(\psi) : \mFDR_H(\psi)\leq t}+o(1).
\end{align*}
The suprema are over all multiple testing procedures $\psi$ satisfying the bound on their mFDR, including oracle procedures allowed knowledge of the parameters $H$.
\end{theorem}
The essence of the proof of \cref{thm:BFDRcontrol} is to show that $\hat{\ell}_i\approx \ell_i$ for most $i\leq N$ (see \cref{lem:HatlVslErrors}, in \cref{sec:proofs:FDRcontrol}) and that consequently $\postFDR_{H}(\hat{\vphi})$ is close to $\postFDR_{\hat{H}}(\hat{\vphi})$. The latter, thanks to our definition of $\hat{\lambda}$, is close $t$.
	
In proving \cref{thm:BFNRoptimality}, there is no \emph{a priori} control of the power analogous to the bound $\postFDR_{\hat{H}}(\hat{\vphi})\leq t$, hence we cannot simply argue by symmetry. Instead, one shows that $\hat{\lambda}$ concentrates around some $\lambda^*\in(0,1]$: see \cref{lem:LambdaHatConcentrates} in \cref{sec:proofs:FDRcontrol}. Then, again using that $\hat{\ell}_i\approx \ell_i$, it follows that $\mTDR_{H}(\hat{\vphi})\approx \mTDR_{H}(\vphi_{\lambda^*,\hat{H}})\approx \mTDR_{H}(\vphi_{\lambda^*,H})$ and similarly that $\mFDR_{H}(\hat{\vphi})\approx \mFDR_{H}(\vphi_{\lambda^*,H})\approx t$. Known optimality results for the class $(\vphi_{\lambda,H} : \lambda \geq 0)$ mean that one is able to show that $\mTDR_H(\vphi_{\lambda^*,H})$ is the largest of procedures with mFDR at most $\mFDR_H(\vphi_{\lambda^*,H})\approx t$ (see \cref{lem:OptimalityOfOracleClass}), so that the same is approximately true of $\mTDR_{H}(\hat{\vphi})$.

See \cref{sec:proofs:FDRcontrol} for the proofs.

\section{Supremum Norm Estimation of Emission Densities}\label[section]{sec:SupNormEstimation}

Of course, \cref{thm:BFDRcontrol,thm:BFNRoptimality} are only useful if one can estimate $H$ at an appropriate rate in the specified norms, and the results of this section ensure that this is indeed possible in a wide range of nonparametric settings. Estimation is possible not only in the two-state setting, and since estimation results are of independent interest we assume in this section that the data $X$ is drawn from model \cref{eqn:def:model} for $Q$ a $J\times J$ matrix and $\pi$ a distribution on $\braces{1,\dots,J}$, with the state vector $\theta$ taking values in $\braces{1,\dots,J}^N$, for some known $J\geq 2$.

\Cref{ass:Qpi,ass:fpi} are designed with the particular FDR context in mind. In the $J$--state estimation setting we instead use the following conditions, designed to ensure a spectral estimation method works.

\begin{manualassumption}{\ref*{ass:Qpi}'}\label[assumption]{ass:Qpi'}
The matrix $Q$ is full rank, the $J$--state Markov chain $(\theta_n)_{n\in \NN}$ is irreducible and aperiodic, and $\theta_1$ follows the invariant distribution. [This is weaker than \cref{ass:Qpi} in general, but equivalent in the two-state setting.]
\end{manualassumption}

\begin{assumption} \label[assumption]{ass:LinearIndependence} The density functions $f_1,\dots f_J$ are linearly independent. [In the two-state setting it suffices to assume $f_0\neq f_1$, which is implied by \cref{ass:fpi}.]
\end{assumption}
Under these assumptions, in a parametric setting a variant of a typical regularity condition suffices to show that estimation is possible at a parametric rate, so that our theorems offer a new proof of the results of \cite{SC09}: see \cref{sec:ParametricSetting}.
Of greater interest here, though, is that \cref{thm:BFDRcontrol} also allows for a nonparametric setting. As noted already, this is a major improvement for applications -- see for instance \cite{YPRH11},  \cite{WSZ19} and \cite{SW20}.
Estimating the Markov parameters $Q$ and $\pi$ consistently up to a permutation at a polynomial rate has already been proved possible (see \cite[Appendix C]{dCGlC17}), and
we therefore focus on estimation, in the supremum norm, of the emissions densities themselves. Note first of all that in a discrete setting estimation is possible at a parametric rate. 
\begin{theorem}\label[theorem]{thm:DiscreteEstimation}
Assume that the dominating measure $\mu$ is the counting measure on $\ZZ$.
Let $M_N$ be a sequence tending to infinity, arbitrarily slowly. Under \cref{ass:Qpi',ass:LinearIndependence}, there exist estimators $\hat{f}_1,\dots, \hat{f}_J$ and a permutation $\tau$ such that
\[ \Pi_H(\norm{\hat{f}_j-f_{\tau(j)}}_\infty\geq M_N N^{-1/2})\to 0.\]
\end{theorem}
The proof is a simplification of that of \cref{thm:LinftyEstimation} (to follow) and so is sketched only: see \cref{sec:ProofOfDiscreteEstimation}.

\smallskip

For the remainder of this \namecref{sec:SupNormEstimation} we assume that the functions $f_1,\dots, f_J$ are densities with respect to the Lebesgue measure on $\RR$, $\mu=\dx$. We demonstrate that consistent estimation of these densities in the supremum norm is indeed possible at a near-minimax rate in the nonparametric setting, under the following typical smoothness condition. 

\begin{assumption} \label[assumption]{ass:smoothness} $f_1,\dots f_J$ belong to $C^s(\RR)$ for some $s>0$, where for $C^0(\RR)$ denoting all bounded continuous functions from $\RR$ to itself (equipped with the usual supremum norm $\norm{\cdot}_\infty$) and writing $j=\floor{s}$ for the integer part of $s$, $C^s(\RR)$ denotes the usual space of (locally) H\"older-continuous functions 
	\begin{alignat*}{3} &C^s(\RR)=\braces{f: f^{(j)}\in C^{s-j}(\RR)}, &&  \quad s\geq 1 \\
		&C^s(\RR)= \braces{f \in C^0(\RR) : \sup_{0<\abs{x-y}\leq 1}\brackets[\Big]{\frac{\abs{f(x)-f(y)}}{\abs{x-y}^s}}<\infty} && \quad s\in(0,1),
		\end{alignat*} equipped with the usual norm
	\begin{alignat*}{3} &\norm{f}_{C^s}=\norm{f^{(\floor{s})}}_{C^{s-\floor{s}}}+\sum_{0\leq i<\floor{s}}\norm{f^{(i)}}_{\infty}, && \quad s\geq 1 \\
	&\norm{f}_{C^s}= \norm{f}_\infty + \sup_{0<\abs{x-y}\leq 1} \frac{\abs{f(y)-f(x)}}{\abs{y-x}^s}, &&\quad 0<s<1. 
	\end{alignat*}  
	\end{assumption}
The results also extend in the usual way to Besov spaces, e.g.\ using results from \cite[Chapter 4]{GN16}.

	\begin{theorem}\label[theorem]{thm:LinftyEstimation} 
	Grant \cref{ass:Qpi',ass:LinearIndependence,ass:smoothness}. Suppose $L_0\to \infty$ as $N\to \infty$, and $L_0^{\max(5,(J+3)/2)}r_N\to 0$, where $r_N=(N/\log N)^{-s/(1+2s)}$. Then there exist estimators $\hat{f}_j,$ $1\leq j\leq J$ (continuous so that the supremum below is measurable) and a permutation $\tau$ such that, for some $C>0$, \begin{equation}\label[equation]{eqn:InProbabilitySupNormResult}\Pi_H\brackets[\Big]{\norm{\hat{f}_{j}-f_{\tau(j)}}_{\infty}\geq C L_0^5 r_N} \to 0. 
	\end{equation}
	Convergence in expectation also holds: for some $C'>0$,
	\begin{equation}\label[equation]{eqn:InExpectationSupNormResult} E_H \norm{\hat{f}_j-f_{\tau(j)}}_\infty \leq C' L_0^5r_n.\end{equation}
\end{theorem}
The proof is given in \cref{sec:proofs:SupNormEstimation}. The parameter $L_0$ has the interpretation of the dimension of a matrix used in the contruction of the estimators (see \cref{alg:Linfty}) and it can be chosen to diverge arbitrarily slowly, so that the upper bound is arbitrarily close to the following lower bound. Such a lower bound is familiar from the i.i.d.\ setting, but does not automatically apply in the current setting. Indeed, the mixture components in a nonparametric mixture model are not identifiable, so that our assumptions necessarily exclude the i.i.d.\ subcase of a HMM. The content of the following \namecref{prop:MinimaxLowerBound} is that these assumptions do not, however, make estimation easier than having i.i.d.\ samples from each of the emission densities. We refer to \cref{sec:ProofOfLowerBound} for a formal statement and proof.
\begin{proposition}[informal statement]\label[proposition]{prop:MinimaxLowerBound}
	The rate $r_N=(N/\log N)^{-s/(1+2s)}$ is a lower bound for the minimax supremum-norm estimation rate for the emission densities in a two--state nonparametric HMM. 
\end{proposition}
The algorithm solving \cref{thm:LinftyEstimation} uses a `spectral' method similar to those of \cite{AHK12,Lehericy18,dCGlC17}. However, \cite{dCGlC17} and \cite{Lehericy18} expand in terms of orthonormal basis functions, and use particular properties of $L^2$--projections which do not straightforwardly adapt to the $L^\infty$ setting. Here, we instead consider spectral \emph{kernel density estimation}: see \cref{alg:Linfty} for a description of the estimating procedure. This approach allows us to directly estimate the values of the density functions at each point $x$ and bypass the need for projection properties. 

Finally, note that \cref{thm:DiscreteEstimation,thm:LinftyEstimation} 
only show that one may estimate the parameters consistently \emph{up to a permutation}. While this is generally sufficient for estimation purpose, since the labelling of the states is usually of no relevance, any multiple testing procedure targeting FDR control necessarily treats the null and the alternative differently, so it is essential that we can identify which of our estimators corresponds to the null state. We will therefore also require the following condition. 
\begin{condition} \label[condition]{ass:TauEstimable}
There exist estimators $\hat{f}_1,\dots,\hat{f}_J$ in \cref{thm:LinftyEstimation} (or \cref{thm:DiscreteEstimation}) for which the permutation $\tau$ is the identity.
\end{condition}
It suffices that there exist $\braces{\hat{f}_1,\dots, \hat{f}_J,\tau}$ as in \cref{thm:LinftyEstimation} for which the permutation $\tau$ can be estimated consistently by some $\hat{\tau}$, since we can define $\check{f}_j=\hat{f}_{\hat{\tau}(j)}$. We give two illustrative assumptions, each plausible in the original two-state FDR setting, under which \cref{ass:TauEstimable} holds. A version of the following \namecref{prop:ExistenceOfEstimators} also holds under such assumptions in the discrete setting, using \cref{thm:DiscreteEstimation} in place of \cref{thm:LinftyEstimation} in the proof, which is given at the end of \cref{sec:proofs:SupNormEstimation}.

\begin{proposition}\label[proposition]{prop:ExistenceOfEstimators}
In the setting of \cref{thm:BFDRcontrol} grant \cref{ass:Qpi,ass:fpi,ass:smoothness}. Then \cref{ass:TauEstimable} is verified, and there exist estimators $\hat{Q},\hat{\pi},\hat{f}_0,\hat{f}_1$ satisfying \cref{eqn:ConsistencyAssumption} for any rate $\eps_N$ slower than $r_N=(N/\log N)^{-s/(1+2s)}$, under either of the following assumptions:
\begin{enumerate}
	\item $\pi_0>\pi_1$.
	\item For some \emph{known} $x^*\in \RR\cup\braces{+\infty}$, $f_1(x)/f_0(x)\to \infty$ as $x\uparrow x^*$. 
\end{enumerate}	
\end{proposition}

\section{Discussion}\label[section]{sec:DiscussionOfMainResults}

\subsection{Applicability of the Results}\label[section]{sec:applicability}
\paragraphi{Generality of the assumptions}
\Cref{ass:Qpi,ass:TauEstimable,ass:fpi,ass:LinearIndependence,ass:smoothness} are not restrictive, so that \cref{thm:BFDRcontrol,thm:DiscreteEstimation,thm:LinftyEstimation} hold in typical nonparametric settings (we discuss the extra assumption of \cref{thm:BFNRoptimality} in \cref{sec:Extensions}).

\Cref{ass:fpi}\,\cref{ass:ExistsX*} is a signal strength assumption, without which the proofs remain valid only for large enough values of $t$. 
It is known that weak signals are a case requiring special attention for multiple testing, discussed for example in a different setting in \cite{HellerRosset20}.  
 
The full rank assumption on $Q$ in \cref{ass:Qpi'} is necessary even for indentifiability up to a permutation in the two-state case (with nonparametric emission densities). For $J>2$ states it is not known whether identifiability holds without this assumption, and full rank is assumed in all papers concerning nonparametric inference of HMM parameters. Irreducibility is essential to ensure all hidden states genuinely influence the data. Aperiodicity is assumed to allow typical Markov chain convergence and concentration results to apply, but in principle it should be possible to avoid this assumption at the expense of requiring specially tailored proofs, since the proofs use empirical averages as a building block.

Consistent estimation of the HMM parameters is possible upon replacing \cref{ass:LinearIndependence} by the weaker assumption that the emission densities are all distinct, see \cite{AHL16} and \cite{L18}. Proving rates under this weaker assumption is much harder and no results exist yet.

\paragraphi{Implementing the method}
Our proposed method for estimating the emission densities can be implemented through Algorithm \ref{alg:Linfty}. Then, given estimators of the parameters, efficient computations of $\ell$--values is easily done using the forward--backward algorithm for HMMs. Indeed the empirical Bayes multiple testing procedure is implemented in 
	\cite{SC09},  \cite{WSZ19} and \cite{SW20}.
[These works use mixture models with unknown number of components to estimate the emission densities, either via fully Bayesian methods or via penalized maximum likelihood (using the EM algorithm).]

\subsection{The Parametric Setting}\label[section]{sec:ParametricSetting}
\cite{BRR98} prove a central limit theorem for the maximum likelihood estimator of the model parameter (which we denote, say, by $h$) under standard regularity conditions, so that it may be estimated at a parametric rate up to label switching. To these, adding the condition that the parametrisation map $h\mapsto (f_{1,h},\dots f_{J,h})$ is Lipschitz continuous with respect to the Euclidean norm and the supremum norm (at least on a neighbourhood of the true parameter), we arrive at the following. 
\begin{proposition}\label[proposition]{prop:ParametricEstimation}
In a parametric model satisfying mild regularity conditions, \cref{ass:Qpi',ass:LinearIndependence} are enough to ensure that there exist estimators $\hat{Q},\hat{\pi},\hat{f}_1,\dots,\hat{f}_J$ such that for some permutation $\tau$ and any $M_N\to\infty$,
	\[\max\brackets[\big]{ \norm{\hat{Q}-Q}_F,\norm{\hat{\pi}-\pi},\norm{\hat{f}_1-f_{\tau(1)}}_\infty,\dots,\norm{\hat{f}_J-f_{\tau(J)}}_\infty}<M_N N^{-1/2}, \]
	with probability tending to 1.	
\end{proposition}
We note that many common parametric families, including Gaussian models, exponential models and Poisson models, satisfy a suitable regularity condition (this can be seen by using standard formulae for exponential families to calculate the derivative of the parametrisation map and bounding).

Under an assumption akin to those of \cref{prop:ExistenceOfEstimators} to ensure that a version of \cref{ass:TauEstimable} holds, we see that \cref{thm:BFDRcontrol,thm:BFNRoptimality} apply in a parametric setting. Except perhaps for the regularity condition, our assumptions are weaker than those of \cite{SC09} (after adapting \cref{thm:BFNRoptimality} slightly -- see \cref{sec:Extensions}), so that we slightly generalise their main results even in the parametric setting.

\subsection{Uniformity in the Parameters}\label[section]{sec:Uniformity}
	The constants of \cref{thm:LinftyEstimation} depend only on quantitative measures (as listed below) of the degree to which  \cref{ass:Qpi',ass:LinearIndependence,ass:smoothness} hold, so that a uniform version of \cref{eqn:InExpectationSupNormResult},
		\[\sup_{H\in\Hh} E_H \norm{\hat{f}_j-f_{\tau(j)}}_\infty \leq C' L_0^5r_n,\] holds if the following bounds are satisfied on the set $\Hh$ (and similarly for \cref{eqn:InProbabilitySupNormResult}). The estimators $\hat{f}_1,\dots \hat{f}_J$ do not depend on knowledge of the bound $M<\infty $, so the result is adaptive in these quantities (though recall that the smoothness $s$ is assumed known -- see also the discussion of adaptation in \cref{sec:Extensions}).
	\begin{itemize}
		\item $\sup_{H\in\Hh}(\kappa(Q))\leq M$, where $\kappa(Q)=\norm{Q}\norm{Q^{-1}}$, the condition number, measures how far $Q$ is from having less than full rank.
		\item $\inf_{H\in\Hh} \gamma_{\textnormal{ps}} \geq M^{-1}$ where $\gamma_{\textnormal{ps}}$ denotes the pseudo spectral gap of the matrix $Q$ as defined in \cite{Paulin2015}. This bound quantitatively measures how far the chain $\theta$ is from being reducible or periodic, and is only used to control the mixing time of the chain $\theta$. It can therefore be replaced by any assumption ensuring a uniform bound on the mixing time; in particular, in the two-state case of \cref{sec:EmpiricalBayesProcedure}, the chain $\theta$ is necessarily reversible and it suffices to assume a uniform lower bound on the absolute spectral gap $\gamma^*$, defined by
		\[\gamma^*=\begin{cases} 	1-\sup\braces{\abs{\lambda} : \lambda \text{ an eigenvalue of $Q$, }\lambda\neq 1} & \text{the eigenvalue 1 of $Q$ has multiplicity 1}, \\ 0 & \text{otherwise.}		
		\end{cases}
			\] 
		\item $\inf_{H\in \Hh} \min_j(\pi_j)>M^{-1}$. This too measures how far the chain is from being reducible.
		\item $\sup_{H\in\Hh}\max_j\norm{f_j}_{C^s}\leq M$.
		\item $\sup_{H\in\Hh}\max(\underline{L},1/C)\leq M$, where $(C,\underline{L})$ are the constants, depending on $H$, from \cref{lem:ExistenceOfSuitableHl} in \cref{sec:lemmasforLinftyEstimation}. Denoting by $\sigma_J(A)$ the $J$th largest singular value of a matrix $A$, these constants control $\sigma_J(O^{L_0})$ where $O^{L_0}=(E[h_l(X_1) \mid \theta_1=j]_{l\leq L_0,j\leq J})$ for some suitably chosen functions $h_l,$ $l\leq L_0$. The \namecref{lem:ExistenceOfSuitableHl} shows that $h_1,\dots,h_{L_0}$ can be chosen in a universal way such that $\max(\underline{L},1/C)<\infty$ whenever $f_1,\dots f_J$ are linearly independent, so these constants quantitively measure the linear independence of these functions. In the case $J=2$, a sufficient (but not necessary) condition for such a uniform bound to hold is that $P_{X\sim f_0}(X\in A)\neq P_{X\sim f_1}(X\in A)$ for some \emph{known} set $A$: one constructs the estimators $\hat{f}_0$, $\hat{f}_1$ using, in \cref{alg:Linfty}, $L_0=2$, $h_1=1$, $h_2=\II_A$.
		\item $\inf_{H\in\Hh}c\geq M^{-1}$ where $c=c(H)$ is the constant of \cref{lem:ExistenceOfHatAHatU} in \cref{sec:lemmasforLinftyEstimation}. The \namecref{lem:ExistenceOfHatAHatU} shows that this constant is positive whenever $f_1,\dots,f_J$ are distinct and so it provides a quantitative measure of the degree of distinctness of these functions. In view of the proof, a sufficient (but not necessary) condition for such a uniform bound to hold is that $f_1,\dots f_J$ can uniformly be separated at a point, i.e.\ that the set $\Hh$ is such that \[ \inf_{H\in \Hh} \sup_{x\in\RR} \min_{j\neq j'} \abs{f_j(x)-f_{j'}(x)}>M^{-1}.\]
	\end{itemize}
	In what follows, we use for example $C=C(\Hh)$ to denote any constant which depends only on the above bounds (i.e.\ on $M<\infty$). We note that the set $\Hh$ over which the upper bound is uniform (under the sufficient conditions of the last two items, with $A=[-1,1]$) includes the set over which the lower bound \cref{prop:MinimaxLowerBound} is proved in \cref{sec:ProofOfLowerBound}, so that the estimation result \cref{thm:LinftyEstimation} can genuinely be viewed as a minimax result. 
		
	\smallskip
	
	The FDR result \cref{thm:BFDRcontrol} is uniform over a large subset $\Ii\subset \Hh$. In particular, in addition to the above constraints, one needs to add the following conditions.  
	\begin{itemize}
		\item $\sup_{H\in\Ii}(\max_j E_{X\sim f_j}\abs{X}^{\nu_0})<\infty$ for some $\nu_0=\nu_0(\Ii)>0$.
		\item $\inf_{H\in\Ii}\Pi_H((f_1/f_0)(X_1)>u)>0$ for each $u>0$.
		\item \Cref{ass:TauEstimable} holds in a uniform way for $H\in \Ii$.
		\item $\inf_{H\in \Ii}\min_{i,j}Q_{ij}>0$. [This is in fact implied already by the bounds on the $\pi_j$ and on the pseudo-spectral gap, since for \cref{thm:BFDRcontrol} we are in the two-state setting.]
	\end{itemize}
 We write $C=C(\Ii)$ to denote any constant which depends only on $\Hh$ and these quantities.

\subsection{Extensions of the Theorems}\label[section]{sec:Extensions}
	\paragraphi{Weakening the assumption of \cref{thm:BFNRoptimality}}
	\Cref{thm:BFNRoptimality} remains true if we replace the assumption on $(f_1/f_0)(X_1)$ with the following; see \cref{lem:liinftyFullSupport} for 
	a proof that this new condition holds under the assumptions of \cref{thm:BFNRoptimality}.
	\begin{condition}\label[condition]{condition:liinfty}
		Viewing the sample $(X_n :1\leq n\leq N)$ as coming from a bi-infinite HMM $(X_n : n\in \ZZ)$, grant that the distribution function of 
		\begin{equation}\label[equation]{eqn:def:ellinfty} \ell_i^\infty(X) := \Pi_H (\theta_i=0 \mid (X_n)_{n\in \ZZ})\end{equation} is continuous and strictly increasing on $[0,1]$.	
	\end{condition}
	This condition is weaker than the `monotone ratio condition' assumed in \cite{SC09}, since the latter implicitly assumes that the distribution function of $\ell_i^\infty$ has a strictly positive derivative. In the discrete context (that is, when the $X_{i}$'s take discrete values), understanding when the distribution of the variables $\ell_i^\infty$ have a density with respect to Lebesgue measure is known to be hard, since it is mostly still an open problem for the closely related stationary filter $ \Phi_i^\infty(X) := \Pi_H (\theta_i=0 \mid (X_n)_{n\leq i})$, see \cite{B57}, \cite{BK15} and references therein.

	Of particular interest, though, is the fact that this new condition is only about the continuity of the distribution function, not about its absolute continuity. Continuity is a weaker property that could be easier to understand and could hold in much more generality, so that	Condition \ref{condition:liinfty} opens up the possiblity that a version of \cref{thm:BFNRoptimality} may hold even in certain discrete settings. Indeed, simulations in \cite{SW20} are suggestive that the conclusions of the \namecref{thm:BFNRoptimality} hold. They compare various multiple testing procedures and provide empirical evidence that the TDR of the empirical Bayes multiple testing method using nonparametric modeling of HMMs roughly matches that of an oracle thresholding procedure and is the best among the multiple testing procedures they compare.

	\paragraphi{Use of marginal FDR and TDR in \cref{thm:BFNRoptimality}}
	The proof of \cref{thm:BFNRoptimality} in fact shows, after some minor adjustments, that
	\[ \TDR_H(\hat{\vphi})\geq \TDR_H(\vphi_{\lambda_{\max},H})-o(1),\]
	where $\lambda_{\max}=\lambda_{\max}(t,H)$ is chosen maximal such that $\FDR_H(\vphi_{\lambda_{\max},H})\leq t$, so that $\hat{\vphi}$ is (asymptotically) optimal for the TDR when restricting to the class of procedures whose TDR and FDR asymptotically coincide with their marginal equivalents. \cite{HellerRosset20} show in a non-Markovian setting that the procedure maximising the TDR among all procedures with controlled FDR is not in this class, but their results leave open the possiblity that \cref{thm:BFNRoptimality} remains true with the full FDR and TDR. Indeed, a main conclusion of their work is that the class $(\vphi_{\lambda,H} : \lambda \geq 0)$ (or rather, the equivalent of this class for their setting) is optimal for the problem of maximing TDR with controlled FDR provided one allows data-driven thresholds -- such as $\hat{\lambda}$ -- whereas the current proof of \cref{thm:BFNRoptimality} uses that for mTDR optimality with mFDR control it suffices to consider the class for non-random thresholds. Furthermore, the difference between the FDR and TDR of the optimal procedure and their marginal versions in the setting of \cite{HellerRosset20} manifests itself for weak signals, so that our signal strength assumption may suffice to rule out such differences.

\paragraphi{Adaptation}
The estimator we construct for \cref{thm:LinftyEstimation} uses knowledge of the smoothness $s$. One can adjust the arguments of \cite{Lehericy18} to show that a careful application of Lepskii's method allows adaptation up to a maximum smoothness $s_{\max}<\infty$ -- and indeed state-by-state adaptation, wherein each state is estimated at a rate adapting to its smoothness parameter $s_j$, rather than requiring $s_j=s$ for all $j$. As usual, the rough idea is to constructs estimators $\hat{f}_j^L,~j\leq J$ for each $L\leq L_{\max}$ and use $\norm{\hat{f}_j^L-\hat{f}_j^{L_{\max}}}_\infty$ as a proxy for the bias, so that one can make a suitable bias-variance tradeoff. In the HMM setting, as noted in \cite{Lehericy18}, one must also use $\hat{f}_j^{L_{\max}}$ to ``align'' the estimators $\hat{f}_j^L$ up to a single permutation $\tau$ rather than needing a different permutation $\tau_L$ for each level $L$; one can show using the triangle inequality that this alignment is successful for all large enough $L\leq L_{\max}$ with probability tending to 1. 

\paragraphi{General measure spaces}
The proofs of \cref{thm:BFDRcontrol,thm:BFNRoptimality} essentially only use the assumption that $\mu$ is Lebesgue measure on $\RR$ or counting measure on $\ZZ$ in showing \cref{lem:ControlOf1/f}, so that versions of these theorems continue to hold on general (metric) measure spaces after adjusting \cref{ass:fpi} appropriately. 
\Cref{thm:DiscreteEstimation} readily generalises to $\mu$ being any discrete measure of known support. The proof of \cref{thm:LinftyEstimation} uses kernel density estimation techniques, and in principle it should be possible to prove a version of this result in any setting where kernel-type estimators with suitable properties exist -- for example, using results from \cite{CGKPP20},  on manifolds.

\section{Proofs}\label[section]{sec:Proofs}
\subsection{Proofs: FDR Control and TDR Optimality}\label[section]{sec:proofs:FDRcontrol}

The following lemma isolates part of the proof of \cref{thm:BFDRcontrol,thm:BFNRoptimality}, showing that $\hat{\ell}_i(X)$ converges to $\ell_i(X)$ at a rate slightly slower than the convergence rate $\eps_N$ of the estimators $\hat{H}$.
\begin{lemma}\label[lemma]{lem:HatlVslErrors}
In the setting of \cref{thm:BFDRcontrol} define $\eps_N'=\eps_N (\log N)^u$, and recall that by definition $u>1+\nu^{-1}$ and by assumption $\eps_N'\to 0$, where $\nu$ is the parameter of \cref{ass:fpi}. Then
	\begin{equation}\label[equation]{eqn:hatlvslerror} \max_{i\leq N} \Pi_H(\abs{\hat{\ell}_i(X)-\ell_i(X)} >\eps_N') \to 0, \quad \text{as~} N\to\infty. \end{equation}
	Consequently, there exists $\delta_N\to 0$ such that
	\[\Pi_H(\#\braces{i\leq N : \abs{\hat{\ell}_i(X)-\ell_i(X)}>\eps_N'}>N\delta_N)\to 0.\]
\end{lemma}
\begin{proof} 
We begin by showing that  $\Pi_H(\abs{\hat{\ell}_i(X)-\ell_i(X)}>M\eps_N')\to 0$ for each fixed $i$, for some constant $M=M(\Ii)$. [Recall that a constant $M(\Ii)$ depends only on certain bounds for the parameter $H=(Q,\pi,f_0,f_1)$ as described in \cref{sec:Uniformity}.]

Let $(E_N)_N$ be a sequence of events with probability tending 1 on which 
\[ \max\brackets[\Big]{\norm{\hat{Q}-Q}_F, \norm{\hat{\pi}-\pi},\max_{j\in \braces{0,1}} \norm{\hat{f}_j-f_j}_{\infty}} \leq \eps_N, \]
and define \begin{equation*}\begin{split} 
&\delta=\min_{i,j} Q_{i,j}, \quad \hat{\delta}=\min_{i,j}\hat{Q}_{i,j}, \\ &\rho=(1-2\delta)/(1-\delta), \quad \hat{\rho}=(1-2\hat{\delta})/(1-\hat{\delta}).\end{split}\end{equation*} Then Proposition 2.2 of \cite{dCGlC17} yields that for some  $C$ depending only on a lower bound for $\delta$, 	\begin{equation}\label[equation]{eqn:dCGlC17Prop2.2} \begin{split} \abs{\hat{\ell}_i(X)-\ell_i(X)}\leq C& \braces[\Big]{\rho^{i-1} \norm{\hat{\pi}-\pi} + \sqbrackets[\big]{(1-\rho)^{-1} + (1-\hat{\rho})^{-1}} \norm{\hat{Q}-Q}_F+ \\ & \quad \sum_{n=1}^N ((\hat{\rho}\vee \rho)^{\abs{n-i}}/f_\pi(X_n)) \max_{j=0,1} \abs{\hat{f}_j(X_n)-f_j(X_n)}}. \end{split} \end{equation}  
	(The proposition there is stated with $c_*(x):=\min_{j=0,1} \sum_{k} Q_{jk} f_k(x)$ in place of $f_\pi(x)$, but we note $c_*(x)$ so defined is lower bounded by $\delta f_\pi(x)$. Also note that \cite{dCGlC17} assume that $f_0,f_1$ are densities with respect to \emph{Lebesgue} measure, but this assumption is not used in the proof of Proposition 2.2 therein.)	
	Recalling we assumed that $\delta$ was (strictly) positive, we see that on $E_N$, for $N$ large enough we have $\hat{\delta}>\tilde{\delta}=\delta/2$, $\hat{\rho} < \tilde{\rho}=(1+\rho)/2$, and we replace $\rho,\hat{\rho}$ and $\delta,\hat{\delta}$ in \cref{eqn:dCGlC17Prop2.2} by $\tilde{\rho}<1$ and $\tilde{\delta}>0$. On the event $E_N$, choosing the constant $M=M(\tilde{\delta},\tilde{\rho},C)=M(\Ii)$ large enough we see by a union bound that 
	\[ \Pi_H(\abs{\hat{\ell}_i(X)-\ell_i(X)} > M\eps_N') \leq  \Pi_H(E_N^c)+\Pi_H\brackets[\bigg]{\eps_N \sum_{n=1}^N \frac{\tilde{\rho}^{\abs{n-i}}}{f_\pi(X_n)} >\eps_N'}. \]
	For $\kappa>0$ to be chosen, define $S_{\kappa,i}=\braces{n\leq N : \abs{n-i}\leq \kappa \log N}.$ We can split the terms in $S_{\kappa,i}$ from those in $S_{\kappa,i}^c$ to see, for $C'=2\sum_{n=0}^\infty \tilde{\rho}^n<\infty$, that
	\[ \begin{split} \sum_{n\leq N}  \frac{\tilde{\rho}^{\abs{n-i}}}{f_\pi(X_n)} \leq  C' \sqbrackets[\Big]{\max_{n\in S_{\kappa,i}}\brackets[\Big]{\frac{1}{f_\pi(X_n)}} + \tilde{\rho}^{\kappa \log N} \max_{n\leq N} \brackets[\Big]{\frac{1}{f_\pi(X_n)}}}, \end{split} \]
	so that again appealing to a union bound, it suffices to show
	\begin{align}\label[equation]{eqn:maxLogNterms} \Pi_H\brackets[\bigg]{\max_{n\in S_{\kappa,i}}\brackets[\Big]{\frac{1}{f_\pi(X_n)}}>\frac{1}{2C'}(\eps_N'/\eps_N)}&\to 0,\text{ and } \\ \label[equation]{eqn:maxAllTerms}  \Pi_H\brackets[\bigg]{ \tilde{\rho}^{\kappa \log N} \max_{n\leq N} \brackets[\Big]{\frac{1}{f_\pi(X_n)}}> \frac{1}{2C'}(\eps_N'/\eps_N)}&\to 0.\end{align} \Cref{lem:ControlOf1/f} (in \cref{sec:LemmasForBFDRcontrol}) tells us that for any $a>1+\nu^{-1}$, with $\nu$ the constant of \cref{ass:fpi}, we have $\Pi_H(\max_{i\leq R} 1/f_\pi(X_i)>R^a)\to 0$ as $R\to \infty$. By stationarity of the process $X$, taking $R=\abs{S_{\kappa,i}}\leq (2\kappa \log N +1)$ we deduce that \[\Pi_H\brackets[\Big]{\max_{n\in S_{\kappa,i}} \frac{1}{f_\pi(X_n)} > (2\kappa \log N +1)^a}\to 0.\] Recalling that $\eps_N'/\eps_N>(\log N)^u$, we see that	\cref{eqn:maxLogNterms} holds for all $\kappa$  if $u>a$.
	Next we apply \cref{lem:ControlOf1/f} with $R=N$ to see \[\Pi_H\brackets[\bigg]{\max_{n\leq N} \brackets[\Big]{\frac{1}{f_\pi(X_n)}}> N^a}\to 0.\] Noting that $\tilde{\rho}^{-\kappa \log N}=N^{\kappa \log (1/\tilde{\rho})}$ and choosing $\kappa>a (\log 1/\tilde{\rho})^{-1}$ yields \cref{eqn:maxAllTerms}. 
 This concludes the proof that for some constant $M$ and each $i\leq N$, $\Pi_H(\abs{\hat{\ell}_i(X)-\ell_i(X)}>M\eps_N')\to 0$.
	
To see that $\max_i \Pi_H(\abs{\hat{\ell}_i(X)-\ell_i(X)}>\eps_N')\to 0$, we note that by initially considering $\eps_n'$ defined for some $u'<u$ we can remove the constant $M$. Thanks to stationarity of the HMM $X$, we further note that \[\max_{i\leq N} \Pi_H\brackets[\Big]{\max_{n\in S_{\kappa,i}} \frac{1}{f_\pi(X_n)}>(2\kappa \log N+1)^a}\to0;\] then, since the other terms in the calculations above do not depend on $i$, we deduce \cref{eqn:hatlvslerror}.
	
	Finally, defining 
	\[\delta_N=\brackets[\Big]{\max_{i\leq N} \Pi_H(\abs{\hat{\ell}_i(X)-\ell_i(X)}>\eps_N')}^{1/2},\] we appeal to Markov's inequality to see that 	\begin{align*} \Pi_H(\#\braces{i\leq N : \abs{\hat{\ell}_i(X)-\ell_i(X)}>\eps_N'}>N\delta_N)& \leq \frac{1}{N\delta_N}\sum_{i=1}^N \Pi_H\brackets[\big]{\abs{\hat{\ell}_i(X)-\ell_i(X)}>\eps_N'} \\ &\leq \delta_N^{-1}\max_{i\leq N}\Pi_H(\abs{\hat{\ell}_i(X)-\ell_i(X)}>\eps_N')= \delta_N,\end{align*} which tends to zero, concluding the proof.
\end{proof} 

\begin{proof}[Proof of \cref{thm:BFDRcontrol}]  Write $\hat{t}=\postFDR_{\hat{H}}{\hat{\vphi}}$ and recall we write $\hat{S}_0$ for the rejection set of $\hat{\vphi}$. We have, for any sequences of positive numbers $\eps_N'$ and of events $F_N$,
	\begin{align*} 
	\abs{\FDR_H(\hat{\vphi})-E_H\hat{t}} 
	& = \abs[\big]{E_{X\sim\Pi_H}[\postFDR_H(\hat{\vphi})-\postFDR_{\hat{H}}(\hat{\vphi})]} \\ 
	&\leq  E_H\sqbrackets[\bigg]{\frac{\sum_{i=1}^N \abs{\ell_i(X)-\hat{\ell}_i(X)}\II\braces{i\in\hat{S}_0}}{1\vee \abs{\hat{S}_0}}} \\
	&\leq \eps_N' + \Pi_H(F_N^c)+E_H\sqbrackets[\Big]{\II_{F_N}\frac{\sum_{i=1}^N \II\braces{\abs{\ell_i(X)-\hat{\ell}_i(X)}>\eps_N'}}{1\vee \abs{\hat{S}_0}}},
	\end{align*} where we have used that $\abs{\ell_i(X)-\hat{\ell}_i(X)}\leq 1$ for all $i$.
\Cref{lem:hatT} in \cref{sec:LemmasForBFDRcontrol} shows that $E_H[\hat{t}]\to \min(t,\pi_0)$, so that it is enough to show the right side tends to zero for suitable $\eps_N'$ and $F_N$.
	
\Cref{lem:O(N)rejections} tells us that $\Pi(\abs{\hat{S}_0}>aN)\to 1$ for some $a>0$. Combining with \cref{lem:HatlVslErrors} by a union bound, we deduce that for suitably chosen $\eps_N'\to 0,$ $\delta_N\to 0$ and $a>0$, we have $\Pi_H(F_N^c)\to 0$ if we define
	\[
	F_N = \braces[\big]{\#\braces[]{i\leq N : \abs{\hat{\ell}_i(X)-\ell_i(X)}>\eps_N'}\leq N\delta_N}\cap \braces[\big]{\abs{\hat{S}_0}>aN}. \]
	Then \[E_H\sqbrackets[\Big]{\II_{F_N}\frac{\sum_{i=1}^N \II\braces{\abs{\ell_i(X)-\hat{\ell}_i(X)}>\eps_N'}}{1\vee \abs{\hat{S}_0}}}
		 \leq \frac{N\delta_N}{aN}\to 0,  \] yielding the result.
		\qedhere
\end{proof}

The following \namecref{lem:LambdaHatConcentrates}, mentioned already in the sketch proof in \cref{sec:TestingGuarantees}, will help us in proving \cref{thm:BFNRoptimality}.
\begin{lemma}\label[lemma]{lem:LambdaHatConcentrates}
Under the assumptions of \cref{thm:BFNRoptimality}, define $\lambda^*\in(t,1]$ implicitly by \[ E[\ell_i^\infty (X) \mid \ell_i^\infty(X)<\lambda^*]= \min(t,\pi_0),\] where $\ell_i^\infty$ is as in \cref{eqn:def:ellinfty} (by stationarity the conditional expectation does not depend on $i$).

	Such a $\lambda^*$ exists; it satisfies, for $\eps>0$, \begin{alignat*}{2}E[\ell_i^\infty (X) \mid \ell_i^\infty(X)<\lambda^*-\eps]&< \min(t,\pi_0), \quad &&
	\\
	E[\ell_i^\infty (X) \mid \ell_i^\infty(X)<\lambda^*+\eps]&> t \qquad  &&\text{if $t<\pi_0$};
\end{alignat*} and we have \begin{equation}\label[equation]{eqn:hatLambdaConverges}\hat{\lambda}\to \lambda^* \quad \text{in probability as }N\to\infty.\end{equation}
\end{lemma}
\begin{proof}
\Cref{lem:liinftyFullSupport} (in \cref{sec:LemmasForBTDRoptimality}) tells us that under the assumptions of \cref{thm:BFNRoptimality}, the distribution function of $\ell_i^\infty$ is continuous and strictly increasing. \Cref{lem:ConditionalExpectationIncreasing} then tells us that the same is true of the map $\lambda\mapsto E[\ell_i^\infty \mid \ell_i^\infty<\lambda]$, and that $E[\ell_i^\infty \mid \ell_i^\infty<t]<t$. Noting also that $E[\ell_i^\infty \mid \ell_i^\infty <1]=E[\ell_i^\infty]=\pi_0$ (since $\ell_i^\infty<1$ with probability 1), we deduce the existence of a unique solution $\lambda^*\in (t,1]$ by the intermediate value theorem. Strict monotonicity of the conditional expectation implies the claimed inequalities when conditioning on $\ell_i^\infty<\lambda^*-\eps$ and on $\ell_i^\infty<\lambda^*+\eps$.
	
For the convergence in probability, we show for $\eps>0$ arbitrary that with probability tending to 1 we have $\postFDR_{\hat{H}}(\vphi_{\lambda^*-\eps,\hat{H}})<t$. We omit the almost identical proof that for $t<\pi_0$ we have $\postFDR_{\hat{H}}(\vphi_{\lambda^*+\eps,\hat{H}})>t$. From these two bounds one deduces that $\hat{\lambda}\in (\lambda^*-\eps,\lambda^*+\eps)$, implying \cref{eqn:hatLambdaConverges}.
	
		By \cref{lem:LHatApproximatedBylinfty}, there exist $\xi_N,\delta_N\to 0$, such that with probability tending to 1
	\[ \#\braces{i: 1\leq i \leq N,~ \abs{\hat{\ell}_i(X)-\ell_i^\infty(X)}>\xi_N}\leq N\delta_N,\] and we observe that
	\begin{align*}
		\postFDR_{\hat{H}}(\vphi_{\lambda^*-\eps,\hat{H}}) &= \frac{\sum \hat{\ell}_i \II\braces{\hat{\ell}_i<\lambda^*-\eps}}{1\vee (\sum \II\braces{\hat{\ell}_i<\lambda^*-\eps})}\\
		& \leq \frac{\sum \hat{\ell}_i \II\braces{\hat{\ell}_i<\lambda^*-\eps,\abs{\hat{\ell}_i-\ell_i^\infty}\leq \xi_N}}{1\vee (\sum \II\braces{\hat{\ell}_i<\lambda^*-\eps,\abs{\hat{\ell}_i-\ell_i^\infty}\leq \xi_N})} + \frac{ \#\braces{i: \abs{\hat{\ell}_i-\ell_i^\infty}>\xi_N}}{\#\braces{i : \hat{\ell}_i<\lambda^*-\eps}} \\
		& \leq \frac{\sum \ell_i^\infty \II\braces{\ell_i^\infty < \lambda^*-\eps + \xi_N}}{1\vee(\sum \II\braces{\ell_i^\infty<\lambda^*-\eps-\xi_N,\abs{\hat{\ell}_i-\ell_i^\infty}\leq \xi_N}}
+\xi_N+ \frac{ \#\braces{i: \abs{\hat{\ell}_i-\ell_i^\infty}>\xi_N}}{\#\braces{i : \hat{\ell}_i<\lambda^*-\eps}}.
	\end{align*}
		Since $\lambda^*> t$ (the proof of) \cref{lem:O(N)rejections} implies that for some $c>0$ and for $\eps>0$ small enough, $\#\braces{i : \hat{\ell}_i<\lambda^*-\eps}>cN$ with probability tending to 1. 
We also lower bound the denominator in the first term of the final line by $\#\braces{i : \ell_i^\infty <\lambda^*-\eps-\xi_N}-\#\braces{i : \abs{\hat{\ell}_i-\ell_i^\infty}>\xi_N}$; for $\eps,\xi_N,c'$ small enough note that $\#\braces{i : \ell_i^\infty <\lambda^*-\eps-\xi_N}>c'N$ with probability tending to 1 by ergodicity (i.e.\ applying \cref{lem:ErgodicTheorems} with $g(x)=\II\braces{x<\lambda^*-\eps-\xi}$ for some $\xi>\xi_N$), using that $\Pi(\ell_i^\infty<\lambda^*-\eps-\xi_N)>0$.  It follows that for an event $C_N$ of probability tending to 1, $\postFDR_{\hat{H}}(\vphi_{\lambda^*-\eps,\hat{H}})$ is upper bounded by
	\begin{align*}  
	& \II_{C_N^c}+\frac{\sum \ell_i^\infty  \II\braces{\ell_i^\infty < \lambda^*-\eps + \xi_N}}{\sum\II\braces{ \ell_i^\infty <\lambda^*-\eps-\xi_N}} \brackets[\Big]{1+O\brackets[\Big]{\frac{\#\braces{i : \abs{\hat{\ell}_i-\ell_i^\infty}>\xi_N}}{\#\braces{i : \ell_i^\infty <\lambda^*-\eps-\xi_N}}}}+ \xi_N+\delta_N/c \\
		 \leq&\frac{\sum \ell_i^\infty  \II\braces{\ell_i^\infty < \lambda^*-\eps + \xi_N}}{\sum\II\braces{ \ell_i^\infty <\lambda^*-\eps-\xi_N}}+o_p(1).\end{align*} 
	Again using the ergodicity result \cref{lem:ErgodicTheorems}, we have that, for fixed $\xi>0$, \begin{align*}
			\frac{1}{N}\sum_{i=1}^N \ell_i^\infty \II\braces{\ell_i^\infty<\lambda^*-\eps+\xi}&\to E_H [\ell_i^\infty\II\braces{\ell_i^\infty <\lambda^*-\eps+\xi}] \quad \text{in probability,}\\
			\frac{1}{N}\sum_{i=1}^N \II\braces{\ell_i^\infty<\lambda^*-\eps-\xi}&\to \Pi_H(\ell_i^\infty<\lambda^*- \eps-\xi)>0 \quad\text{in probability,} 
		\end{align*} 
so that we may apply Slutsky's lemma (e.g.\ Lemma 2.8 in \cite{vdV98}) to deduce that for $N$ large enough
	\begin{equation*}\begin{split} \frac{\sum_{i=1}^N \ell_i^\infty \II\braces{\ell_i^\infty<\lambda^*-\eps+\xi_N}}{\sum_{i=1}^N \II\braces{\ell_i^\infty<\lambda^*-\eps-\xi_N}} &\leq \frac{\sum_{i=1}^N \ell_i^\infty \II\braces{\ell_i^\infty<\lambda^*-\eps+\xi}}{\sum_{i=1}^N \II\braces{\ell_i^\infty<\lambda^*-\eps-\xi}} \\ &\leq    \frac{E_H [\ell_i^\infty\II\braces{\ell_i^\infty <\lambda^*-\eps+\xi}]}{\Pi_H(\ell_i^\infty<\lambda^*-\eps-\xi)}+o_p(1).\end{split} \end{equation*}
	Finally we note that 
	\[ 	\frac{E_H [\ell_i^\infty\II\braces{\ell_i^\infty <\lambda^*-\eps+\xi}]}{\Pi_H(\ell_i^\infty<\lambda^*-\eps-\xi)} =E_H\sqbrackets[\big]{\ell_i^\infty \mid \ell_i^\infty <\lambda^*-\eps+\xi} \frac{\Pi_H(\ell_i^\infty<\lambda^*-\eps+\xi)}{\Pi_H(\ell_i^\infty <\lambda^*-\eps-\xi)}.	\] 
	Uniformly in $\xi$ satisfying $0<\xi<\eps/2$, we have by monotonicity
	\[ E_H[\ell_i^\infty \mid \ell_i^\infty <\lambda^*-\eps+\xi] \leq E_H[\ell_i^\infty \mid \ell_i^\infty <\lambda^*-\eps/2]<\min(t,\pi_0).\] Observe also that $\Pi_H(\lambda^*-\eps-\xi\leq  \ell_i^\infty <\lambda^*-\eps+\xi)\to 0$ as $\xi\to 0$ by the continuity of the distribution function of $\ell_i^\infty$, while $\Pi_H(\ell_i^\infty <\lambda^*-\eps-\xi)$ is bounded away from zero for $\lambda^*-\eps-\xi$ bounded away from zero. It follows that by choosing $\xi=\xi(\eps)$ small enough we may ensure that 
		\[E_H[\ell_i^\infty \mid \ell_i^\infty <\lambda^*-\eps+\xi] \frac{\Pi_H(\ell_i^\infty<\lambda^*-\eps+\xi)}{\Pi_H(\ell_i^\infty <\lambda^*-\eps-\xi)}<\min(t,\pi_0). \] We conclude, as claimed, that $\postFDR_{\hat{H}}(\vphi_{\lambda^*-\eps,\hat{H}})<\min(t,\pi_0)$ with probability tending to 1.
\end{proof}
\begin{proof}[Proof of \cref{thm:BFNRoptimality}] 
Define $\lambda^*$ as in \cref{lem:LambdaHatConcentrates}. In the case $\lambda^*=1$, one shows 
that $\hat{\vphi}$ rejects all but $o_p(N)$ of the hypotheses. It follows that, asymptotically, its mTDR is close to that of the procedure which rejects all null hypotheses, which trivially has the best mTDR of any procedure. We omit the proof details in this case and henceforth assume that $\lambda^*<1$, or equivalently (in view of \cref{lem:LambdaHatConcentrates}) that $t<\pi_0$.

We compare $\hat{\vphi}$ to the `oracle' procedure $\vphi_{\lambda^*,H}$, which we will argue has optimal multiple testing properties. For $\eps_N>0$ we may decompose
\[ \II\braces{\ell_i< \lambda^*}\leq \II\braces{\lambda^*-\eps_N\leq \ell_i <\lambda^*} +\II\braces{\hat{\ell_i}< \hat{\lambda}}+ \II\braces{\hat{\lambda}<\lambda^*-\eps_N/2}+\II\braces{\hat{\ell}_i-\ell_i>\eps_N/2}.\]
\Cref{lem:LambdaHatConcentrates} tells us that $\hat{\lambda}$ tends to $\lambda^*$ in probability, so that $\Pi(\hat{\lambda}<\lambda^*-\eps_N/2)\to 0$ for $\eps_N$ tending to zero slowly enough,
and \cref{lem:HatlVslErrors} tells us that $\#\braces{i: \abs{\hat{\ell}_i-\ell_i}>\eps_N/2}/N\to 0$ in probability, again for $\eps_N$ tending to zero slowly enough. \Cref{lem:LHatApproximatedBylinfty} tells us that there exist $\xi_N\to 0$ such that $\#\braces{i : \abs{\ell_i-\ell_i^\infty}>\xi_N}/N\to 0$ in probability, and \cref{lem:liinftyFullSupport} tells us that under the conditions of \cref{thm:BFNRoptimality} the distribution function of $\ell_i^\infty$ is continuous, so that as $N\to\infty$
\[ E\#\braces{i : \lambda^*-\eps_N\leq \ell_i<\lambda^*}/N\leq E\#\braces{i : \abs{\ell_i-\ell_i^\infty}>\xi_N}/N+\Pi( \lambda^*-\eps_N-\xi_N\leq \ell_i^\infty<\lambda^*+\xi_N)\to 0.\]  
We deduce that
\[E[\#\braces{i: \theta_i=1,\hat{\ell}_i< \hat{\lambda}}]\geq E[\#\braces{i: \theta_i=1,\ell_i< \lambda^*}]-o(N),\] so that, dividing each side by $E\#\braces{i : \theta_i=1}=N\pi_1$,
\[ \mTDR_H(\hat{\vphi})\geq \mTDR_H(\vphi_{\lambda^*,H})-o(1).\]
Next we consider the mFDR. A similar decomposition to those above yields that
\begin{align*}
	E[\#\braces{i: \theta_i=0,\hat{\ell}_i<\hat{\lambda}}]&\leq E[\#\braces{i : \theta_i=0,\ell_i<\lambda^*}]+o(N),\\
	E[\#\braces{i : \hat{\ell}_i<\hat{\lambda}}]&\geq E[\#\braces{i: \ell_i<\lambda^*}]-o(N).
\end{align*} 
One also has (by comparison to $\ell_i^\infty$ as above, or by adapting the proof of \cref{lem:O(N)rejections}) that $E\#\braces{i: \ell_i<\lambda^*+\eps_N}\geq cN$ for some $c>0$, so that a Taylor expansion yields 
\[ \mFDR_H(\hat{\vphi})\leq \frac{E\#\braces{i: \theta_i=0,\ell_i<\lambda^*}+o(N)}{E\#\braces{i : \ell_i<\lambda^*}-o(N)}\leq \mFDR_H(\vphi_{\lambda^*,H})+o(1).\]

Define $g(x)=\sup\braces{\mTDR_H(\psi) : \mFDR_H(\psi)\leq x}$. Trivially $\mTDR_H(\hat{\vphi})\leq g(\mFDR_H(\hat{\vphi}))$, and hence the following chain of equalities (justified below) proves the first claim of the \namecref{thm:BFNRoptimality}:
\begin{align*} \mTDR_H(\hat{\vphi})&\geq \mTDR_H(\vphi_{\lambda^*,H})-o(1) \\
	&\geq  g\brackets[\big]{\mFDR_H(\vphi_{\lambda^*,H})}-o(1)\\
	&\geq g\brackets[\big]{\mFDR_H(\hat{\vphi})-o(1)}-o(1)\\ 
	&\geq g\brackets[\big]{\mFDR_H(\hat{\vphi})}-o(1).
\end{align*}
The first line was proved above. The second is a consequence of an optimality property for the class $(\vphi_{\lambda,H}:\lambda \in [0,1])$ given by \cref{lem:OptimalityOfOracleClass} in \cref{sec:LemmasForBTDRoptimality}. The third line then follows from what was proved above, and the final line follows by a continuity-type result for $g$ given by \cref{lem:OptimalTestNearContinuity}.

It remains to prove the second claim of the theorem. This will follow, with the same arguments as above, from proving that $\mFDR_H(\vphi_{\lambda^*,H})\geq t-o(1)$. 
Observe that, using \cref{lem:LHatApproximatedBylinfty} as above, one can show
\begin{align*} 
	E[\sum_{i\leq N}  \ell_i\II\braces{\ell_i<\lambda^*}] &= E\sum_{i\leq N}[\ell_i^\infty \II \braces{\ell_i^\infty<\lambda^*}]+o(N) \\
	E[\sum_{i\leq N} \II\braces{\ell_i<\lambda^*}]&=E[\sum_{i\leq N} \II\braces{\ell_i^\infty <\lambda^*}]+o(N).
\end{align*}
Stationarity of the HMM implies that  \begin{align*}E[\sum_{i\leq N} \ell_i^\infty \II\braces{\ell_i^\infty<\lambda^*}]&=NE[\ell_1^\infty \mid \ell_1<\lambda^*]\Pi(\ell_1<\lambda^*),\\
	E[\sum_{i\leq N} \II\braces{\ell_i^\infty<\lambda}]&=N\Pi(\ell_1^\infty<\lambda^*), \end{align*}
and hence by definition of $\lambda^*$ (recall we have assumed $t<\pi_0$) \[ E\sum_{i\leq N} (\ell_i^\infty - t)\II \braces{\ell_i^\infty<\lambda^*}=0.\]
Returning to the $\ell$--values themselves and using also \cref{lem:liinftyFullSupport} to see that $\Pi(\ell_1^\infty<\lambda^*)>0$, we deduce that
\begin{align*} N^{-1}E[\sum_{i\leq N}  (\ell_i -t)\II\braces{\ell_i<\lambda^*} ] &\to 0,\\
	N^{-1} E\sum_{i\leq N} \II\braces{\ell_i<\lambda^*} \to \Pi(\ell_1^\infty<\lambda^*)>0,
\end{align*}
and we may rearrange to see that $\mFDR_H(\vphi_{\lambda^*,H})\geq t-o(1)$.
\end{proof}

\subsection{Proofs: Supremum Norm Estimation of HMM Parameters}\label[section]{sec:proofs:SupNormEstimation}
We construct the estimators of \cref{thm:LinftyEstimation} using a spectral kernel density estimation method.
Let $K$ be a bounded Lipschitz-continuous function, supported in $[-1,1]$, such that if we define
 \begin{equation} \label[equation]{eqn:def:KL}  \begin{split} K_L(x,y)&=2^L K(2^L(x-y)), \\ 
		K_L[f](x)&=\int K_L(x,y) f(y) \dy, \end{split} \end{equation}
then we have, for any $f\in C^s(\RR)$,
	\begin{equation}\label[equation]{eqn:KLfapproximatesf}
		\norm{f-K_L[f]}_{\infty} \leq C\norm{f}_{C^s} 2^{-Ls}. \end{equation}
Note that such a function, a `bounded convolution kernel of order $s$', exists, see \cite{Tsybakov09} (in particular, to ensure $K$ is Lipschitz, one builds the kernel using a Gegenbauer basis with parameter $\alpha>1$ as in Section 1.2.2 thereof). We also note here that for some $C=C(\Hh)$,
\begin{equation} \label[equation]{eqn:lambdamaxBounded} \max_j \norm{K_L[f_j]}_\infty \leq 2\norm{K}_\infty \max_j\norm{f_j}_\infty\leq C\end{equation} since $\int_{-1}^1 \abs{K(x)}\dx \leq 2\norm{K}_\infty$.
[Recall that a constant $C=C(\Hh)$ depends only on certain bounds for the parameter $H=(Q,\pi,f_1,\dots,f_J)$ as described in \cref{sec:Uniformity}. In fact, as with the above, we allow such a constant to also depend on the kernel $K$ since this kernel can be chosen independent of $H$. Similarly, we will permit such a constant $C$ to depend on the functions $h_1,\dots,h_{L_0}$ and the sets $\DD_N$ of \cref{alg:Linfty}.]

The premise of the estimation algorithm comes from the following \namecref{lem:AlgorithmIntuition}, which adapts ideas found in \cite{AHK12} and \cite{Lehericy18}.
 \begin{lemma}\label[lemma]{lem:AlgorithmIntuition} For $L_0\in\NN$, let $h_1,\dots h_{L_0}$ be arbitrary functions. Define, for data $X$ from the HMM \cref{eqn:def:model},
	\begin{align}\label[equation]{eqn:def:Mx} M^x&\equiv M^{x,L_0,L}:=(E_H[h_l(X_1)K_L(x,X_2)h_m(X_3)]_{l,m\leq L_0})\in \RR^{L_0\times L_0}, \\
		P&\equiv P^{L_0}:= (E_H[h_l(X_1)h_m(X_3)]_{l,m\leq L_0}) \in \RR^{L_0\times L_0},\\
		D^x&\equiv D^{x,L}:=\diag\brackets[\big]{ K_L[f_j](x)_{j\leq J}} \in \RR^{J\times J}, \\
		O&\equiv O^{L_0}:=(E_H[h_l(X_1) \mid \theta_1=j]_{l\leq L_0, j\leq J})\in\RR^{L_0\times J}.
	\end{align} 
Then  
\begin{align*} M^x&=O \diag (\pi) Q D^x QO^\transpose, \quad \text{and} \\
	 P&= O \diag (\pi) Q^2 O^\transpose. \end{align*}
	If $V\in\RR^{L_0\times J}$ is such that $V^\transpose PV$ is invertible (it suffices to assume $PV$ has rank $J$, which holds under the assumption that $P$ has rank $J$ if the columns of $V$ consist of orthonormal right singular vectors of $P$, or any other orthonormal basis of the column space of $P$) then the matrix 
	\begin{equation} B^x\equiv B^{x,L_0,L}:=  (V^\transpose PV)^{-1}  V^\transpose M^x V  \end{equation}
satisfies 
\begin{equation} \label[equation]{eqn:BsimultaneouslyDiagonalisable} B^x=(QO^\transpose V)^{-1} D^x (QO^\transpose V),	\end{equation}
	so that the matrices $(B^x : x\in \RR)$ are diagonalisable \emph{simultaneously}, with $B^x$ having eigenvalues $(D^x_j: j\leq J) =(K_L[f_j](x): j\leq J).$ 
\end{lemma}
\begin{proof}
	Conditioning on $(\theta_1,\theta_2,\theta_3)$, we see 
	\begin{align*}
	M^x_{l,m}&=\sum_{a,b,c}\Pi_H(\theta_1=a,\theta_2=b,\theta_3=c)E_H[h_l(X_1) K_L(x,X_2) h_m(X_3) \mid \theta_1=a,\theta_2=b,\theta_3=c]  \\
	&= \sum_{a,b,c} \pi_a Q_{a,b} Q_{b,c} O_{l,a}  O_{m,c}E_{X\sim f_b}[K_L(x,X)] \\
	&= (O \diag (\pi) Q D^x QO^\transpose)_{l,m},
	\end{align*}
	and similarly we have
	\begin{align*}
	P&=(\sum_{a,b,c} \pi_a Q_{a,b}Q_{b,c} O_{l,a} O_{m,c})_{l,m}=O \diag (\pi) Q^2 O^\transpose.
	\end{align*}
	
	Next, note that if $V^\transpose PV$ is invertible then so is $QO^\transpose V$ (since  $V^\transpose PV= V^\transpose O \diag (\pi) Q (Q O^\transpose V),$ and a product $AB$ of square matrices is invertible if and only if each of $A$ and $B$ are). The result \cref{eqn:BsimultaneouslyDiagonalisable} then follows from the expressions for $P$ and $M^x$. \end{proof}
\Cref{lem:AlgorithmIntuition} suggests estimating the eigenvalues $K_L[f_j](x)$ of $B^x$ by using empirical versions of $V$, $P$ and $M^x$, an idea which is implemented in the following algorithm. 
The algorithm requires as inputs functions $h_1,\dots h_{L_0}$ and sets $\DD_N$ with certain properties; the existence of suitable inputs is discussed in the remarks thereafter. We introduce notation for the ``eigen-separation'' of a matrix $B\in \RR^{J\times J}$ with eigenvalues $\lambda_1,\dots, \lambda_J$: 
\begin{equation}
	\label[equation]{eqn:def:sep}
	\sep(B)= \min_{i\neq j}\abs{\lambda_i-\lambda_j}.
\end{equation}
Recall that $\sigma_J(B)$ denotes the $J$th largest singular value of $B$. 
\begin{algorithm}[H] 
	\caption{Kernel density estimator}
	\label[algorithm]{alg:Linfty}
	\begin{algorithmic}
		\Input{\begin{itemize} \item Data $(X_n : n\leq N+2)$ drawn from the HMM \cref{eqn:def:model}.
			\item Functions $h_1,\dots h_{L_0}$, uniformly bounded, such that $O=(E[h_l(X_1) \mid \theta_1 = j]_{l\leq L_0, j\leq J})$ is of rank $J$, with $\sigma_J(O)$ bounded away from 0 uniformly in $N$, at least for $N$ large enough. 
			\item Finite sets $\DD_N\subseteq \braces{(a,u)\in\RR^{J(J-1)/2}\times\RR^{J(J-1)/2} : \sum \abs{a_i} \leq 1}$ such that $\max_{(a,u)\in\DD_N} \sep(B^{a,u})$ is bounded away from 0 uniformly in $N$, at least for $N$ large enough, where $B^{a,u}=\sum a_i B^{u_i}$ for $B^x$ as in \cref{lem:AlgorithmIntuition} for some $V$.\end{itemize}
		}
		\vspace{5pt}
		\Estimate {the matrices $P,$ $(M^x,~x\in\RR)$ of \cref{lem:AlgorithmIntuition} by taking empirical averages:
			\begin{align*} \hat{P}=\hat{P}^{L_0}&=(N^{-1}\sum\nolimits_{n\leq N} h_l(X_n)h_m(X_{n+2}))_{l,m\leq L_0},\\ \hat{M}^x=\hat{M}^{x,L_0,L}&=(N^{-1}\sum\nolimits_{n\leq N} h_l(X_n)K_L(x,X_{n+1})h_m(X_{n+2}))_{l,m\leq L_0}.\end{align*}
			Let $\hat{V}=\hat{V}^{L_0}\in \RR^{L_0\times J}$ be a matrix of orthonormal right singular vectors of $\hat{P}$ (fail if $\hat{P}$ is of rank less than $J$).} 
		\State {\textbf{set}, for $x\in\RR$ and for $a,u\in \RR^{J(J-1)/2}$ \[\hat{B}^x=\hat{B}^{x,L_0,L}:= (\hat{V}^\transpose \hat{P}\hat{V})^{-1} \hat{V}^\transpose \hat{M}^x\hat{V}, \quad  \hat{B}^{a,u}:=\sum a_i \hat{B}^{u_i}.\]}
		\State {\textbf{choose} $\hat{R}$ of normalised columns diagonalising $\hat{B}^{\hat{a},\hat{u}}$, where $(\hat{a},\hat{u})\in\argmax_{\DD_N} \sep(\hat{B}^{a,u})$ (fail if $\hat{B}^{\hat{a},\hat{u}}$ is not diagonalisable).}
		\Output{$(\hat{f}_j : j\leq J)$, where, defining
			\[ \tilde{f}^L_j(x)=(\hat{R}^{-1}\hat{B}^x \hat{R})_{jj}, \] we set \[\hat{f}_j(x)=\begin{cases} \tilde{f}_j^L(x) & \abs{\tilde{f}_j^L(x)}\leq N^\alpha \\ N^\alpha \sign(\tilde{f}_j^L(x)) & \text{otherwise}, \end{cases} \] for $\alpha>0$ arbitrary and for $L$ such that $2^L \asymp (N/\log N)^{1/(1+2s)}.$ [The in-probability result \cref{eqn:InProbabilitySupNormResult} also holds for $\tilde{f}^L_j$.]} 
	\end{algorithmic}
\end{algorithm} 
\begin{remarks}
	\begin{enumerate}[i.]
		\item  For notational convenience, we have considered observing $N+2$ data points $X_1,\dots,X_{N+2}$ so that we can form $N$ triples of consecutive observations; the proofs go through for the original $N$ data points by adjusting constants. 
	\item Under \cref{ass:LinearIndependence}, $h_1,\dots h_{L_0}$ can be chosen without knowledge of the parameters, for example by letting $L_0$ tend to infinity arbitrarily slowly and taking the $h_l$ to be indicator functions of the first $L_0$ of a countable collection of sets generating the Borel $\sigma$--algebra (see \cref{lem:ExistenceOfSuitableHl}, in \cref{sec:WellDefinednessOfEstimators}). In principle, $L_0=J$ is sufficient to achieve $O$ of rank $J$, but without further assumptions,  the appropriate functions $h_1,\dots h_J$ will necessarily depend on the unknown parameters. In the case $J=2$, it suffices to assume in addition to the other conditions of \cref{thm:LinftyEstimation} that $P_{X\sim f_1}(X\in A)\neq P_{X\sim f_2} (X\in A)$ for some known $A$, by taking $h_1=1,h_2=\II_A$. 
	\item \Cref{lem:AlgorithmIntuition} implies that the condition on $\DD_N$ is independent of $V$ provided $V$ is such that $V^\transpose P V$ is invertible. \Cref{lem:ExistenceOfHatAHatU}, the proof of which uses only that $f_1,\dots, f_J$ are distinct, shows that the choice $V=\hat{V}$ is suitable with probability tending to 1 and that $\DD_N$ can be chosen independent of the parameters, for example by taking a cartesian product of increasing dyadic sets of rationals. In the case $J=2$, the description of the algorithm simplifies, in that necessarily $\hat{a}=1\in\RR^1$. A corresponding simplification also works in the general $J$ state case if one is willing to assume that there exists $x_0\in\RR$ for which the values $f_j(x_0),~j\leq J$ are all distinct, in that one may define $\hat{R}$ as diagonalising $\hat{B}^{\hat{x}}$ where $\hat{x}$ maximises $\sep(\hat{B}^x)$ over $x$ in (some finite increasing sieve in) $\RR$.
	\item 	\Cref{lem:EventA,lem:EmpiricalEigengap} prove that with probability tending to 1 $\hat{P}$ has rank $J$ and that $\hat{B}^{\hat{a},\hat{u}}$ is diagonalisable, and hence that the outputs $\hat{f}_j$ are well-defined.
	\item 	Since the $f_j$ are assumed H\"older continuous, and satisfy tail bounds, one could in fact calculate $\hat{f}_j(x)$ only for $x$ in some finite set, then construct estimators $\check{f}_j$ via interpolation, in order to ease computation.
	\end{enumerate}
\end{remarks}
 
\begin{proof}[Proof of \cref{thm:LinftyEstimation}]
 Construct $\hat{f}_j,\tilde{f}_j^L$ using \cref{alg:Linfty}.	Continuity of these functions follows from continuity of the map $x\mapsto \hat{B}^x$, which in turn follows from that of the map $x\mapsto \hat{M}^x$, proved in \cref{lem:EventA}.
	Observe also that $\norm{f_j}_\infty<\infty$ for all $j\leq J$, so that for $N$ large enough that $\norm{f_j}_\infty \leq N^a$ we have
	\[ \norm{\hat{f}_j-f_{\tau(j)}}_\infty \leq \norm{\tilde{f}_j^L-f_{\tau(j)}}_\infty,\] hence for the in-probability result it suffices to prove \cref{eqn:InProbabilitySupNormResult} with $\tilde{f}_j=\tilde{f}^L_j$ in place of $\hat{f}_j$.
	
For a constant $c>0$, define the event
	\begin{equation}\label[equation]{eqn:eventA}
		\Aa= \braces{ \norm{\hat{P}-P}\leq cL_0 r_N,~ \norm{\hat{M}^x-M^x}\leq cL_0^2 r_N ~\forall x\in \RR}.
	\end{equation} This is indeed a measurable event, and for suitable $c=c(\kappa,\Hh)$ it has probability at least $1-N^{-\kappa}$, by \cref{lem:EventA}, which also tells us that $\hat{V}^\transpose P\hat{V}$ is invertible on $\Aa$ and that, defining \begin{equation*} \tilde{B}^x:=(\hat{V}^\transpose P\hat{V})^{-1} \hat{V}^\transpose M^x \hat{V}, \end{equation*}
 we have, for some $C$ depending on $\Hh$ and on the constant $c$ of event $\Aa$, \begin{equation*}
\II_\Aa	\sup_{x\in \RR} \norm{\tilde{B}^x-\hat{B}^x}\leq CL_0^2 r_N.
\end{equation*}

\Cref{lem:AlgorithmIntuition} tells us (on $\Aa$) that $\tilde{B}^x=(QO^\transpose \hat{V})^{-1} D^x QO^\transpose \hat{V}$, and we write $\tilde{R}$ for a matrix whose columns are those of $QO^\transpose \hat{V}$ but scaled to have unit Euclidean norm, which thus diagonalises $\tilde{B}^x$ for all $x$. By \cref{lem:RhatEstimatesR}, $\norm{\hat{R}-\tilde{R}_\tau}\leq CL_0^{7/2} r_N$ on $\Aa$ for some permutation $\tau$, where $\tilde{R}_\tau$ is obtained by permuting the columns of $\tilde{R}$ according to $\tau$. Next we apply \cref{lem:C.4AHK12} with $\Tt=\RR$, $A_x=\tilde{B}^x,$ $\hat{A}_x=\hat{B}^x$, $R=\tilde{R}$. Noting that $\norm{\tilde{R}^{-1}}\leq C'L_0^{1/2}$ and $\kappa(\tilde{R})\leq C'L_0$ for some constant $C'=C'(\Hh)$ (see  \cref{lem:kappaQObound}), and that the constant $\lambdamax$ of the \namecref{lem:C.4AHK12} is bounded by a constant depending only on $\Hh$ (see \cref{eqn:lambdamaxBounded}), we deduce that
\[\sup_x\max_j \abs{\tilde{f}_j^L(x)-K_L[f_{\tau(j)}](x)} \leq c'L_0[L_0^2 r_N +L_0^{1/2}L_0^{7/2} r_N]\leq c''L_0^5  r_N,\] for some constants $c',c''$.   
 The in-probability result \cref{eqn:InProbabilitySupNormResult} follows, since the choice of $L$ ensures by \cref{eqn:KLfapproximatesf} that $\norm{f_{\tau(j)}-K_L[f_{\tau(j)}]}_\infty \leq C'' r_N$ on $\Aa$ for some $C''$,  so that for a suitable constant $C$, \[\Pi_H(\norm{\tilde{f}_j-f_{\tau(j)}}_\infty >CL_0^5 r_N) \leq \Pi(\Aa^c)\leq N^{-\kappa}\to 0.\]

For the in-expectation result \cref{eqn:InExpectationSupNormResult}, observe that by truncating at $\pm N^\alpha$ we have ensured that
\[ E_H \norm{\hat{f}_j-f_{\tau(j)}}_\infty \leq CL_0^5 r_N + 2N^\alpha \Pi_H(\Aa^c).\]
Choosing $c=c(\kappa,\Hh)$ in the definition of the event $\Aa$ corresponding to some $\kappa \geq s/(1+2s) + \alpha$ concludes the proof.
\end{proof}

\begin{proof}[Proof of \cref{prop:ExistenceOfEstimators}]
	Let $\hat{f}_0,\hat{f}_1,\hat{Q},\hat{\pi}$ be estimators which satisfy \begin{equation}\label[equation]{eqn:AssumedConsistentUpToPerm}\Pi_H\brackets{\norm{\hat{f}_{0}-f_{\tau(0)}}+\norm{\hat{f}_{1}-f_{\tau(1)}}+\norm{\hat{Q}-Q_{\sigma,\sigma}}_F+\norm{\hat{\pi}-\pi_\sigma}>C\eps_N}\to 0\end{equation} 
	for some permutations $\tau,\sigma$ and a constant $C>0$, with $Q_{\sigma,\sigma}$ defined by permuting the rows and columns of $Q$, and $\pi_\sigma$ defined similarly.	 The existence of suitable $\hat{f}_0,\hat{f}_1$ is given by \cref{thm:LinftyEstimation}, and the existence of suitable $\hat{Q},\hat{\pi}$ is proved by results in \cite[Appendix C]{dCGlC17} (and by arguments as in \cite[Section 8.6]{dCGL16} to accelerate the possibly slow rate).
	Moreover, the estimators of \cite{dCGlC17} are constructed using a spectral method, so that one may in fact assume $\sigma=\tau$. [One could also ``align'' $\sigma$ and $\tau$ by hand, by noting that by ergodicity the invariant density $f_\pi$ can be estimated at the rate $r_N$ using a standard kernel density estimator, and permuting rows and columns of $\hat{Q}$ and $\hat{\pi}$ so that $\sum \hat{\pi}_i \hat{f}_i$ is close to this kernel density estimator; linear independence of the $f_i$ ensures that this alignment method works.]	
	
	Next, under the assumption $\pi_0>\pi_1$, define $\check{f}_j=\hat{f}_{\hat{\tau}(j)},$ $\check{Q}=\hat{Q}_{\hat{\tau},\hat{\tau}}$ and $\check{\pi}=\hat{\pi}_{\hat{\tau}}$, where $\hat{\tau}(0)=1-\hat{\tau}(1)=\II\braces{\hat{\pi}_1>\hat{\pi}_0}$. Consistency of $\hat{\pi}$ implies that $\hat{\tau}$ consistently estimates the permutation $\tau=\sigma$ of \cref{eqn:AssumedConsistentUpToPerm}, hence
	\[\begin{split} &\Pi_H\brackets{\norm{\check{f}_{0}-f_{0}}+\norm{\check{f}_{1}-f_{1}}+\norm{\check{Q}-Q}_F+\norm{\check{\pi}-\pi})>C\eps_N}  \\  & \quad \leq \Pi_H (\hat{\tau}\neq \tau)+\Pi_H\brackets{\norm{\hat{f}_{0}-f_{\tau(0)}}+\norm{\hat{f}_{1}-f_{\tau(1)}}+\norm{\hat{Q}-Q_{\tau,\tau}}_F+\norm{\hat{\pi}-\pi_\tau}>C\eps_N}\to 0.\end{split} \]

	For the second case, we want to define $\hat{\tau}(0)=\II\braces{\limsup_{x\uparrow x^*} (\hat{f}_0/\hat{f}_1)(x)>1}$ and proceed similarly, but the compact support of $K$ means that $\hat{f}_1(x)=\hat{f}_0(x)=0$ for $x>2^{-L}+\max_k X_k,$ and the right side may be strictly smaller than $x^*$. Instead, noting that necessarily $\Pi(X_1\leq x^*)>0$ and assuming without loss of generality that $x^*>0$, we set $\tilde{X}_n=X_n\II\braces{X_n\leq x^*}$ and define 
	\begin{align*} \hat{\tau}(0)=1-\hat{\tau}(1)&= \II\braces{\hat{f}_0(M_N)>\hat{f}_1(M_N)}, \\ 
		M_N &=\max_{i\leq \log N}(\tilde{X}_i);
	\end{align*} 
note that by construction we have $\hat{f}_{\hat{\tau}(1)}(M_N)\geq \hat{f}_{\hat{\tau}(0)}(M_N)$.
We show that $\norm{\hat{f}_{\hat{\tau}(1)}-f_0}_\infty>C\eps_N$ on an event $A_N$ of probability tending to 1; it will follow from \cref{eqn:AssumedConsistentUpToPerm} that $\hat{\tau}\equiv \hat{\tau}^{-1}=\tau$ on $A_N$, and the result will follow.
	
	The variables $\tilde{X}_i,~i\leq N$ have a density with respect to the measure $\mu$ defined by adding an atom at 0 to Lebesgue measure. Let $u$ be as in \cref{thm:BFDRcontrol}, so that $u>1+\nu^{-1}$ and $\eps_N(\log N)^u\to 0$ for $\nu$ as in \cref{ass:ExistsCNu}. The proof of \cref{lem:ControlOf1/f} shows that with probability tending to 1 we have $f_1(M_N)\geq \min_{i\leq \log N} (f(\tilde{X}_i))\geq (\log N)^{-u}$, hence $f_1(M_N)>3C\eps_N$. We also note that $M_N\uparrow x^*$ almost surely, so that $f_1(M_N)>3f_0(M_N)$ for all $N$ large enough.
	
	Let $A_N$ be an event of probability tending to 1 on which \[f_1(M_N)>3C\eps_N, \quad f_1(M_N)>3f_0(M_N),\quad \norm{\hat{f}_0-f_{\tau(0)}}_\infty\leq C\eps_N, \quad   \norm{\hat{f}_1-f_{\tau(1)}}_\infty\leq C\eps_N,\] whose existence we have just demonstrated. On $A_N$ we have both $\hat{f}_1(M_N)\geq f_{\tau(1)}(M_N)-C\eps_N$ and $\hat{f}_0(M_N)\geq f_{\tau(0)}(M_N)-C\eps_N$ hence
	\[ \begin{split} \hat{f}_{\hat{\tau}(1)}(M_N)&=\max(\hat{f}_{0}(M_N),\hat{f}_{1}(M_N))\geq \max_j (f_j(M_N)-C\eps_N)=f_1(M_N)-C\eps_N>\tfrac{1}{3}f_1(M_N)+C\eps_N\\ &>f_0(M_N)+C\eps_N, \end{split} \]
so that $\norm{\hat{f}_{\hat{\tau}(1)}-f_0}_\infty >C\eps_N$ on $A_N$ as claimed.\qedhere	
\end{proof}

\acks{
The authors would like to acknowledge \'Etienne Roquain, Gloria Buritica and Ramon van Handel for fruitful discussions about this work. 
K.A.\ was supported in this work by grants from the Fondation Mathématique Jacques Hadamard (FMJH). I.C.\ and E.G.\ would like to acknowledge support for this project
	from Institut Universitaire de France. I.C.\ is partly supported by ANR-17-CE40-0001 grant (BASICS).
}

\appendix

\section{Auxiliary Results for \cref{sec:EmpiricalBayesProcedure}} \label[appendix]{sec:LemmasForMultipleTesting}
\subsection{Lemmas for \cref{thm:BFDRcontrol}}\label[appendix]{sec:LemmasForBFDRcontrol}
	Recall $f_\pi=\pi_0f_0+\pi_1f_1$ is the density of each $X_i$, $i\leq N$, in the HMM model \cref{eqn:def:model}.
	\begin{lemma}\label[lemma]{lem:ControlOf1/f}
		Under \cref{ass:fpi} we have, for any $a>1+\nu^{-1}$, \[\Pi_H(\max_{i\leq R} 1/f_\pi(X_i) > R^{a})\to 0\quad \text{as }R\to\infty.\] 
	\end{lemma}
\begin{proof}
For $A=R^a,B=R^b$ with $a,b>0$ to be chosen, we have by a union bound and stationarity
\begin{align*} \Pi_H\brackets[\Big]{\max_{i\leq R} \frac{1}{f_\pi(X_i)} > A} &\leq R\Pi_H\brackets[\big]{f_\pi(X_1)<A^{-1}} \\
&\leq R\int_{-B}^B \II\braces[\big]{f_\pi(x)<A^{-1}}f_\pi(x)\dmu(x) + R\Pi_H\brackets[\big]{\abs{X_1}> B} \\
&\leq R\mu([-B,B])/A + R\Pi_H\brackets[\big]{\abs{X_1}> B}. 
\end{align*}
Since $f_\pi$ is a mixture of the densities $f_0,f_1$, an application of Markov's inequality yields \[\Pi_H\brackets[\big]{\abs{X_1} >B}\leq \max_j P_{X\sim f_j}\brackets[\big]{\abs{X}> B}\leq B^{-\nu} \max_j E_{X\sim f_j} \abs{X}^\nu,\] which is at most a constant times $B^{-\nu}$ by the assumption. Choosing $b>1/\nu$, we have $R\Pi_H(\abs{X_1}\geq B)\to 0$.

Since $B=R^b\geq 1$ and $\mu$ is equal to either to Lebesgue or counting measure, $\mu([-B,B])\leq 2B+1\leq 3B$. Then
\[ R \mu([-B,B])/A \leq 3R^{1+b-a},\]
which tends to zero for $a>1+b$, so that any $a>1+\nu^{-1}$ is permissible. 
\end{proof}

For the following two lemmas recall the definition $\hat{S}_0=\hat{S}_0(t)=\braces{i : \hat{\vphi}_i=1}$, where $\hat{\vphi}$ is as in \cref{def:HatPhi}, so that $\hat{K}=\abs{\hat{S}_0}$ is characterised by \[\frac{1}{\hat{K}} \sum_{i=1}^{\hat{K}} \hat{\ell}_{(i)} \leq t < \frac{1}{\hat{K}+1} \sum_{i=1}^{\hat{K}+1} \hat{\ell}_{(i)}.
\] where, by convention, the left inequality holds if $\hat{K}=0$, and $\hat{\ell}_{(N+1)}=\infty$ so that the right inequality holds if $\hat{K}=N$. Recall the definition \[\hat{t}:=\postFDR_{\hat{H}}(\hat{\vphi})=\frac{1}{\hat{K}} \sum_{i=1}^{\hat{K}} \hat{\ell}_{(i)}.\]
\begin{lemma}\label[lemma]{lem:hatT}
	In the setting of \cref{thm:BFDRcontrol},
$E_H\hat{t} \to \min(t,\pi_0)$. 
	\end{lemma}

\begin{proof}
Since $0\leq \hat{t}\leq 1$, it's enough to show that $\hat{t}\to \min(t,\pi_0)$ in probability. By \cref{lem:LHatValsErgodic}, we have \begin{equation}\label[equation]{eqn:postFDR1toPi0} \frac{1}{N}\sum_{i=1}^N \hat{\ell}_i(X)\to \pi_0 \quad \text{in probability.}\end{equation} 
 By monotonicity of the average of increasing numbers, we have
\[ \hat{t}\leq \frac{1}{N}\sum_{i=1}^N \hat{\ell}_{(i)}=\frac{1}{N}\sum_{i=1}^N \hat{\ell}_i,\] and by construction we note also that $\hat{t}\leq  t$, hence $\hat{t}\leq \min(t,\pi_0)+o_p(1)$.
	
To obtain a matching lower bound, we decompose relative to the event $\Cc=\braces{\hat{K}=N}$. Observe, using \cref{eqn:postFDR1toPi0}, that
		\[ \hat{t}\II_{\Cc} = \II_\Cc \frac{1}{N} \sum_{i=1}^N \hat{\ell}_i\geq \II_\Cc \pi_0- o_p(1).\]
	By definition of $\hat{K}$ we also have
	\[ t\II_{\Cc^c} < \frac{1}{\hat{K}+1} \sum_{i=1}^{\hat{K}+1} \hat{\ell}_{(i)}\II_{\Cc^c}= \frac{\hat{K}}{\hat{K}+1} \hat{t} \II_{\Cc^c}+\frac{\hat{\ell}_{(\hat{K}+1)}}{\hat{K}+1}\II_{\Cc^c},\]
	hence, since $\hat{\ell}_{(\hat{K}+1)}\leq 1$ on $\Cc^c$,
	\[\hat{t}\II_{\Cc^c}>\frac{\hat{K}+1}{\hat{K}}t \II_{\Cc^c}-\frac{\hat{\ell}_{(\hat{K}+1)}}{\hat{K}}\II_{\Cc^c}>t\II_{\Cc^c}-\frac{1}{\hat{K}}.\]
	By \cref{lem:O(N)rejections}, $\hat{K}\to \infty$ in probability for any $t>0$, so that the above display implies
	$ \hat{t}\II_{\Cc^c} > t\II_{\Cc^c} -o_p(1)$ and
	hence
	\[ \hat{t} > t\II_{\Cc^c} + \pi_0\II_{\Cc} -o_p(1)\geq \min(t,\pi_0)-o_p(1), \] proving the lower bound.\end{proof}

Recall the definition of constants $C=C(\Ii)$ from \cref{sec:Uniformity}.
\begin{lemma}\label[lemma]{lem:O(N)rejections} 
	In the setting of \cref{thm:BFDRcontrol}, for all $t>0$, there exists $a=a(t,\Ii)>0$ such that \[ \Pi_H(\abs{\hat{S}_0}>aN)\to 1. \]
\end{lemma}

\begin{proof} 
The definition of $\hat{\lambda}$ trivially implies $\hat{\lambda}\geq t$, so that
\[\braces{i: \hat{\ell}_i<t}\subseteq \braces{i: \hat{\ell}_i<\hat{\lambda}}\subseteq\hat{S}_0.\]	
For $A\in \NN$ write  \[ \ell_i'(X) := \Pi_H(\theta_i=0 \mid X_{i-A},\dots, X_{i+A}), \quad A<i\leq N-A.\] By \cref{lem:LHatApproximatedByl'}, there exist $A=A(t)$ and events $G_N$ of probability tending to 1 such that
\[\braces[\Big]{\#\braces[\big]{i: A<i\leq N-A,~ \abs{\hat{\ell}_i(X)-\ell_i'(X)}>t/2}\leq N\delta_N},\] for some $\delta_N\to 0$.
On $G_N$, we observe that \[\#\braces{i\leq N : \hat{\ell}_i<t} \geq \#\braces{i : A<i\leq N-A,~\ell_i'<t/2}-N\delta_N,\] hence it suffices to show that there exists $c>0$ such that $\#\braces{i: A<i\leq N-A : \ell_i'<t/2}>cN$ with probability tending to 1.

By ergodicity (i.e.\ applying \cref{lem:ErgodicTheorems} with $g(x)=\II\braces{x<t/2}$) we have for any $\eps>0$
\[\Pi_H\brackets[\Big]{\#\braces[\big]{i: A<i\leq N-A : \ell_i'<t/2}>\brackets[\big]{N-2A}\brackets[\big]{\Pi_H(\ell_i'<t/2)-\eps}}\to 1, \] hence it suffices to show that $\Pi_H(\ell'_i<t/2)\neq 0$.
	
Fix $i$ satisfying $A<i\leq N-A$. For $\alpha,\beta\in \braces{0,1}^{A}$ write \[\eta_{\alpha,\beta}=\pi_{\alpha_1}\prod_{a< A} Q_{\alpha_{a},\alpha_{a+1}}Q_{\beta_a,\beta_{a+1}}.\] Introducing the notation  $\theta_{a}^{b}=(\theta_a,\theta_{a+1},\dots,\theta_b)\in \RR^{b+1-a}$, we note that \[\Pi_H(\theta_{i-A}^{i+A}=(\alpha,0,\beta))=\eta_{\alpha,\beta} Q_{\alpha_A,0}Q_{0,\beta_1}, \quad \Pi_H(\theta_{i-A}^{i+A}=(\alpha,1,\beta))=\eta_{\alpha,\beta} Q_{\alpha_A,1}Q_{1,\beta_1}.\]  Define
\begin{align*} p_0&=\sum_{\alpha,\beta \in \braces{0,1}^{A}} Q_{\alpha_A,0}f_0(X_i)Q_{0,\beta_1} \eta_{\alpha,\beta}\prod_{a\leq A} f_{\alpha_a}(X_{i-A+a-1}) f_{\beta_{A+a-1}}(X_{i+a}) \\
	p_1&=\sum_{\alpha,\beta \in \braces{0,1}^{A}} Q_{\alpha_A,1}f_1(X_i)Q_{1,\beta_1} \eta_{\alpha,\beta}\prod_{a\leq A} f_{\alpha_a}(X_{i-A+a-1}) f_{\alpha_{A+a-1}}(X_{i+a}),
\end{align*}
and observe that \[\ell'_i(X)=\Pi_H(\theta_i=0\mid X_{i-A}^{i+A}) = \frac{p_0}{p_0+p_1}.\] Note that each term in the sum defining $p_1$ is at least $\delta^2 f_1(X_i)/f_0(X_i)$ times the corresponding term in the sum defining $p_0$, with $\delta>0$ as in \cref{ass:Qpi}, hence 
\[ p_1\geq p_0\delta^2 \frac{f_1(X_i)}{f_0(X_i)}.\]
In view of \cref{ass:fpi}, assume without loss of generality that there exists $x^*\in\RR\cup\braces{\pm \infty}$ such that $f_1(x)/f_0(x)\to \infty$ as $x\uparrow x^*$. Then we deduce that for some $u=u(t,\delta)>0$, \[\Pi_H(\ell'_i<t/2)\geq \Pi_H\brackets[\Big]{\frac{f_1(X_i)}{f_0(X_i)}> \frac{2-t}{t\delta^2}} \geq \pi_1 P_{X\sim f_1}(x^*-u\leq X\leq x^*)>0,\] as required. \end{proof}

\begin{lemma}\label[lemma]{lem:LHatValsErgodic} In the setting of \cref{thm:BFDRcontrol},
	\[\frac{1}{N}\sum_{i=1}^N \hat{\ell}_i(X) \to \pi_0\] in probability as $N\to \infty$.
\end{lemma}

\begin{proof}
	It is required to prove, for $\eps>0$ arbitrary, that
	\[ \Pi_H\brackets[\Big]{\abs[\Big]{ \frac{1}{N} \sum_{i=1}^N \hat{\ell}_i(X) - \pi_0}>\eps}\to 0.\]
	By \cref{lem:LHatApproximatedByl'}, defining \[ \ell_i'(X) = \Pi_H(\theta_i=0 \mid X_{i-A},\dots, X_{i+A}), \quad A<i\leq N-A,\] there exists $A=A(\eps)$ for which, with probability tending to 1,
	\[ \#\braces{i : A<i\leq N-A,~\abs{\hat{\ell}_i(X)-\ell_i'(X)}>\eps/2}\leq N\delta_N.\]
	On the event on which the last line holds we can decompose:
	\[ \abs[\Big]{\frac{1}{N}\sum_{i=1}^N \hat{\ell}_i(X) - \pi_0}\leq \frac{2A}{N} +\eps/2 +\delta_N+ \frac{1}{N}\abs[\Big]{ \sum_{i=A+1}^{N-A} (\ell_i'(X)-\pi_0) }.\]
Finally, by ergodicity of $\ell_i'(X)$ (see \cref{lem:ErgodicTheorems}) we have
	\[\Pi_H\brackets[\Big]{\frac{1}{N} \abs[\Big]{\sum_{i=A+1}^{N-A} (\ell_i'(X) -\pi_0)} >\eps/4}\leq \Pi_H\brackets[\Big]{\frac{1}{N-2A} \abs[\Big]{\sum_{i=A+1}^{N-A} (\ell_i'(X) - E_H[\ell_i'(X)])} >\eps/4}\to 0,\]
	where we have used that \[E_H[\ell_i'(X)]=E_H \Pi_H(\theta_i=0 \mid X_{i-A},\dots, X_{i+A})=\Pi_H(\theta_i=0)=\pi_0.\] The result follows.
\end{proof}
\begin{lemma}\label[lemma]{lem:LHatApproximatedByl'}
	For $A\in \NN$, define 
	\[ \ell_i'(X) = \Pi_H(\theta_i=0 \mid X_{i-A},\dots, X_{i+A}), \quad A<i\leq N-A. \] 
	 For any fixed $\eps>0$, there exists $A=A(\eps)$ and $\delta_N \to 0$ such that
	\[ \#\braces{i : A<i\leq N-A,~ \abs{\hat{\ell}_i(X)-\ell_i'(X)}> \eps}\leq N\delta_N, \quad \text{with probability tending to 1.}\] 
\end{lemma}
A similar result holds in the limit $A\to \infty$, see \cref{lem:LHatApproximatedBylinfty} below.\\

\begin{proof} 
Essentially, this is a consequence of \cref{lem:HatlVslErrors} and exponential mixing -- hence forgetfulness -- of the Markov chain $\theta$. Precisely, \cref{lem:HatlVslErrors} tells us that there exist events $G_N$ of probability tending to 1 on which 
\[\braces[\Big]{\#\braces[\big]{i\leq N : \abs{\hat{\ell}_i(X)-\ell_i(X)}>\eps_N'}\leq N\delta_N},\] for some $\eps_N'\to 0$; in particular note $\eps_N'<\eps/2$ for $N$ large. Next, we apply \cite[Proposition 4.3.23iii]{CMR05}. Our \cref{ass:Qpi} implies that Assumption 4.3.24 therein holds, so by the consequent Lemma 4.3.25 one sees that the $\rho_0(y)$ in the proposition can be replaced by $\rho=(1-2\delta)/(1-\delta)$. 
Applying the proposition with $j=k-A$ yields
\[ \abs{\Pi_H(\theta_k=0 \mid X_1,\dots,X_n) - \Pi_H(\theta_k=0 \mid X_{k-A},\dots,X_n)}<2\rho^A,\quad k>A.\] Any two-state Markov chain is reversible, hence by time-reversal we similarly obtain 
\[ \abs{\Pi_H(\theta_k=0\mid X_{k-A},\dots,X_n)-\Pi_H(\theta_k=0 \mid X_{k-A},\dots,X_{k+A})}<2\rho^A,\] and hence 
\[ \abs{\ell_k(X)-\ell_k'(X)}<4\rho^A, \quad A<k\leq N-A.\] Choose $A=A(\eps)$ so that $4\rho^A<\eps/2$; then, on $G_N$ and for $N$ large, an application of the triangle inequality yields
\[ \#\braces{i: A<i\leq N-A,~ \abs{\hat{\ell}_i(X)-\ell_i'(X)}>\eps}\leq N\delta_N,\] and the result follows.
\end{proof}

\subsection{Lemmas for \cref{thm:BFNRoptimality}} \label[appendix]{sec:LemmasForBTDRoptimality}We may concretely define $\ell_i^\infty$ as the almost sure limit
\begin{equation}\label{eqn:def:ellinftyLimit} \ell_i^\infty(X)= \lim_{K\to \infty} \Pi_H(\theta_i=0 \mid X_{-K},\dots,X_K);\end{equation} this limit is well defined by a standard martingale convergence theorem. 

\begin{lemma}\label[lemma]{lem:liinftyFullSupport}
	In the setting of \cref{thm:BFDRcontrol}, assume that the distribution function of the variable $f_1(X_1)/f_0(X_1)$ is continuous and strictly increasing on $(0,\infty)$. Then the distribution function of $\ell_i^\infty(X)$ is continuous and strictly increasing on $[0,1]$.
\end{lemma}
Note that atomicity of $\ell_i(X)$ relates to that of $f_1(X_i)/f_0(X_i)$, rather than that of $X_i$ itself, since for example the distribution of $\ell_1$ is atomic when $N=1$ if $\Pi_H(f_1(X_1)/f_0(X_1)=c)>0$ for some constant $c$. It is therefore unsurprising that the key properties of the distribution of $\ell_i^\infty(X)$ depend on the distribution of the ratio $f_1(X_i)/f_0(X_i)$. \\

\begin{proof} 
	Let $G_0$ denote the distribution function of $(f_{1}/f_{0})(X_1)$ when $X_1 \sim f_0 \mu$ and $G_1$ the distribution function of $(f_{1}/f_{0})(X_1)$ when  $X_1 \sim f_1 \mu$. \\
	Define the stationary filter sequence  $( \Phi_i^\infty(X))_{i\in \ZZ}$ by
	\begin{equation}\label[equation]{eqn:def:filter} \Phi_i^\infty(X) := \Pi_H (\theta_i=0 \mid (X_n : n\in\ZZ, n\leq i)).
	\end{equation}
	Using the usual forward-backward equations, see \cite{BPSW70}, and taking almost-sure limits one obtains the following forward equation: for each $i$,
	$$
	\Phi_i^\infty(X) =\frac{[(1-p)\Phi_{i-1}^\infty(X)+q (1-\Phi_{i-1}^\infty(X))]f_{0}(X_{i})}{((1-p)f_{0}(X_{i})+ pf_{1}(X_{i}) )\Phi_{i-1}^\infty(X)+(qf_{0}(X_{i})+(1-q) f_{1}(X_{i})) (1-\Phi_{i-1}^\infty(X))}
	$$
	where $p=Q_{01}$ and $q=Q_{10}$, leading to
	\begin{equation}\label[equation]{eqn:forward} \Phi_i^\infty(X) =\frac{(1-p)\Phi_{i-1}^\infty(X)+q (1-\Phi_{i-1}^\infty(X))}{(1-p+ p(f_{1}/f_{0})(X_{i}) )\Phi_{i-1}^\infty(X)+(q+(1-q) (f_{1}/f_{0})(X_{i})) (1-\Phi_{i-1}^\infty(X))}.
	\end{equation}
	That is, if we define $A(\Phi)=(1-p)\Phi + q (1-\Phi)$, then
	\begin{equation}\label[equation]{eqn:forward2} \Phi_i^\infty(X) =\frac{A(\Phi_{i-1}^\infty(X))}{A(\Phi_{i-1}^\infty(X))+(f_{1}/f_{0})(X_{i}) (1-A(\Phi_{i-1}^\infty(X)))}.
	\end{equation}
	Since conditional on $\Phi_{i-1}^\infty(X)$, $X_{i}$ has distribution $\left[A(\Phi_{i-1}^\infty(X)) f_{0}(x) + (1-A(\Phi_{i-1}^\infty(X)) f_{1}(x)\right]\mu$, 	we deduce that $(\Phi_{i}^\infty(X))_{i\in\ZZ}$ is a stationary Markov chain with transition kernel $K(\Phi, d\Phi')$ given by
	\begin{eqnarray*}
		K(\Phi, d\Phi')&=&\int \delta_{g(\Phi, x)} (d\Phi') \left[(\Phi (1-p)+(1-\Phi)q) f_{0}(x) + (\Phi p +(1-\Phi)(1-q)) f_{1}(x)\right]d\mu (x)\\
		&=&\int \delta_{g(\Phi, x)} (d\Phi') \left[A(\Phi) f_{0}(x) + (1-A(\Phi)) f_{1}(x)\right]d\mu (x),
	\end{eqnarray*}
	where
	\[
	g(\Phi,x)=\frac{A(\Phi)}{A(\Phi)+(f_{1}/f_{0})(x)(1-A(\Phi))}.
	\]
	Then,  for each $t\in (0,1)$, we have
	\[
	\Pi_H \left(\Phi_{i}^\infty(X) \leq t \vert \Phi_{i-1}^\infty(X) \right)=\Pi_H \left((f_{1}/f_{0})(X_{i}) \geq  \frac{A(\Phi_{i-1}^\infty(X))}{1-A(\Phi_{i-1}^\infty(X))}(1/t -1)
	\mid \Phi_{i-1}^\infty(X) \right).
	\]
	Recall that  $\pi_{0}G_0 + \pi_{1}G_1$ is assumed to be continuous and strictly increasing on $(0,+\infty)$, and that $\pi_{0}>0$ and $\pi_{1}>0$, so that $G_{0}$ and $G_{1}$ are both continuous, and on the set where $G_0$ is not strictly increasing, $G_{1}$ is strictly increasing and vice versa.		
	We deduce that
	\begin{multline*}
		\Pi_H\left(\Phi_{i}^\infty(X) \leq t \mid \Phi_{i-1}^\infty(X) \right)=A(\Phi_{i-1}^\infty(X))
		\left[1-
		G_{0}\left( \frac{A(\Phi_{i-1}^\infty(X))}{1-A(\Phi_{i-1}^\infty(X))}\left(\frac{1}{t} -1\right)\right)\right]\\
		+\left(1-A(\Phi_{i-1}^\infty(X))\right)\left[1-
		G_{1}\left(\frac{A(\Phi_{i-1}^\infty(X))}{1-A(\Phi_{i-1}^\infty(X))}\left(\frac{1}{t} -1\right)\right)\right].
	\end{multline*}
	Then $\Phi_{i}^\infty(X)$ has, conditionally on $\Phi_{i-1}^\infty(X)$, a continuous and strictly increasing distribution function on $(0,1)$. 
	The same holds for $\Phi_{i}^\infty(X)$ since for all $t$, 
	$$
	\Pi_H\left(\Phi_{i}^\infty(X) \leq t  \right)=E_H[\Pi_H\left(\Phi_{i}^\infty(X) \leq t \mid \Phi_{i-1}^\infty(X) \right)].
	$$
	That is, $\Phi_{i}^\infty(X)$ has (conditionally on $\Phi_{i-1}^\infty(X)$ and unconditionally) no atoms and support $(0,1)$.
	
	The same ideas used to derive \cref{eqn:forward,eqn:forward2} allow us to show that for all $i$,
	\begin{equation}\label[equation]{eqn:backward2} \ell_i^\infty(X) =\frac{(1-p)\Phi_{i}^\infty(X)\ell_{i+1}^\infty(X)}{A(\Phi_{i}^\infty(X))}
		+\frac{p\Phi_{i}^\infty(X)(1-\ell_{i+1}^\infty(X))}{1-A(\Phi_{i}^\infty(X))}.
	\end{equation}
	Let  $C(\Phi)=\frac{\Phi (1-\Phi)}{A(\Phi)(1-A(\Phi))}$ and notice that for any $p,q\in (0,1)$ there exists $a=a(p,q)<1$ such that for all $\Phi \in (0,1)$, $|1-p-q | C(\Phi) \leq a$. Then an easy recursion yields
	$$
	\ell_{i}^\infty(X)=\frac{p\Phi_{i}^\infty(X)}{1-A(\Phi_{i}^\infty(X))} +\sum_{k\geq 1}(1-p-q)^{k}C(\Phi_{i}^\infty(X))C(\Phi_{i+1}^\infty(X))\cdots C(\Phi_{i+k-1}^\infty(X))\frac{p\Phi_{i+k}^\infty(X)}{1-A(\Phi_{i+k}^\infty(X))}.
	$$
	Indeed, since for any $\Phi \in (0,1)$, $|1-p-q | C(\Phi)\leq a(p,q) <1$, the series converges almost surely.
	We see that for each $i$, $\ell_{i}^\infty(X)$ is a function of $(\Phi_{k}^\infty(X))_{k\geq i}$, and we have
	\[
	\ell_{i}^\infty(X)  = \frac{p\Phi_{i}^\infty(X)}{1-A(\Phi_{i}^\infty(X))} + (1-p-q)C(\Phi_{i}^\infty(X))\ell_{i+1}^\infty(X).
	\]
	It follows that for all $t$, 
	\begin{equation}
		\label{eq:repfell}
		\begin{split}
			&\Pi_H \brackets{ \ell_{i}^\infty(X) \leq t \vert \Phi_{i-1}^\infty(X) } \\  =&E_H\sqbrackets[\Big]{ \Pi_H \brackets[\Big]{  (1-p-q)C(\Phi_{i}^\infty(X))\ell_{i+1}^\infty(X) \leq t -\frac{p\Phi_{i}^\infty(X)}{1-A(\Phi_{i}^\infty(X))}\mid \Phi_{i}^\infty(X) } \mid  \Phi_{i-1}^\infty(X) }.
		\end{split}
	\end{equation}
	Define the function $F_{\ell}$ by
	\[	F_{\ell}(t ; \Phi_{i-1}^\infty(X))=\Pi_H \left(\ell_{i}^\infty(X) \leq t \vert \Phi_{i-1}^\infty(X)\right);
	\]
	note that by stationarity $F_\ell$ does not depend on $i$. Then by (\ref{eq:repfell}), if $(1-p-q)> 0$, we have
	\[
	F_{\ell}(t ; \Phi_{i-1}^\infty(X))=E_H\sqbrackets[\bigg]{ F_{\ell}\brackets[\bigg]{\frac{1}{(1-p-q)C(\Phi_{i}^\infty(X))}\brackets[\Big]{  t -\frac{p\Phi_{i}^\infty(X)}{1-A(\Phi_{i}^\infty(X))}}; \Phi_{i}^\infty(X)}
		\mid \Phi_{i-1}^\infty(X) };
	\]
	that is, for any $t$ and any $\Phi \in (0,1)$,
	\begin{equation}
		\label{eq:Fsmoother1}
		F_{\ell}(t ; \Phi)=\int F_{\ell}\brackets*{\frac{1}{(1-p-q)C(x)}\left(t-\frac{px}{1-A(x)} \right); x}K(\Phi,dx).
	\end{equation}
	Similarly, if $(1-p-q)< 0$, defining the function $\tilde{F_{\ell}}$ by $\tilde{F_{\ell}}(t,\Phi)=\lim_{s\rightarrow t, s<t}F_{\ell}(t ; \Phi)$,
	\begin{equation}
		\label{eq:Fsmoother2}
		F_{\ell}(t ; \Phi)=\int \left[1-\tilde{F_{\ell}}\left(\frac{1}{(1-p-q)C(x)}\left(t-\frac{px}{1-A(x)} \right); x\right)\right]K(\Phi,dx).
	\end{equation}
	Note that under \cref{ass:Qpi}, $(1-p-q)\neq 0$.\\
	Finally, the fact that  $\Phi_{i}^\infty(X)$ has no atoms and support $(0,1)$  (both conditionally on $\Phi_{i-1}^\infty(X)$ and unconditionally) implies, together with equations \cref{eq:Fsmoother1} and  \cref{eq:Fsmoother2}, that whatever the sign of $(1-p-q)$, the function  $t\mapsto E_H[F_{\ell}(t ; \Phi_{i-1}^\infty(X))]$ is continuous and strictly increasing, which is to say that the distribution function of $\ell_{i}^\infty(X)$ is continuous and strictly increasing. 
\end{proof}

\begin{lemma}\label[lemma]{lem:ConditionalExpectationIncreasing}
	Under the conditions of \cref{thm:BFNRoptimality}, writing $\ell_i^\infty(X)=\Pi_H(\theta_i=0\mid (X_n)_{n\in\ZZ}),$
	the function $m$ defined by 
	\[ m(\lambda)= E[\ell_i^\infty(X) \mid \ell_i^\infty(X)<\lambda]\] 
	is continuous and strictly increasing on $(0,1)$, and $m(\lambda)<\lambda$ for all $\lambda\in(0,1)$.
	\end{lemma}
\begin{proof}
	For any random variable $U$ and any $a<b$ such that $P(U<a)>0$, we have
\[\begin{split} E[ U \mid U<b]
		&= E[ U \mid  U<a] P(U<a \mid  U<b) + E[ U \mid a\leq U<b] P(U\geq a \mid U<b) \\ 
		&=E[U\mid U<a](1-P(U\geq a \mid U<b)) + E[U\mid a\leq U<b] P(U\geq a\mid U<b),\end{split} \]
		hence
		\begin{equation} \label{eqn:ConditionalExpectationU} E[U\mid U<b]-E[U\mid U<a] = \frac{P(a\leq U<b)}{P(U<b)}\brackets[\Big]{E[U\mid a\leq U<b]-E[U\mid U<a]}. \end{equation}
		Note now that $E[U\mid U<a]< a$: indeed, if \[V\overset{d}{=} (U-a)\mid U<a,\] then $V\leq 0$ and $V$ is strictly negative with positive probability, hence $E[V]<0$. [For $U=\ell_i^\infty$ this yields that $m(\lambda)<\lambda$ as claimed.]
We similarly note that $E[U \mid a\leq U<b]\geq a$, so that, using also that $U$ is bounded, so that $E[U\mid U<a]\geq -c$ for some $c<\infty$,
\[0< E[U \mid a\leq U <b]- E[ U \mid U<a]<b+c.\]
Consequently, returning to \cref{eqn:ConditionalExpectationU}, to see that 
$E[U\mid U<x]$ is strictly increasing on $\braces{x : P(U<x)>0}$ it suffices to show that $P(a\leq U<b)>0$ for all $a,b$, and to show it is continuous it suffices to show that $P(a\leq U<b)\to 0$ as $b-a\to 0$. Taking $U=\ell_i^\infty$, we conclude by \cref{lem:liinftyFullSupport}, which tells us that the distribution function of $\ell_i^\infty$ is continuous and strictly increasing and also implies that $\Pi_H(\ell_i^\infty<\lambda)>0$ for all $\lambda>0$.
\end{proof}

\begin{lemma}\label[lemma]{lem:LHatApproximatedBylinfty}
Recall the definition
	\[
	\ell_i^\infty(X) = \Pi_H (\theta_i=0 \mid (X_n : n\in \ZZ)).	
	\] 
	There exist $\delta_N,\xi_N,\xi_N' \to 0$ such that with probability tending to 1,
	\begin{align*} &\#\braces{i : 1\leq i \leq N,~ \abs{\ell_i(X)-\ell_i^\infty(X)}>\xi_N} \leq N\delta_N \\ 
	 &\#\braces{i : 1\leq i \leq N,~ \abs{\hat{\ell}_i(X)-\ell_i^\infty(X)}>\xi_N'}\leq N\delta_N.\end{align*}
\end{lemma}

\begin{proof}
	Define $\ell_i'(X)=\Pi_H(\theta_i=0\mid X_{i-A_N},\dots,X_{i+A_N}).$ As in \cref{lem:LHatApproximatedByl'}, we may argue using \cite[Proposition 4.3.23iii]{CMR05} that for a suitable sequence $A_N\to \infty$ satisfying $A_N/N\to 0$, that 
	\[ \#\braces{i\leq N : \abs{\ell_i(X)-\ell_i'(X)}>4\rho^{A_N}}\leq 2A_N.\]
	Recalling from \cref{eqn:def:ellinftyLimit} that $\ell_i^\infty(X)$ is formally defined as an almost sure limit of $\ell_i'(X)$ as $A_N\to \infty$, so that $\ell_i'\to \ell_i^\infty$ in probability also, this proves the first bound. The second bound follows similarly after an appeal to \cref{lem:HatlVslErrors}.
\end{proof}
\begin{lemma}[Ergodic theorems]\label[lemma]{lem:ErgodicTheorems}
The sequences $\ell_i'$ and $\ell_i^\infty$, defined for $A\in \NN$ by \begin{align*} \ell_i'(X) &= \Pi_H(\theta_i=0 \mid X_{i-A},\dots, X_{i+A}), \quad A<i\leq N-A, \\
		\ell_i^\infty(X) &= \Pi_H (\theta_i=0 \mid (X_n : n\in \ZZ)),	
	\end{align*}  are ergodic, so that for any bounded function $g$,
	\[ \frac{1}{N} \sum_{i=1}^N g(\ell_i')\to E_\pi[g(\ell_1')], \quad \text{a.s.\ (hence also in probability)},\] and similarly for $\ell_i^\infty$.
\end{lemma}
\begin{proof}
	These are standard ergodicity results for functions of Markov chains; see for example \cite[Chapter 6]{Durrett19}. In the case of $\ell_i'$ one can also note that $g(\ell_i'(X))$ is a function of the Markov chain $(\theta_{i-A},\dots,\theta_{i+A},X_{i-A},\dots,X_{i+A})$ to reduce to the ergodic theorem for Markov chains themselves. 
\end{proof}

\begin{lemma}\label[lemma]{lem:OptimalityOfOracleClass}
In the setting of \cref{thm:BFDRcontrol}, define the class $(\vphi_{\lambda,H} : \lambda\in[0,1])$ as in \cref{eqn:def:PhiLambda}, and define the mTDR and mFDR as in \cref{eqn:def:mFDR,eqn:def:mTDR}.
	Then for each $\lambda\in(0,1)$ we have
		\begin{equation*}
			 \mTDR_H(\vphi_{\lambda,H}) = \sup\braces{\mTDR_H(\psi): \mFDR_H(\psi)\leq \mFDR_H(\vphi_{\lambda,H})}.\end{equation*}
	\end{lemma}
\begin{remarks}
	\begin{enumerate}[i.]
		\item 
A version of this result in the HMM setting originates in \cite{SC09}, but to avoid a monotonicity property needed therein we instead adapt the proof of \cite[Lemma 9.2]{RRV19} (see also the proof of \cite[Theorem 1]{CSW19}). The proof is valid for $\ell$--value procedures in any (correctly specified) model, not just the hidden Markov model \cref{eqn:def:model}. 
\item The result does not in general hold for $\lambda=0$, since $\mFDR_H(\psi)=0$ whenever $E_H[\ell_i(X)\psi_i(X)]=0$, so that if $\Pi_H(\ell_i(X)=0)>0$, the test $\psi$ defined by $\psi_i(X)=\II\braces{\ell_i(X)=0}$ has $\mFDR_H(\psi)=0$ and $\mTDR_H(\psi)=\Pi(\ell_i(X)=0 \mid \theta_i=1)=\Pi(\ell_i(X)=0)>0=\mTDR_H(\vphi_{0,H})$. 
\item  In general, $\braces{\mFDR_H(\vphi_{\lambda,H}) : \lambda\in [0,1]}$ is a proper subset of $[0,1]$, and consequently the class $\vphi_{\lambda,H}$ need not be optimal for every threshold.  In particular, the supremum of the set is generally strictly smaller than one, and -- especially in discrete data settings -- there may be jump discontinuities in the function $\lambda\mapsto \mFDR_H(\vphi_{\lambda,H})$. The first of these does not cause any issues, since $\mTDR_H(\vphi_{1,H})=1=\sup_\psi {\mTDR_H(\psi)}$, while \cref{lem:OptimalTestNearContinuity} overcomes the issues raised in the second case in the setting of \cref{thm:BFNRoptimality}.
\end{enumerate}
\end{remarks}

\begin{proof}
Fix $\lambda>0$; write $\vphi$ for $\vphi_{\lambda,H}$ and let $a=\mFDR_H(\vphi)$. Observe that for any multiple testing procedure $\psi$, $\mFDR_H(\psi)\leq a$ if and only if 
	\[ E \sum_{i\leq N} (\ell_i-a)\psi_i \leq 0,\] with equality in one implying equality in the other. It follows that if $\mFDR_H(\psi)\leq a$ then 
	\begin{equation}\label{eqn:mFDRdifference} E \sum_{i\leq N} (\ell_i-a)(\vphi-\psi) \geq 0.\end{equation} 
We note also that $a<\lambda$. Indeed, if $\vphi=0$ almost surely then this is true by definition (recall the convention that 0/0=0 in the definition \cref{eqn:def:mFDR} of the mFDR). Otherwise, there exists $k$ such that $\vphi_k=1$ (and hence $\ell_k<\lambda$) with positive probability; then $U=(\ell_k-\lambda)\vphi_k$ satisfies $U\leq 0$ and $\Pi(U<0)>0$, which together imply that $E[U]<0$ and hence \[ E\sum_{i\leq N} (\ell_i-\lambda)\vphi_i <0.\]

We next show that, for all $i$,
	\begin{equation} \label[equation]{eqn:FromMFDRtoMTDR}(\ell_i-a)(\vphi_i-\psi_i)\leq \frac{\lambda-a}{1-\lambda}(1-\ell_i)(\vphi_i-\psi_i).\end{equation}
	Indeed, if 	$\vphi_i=1$, then $\ell_i<\lambda$, so that
	\[ \ell_i-a < \frac{1-\ell_i}{1-\lambda} (\lambda-a),\] and multiplying by $\vphi_i-\psi_i\geq 0$ yields the inequality, while if $\vphi_i=0$, then $\ell_i\geq \lambda>a$, so that
	\[ \ell_i-a\geq  \frac{1-\ell_i}{1-\lambda}(\lambda-a),\] and multiplying by $\vphi_i-\psi_i\leq 0$ yields the inequality. 
	
	Now, since $a<\lambda<1$, so that also $(1-\lambda)/(\lambda-a)> 0$, we deduce from \cref{eqn:mFDRdifference,eqn:FromMFDRtoMTDR} that 
	\[ E \sum_{i\leq N} (1-\ell_i)(\vphi_i-\psi_i)\geq 0.\] 
	Finally, by definition,
	\[ \mTDR_H(\vphi)=\frac{E[\sum_{i\leq N} (1-\ell_i)\vphi_i]}{N\pi_1}, \quad \mTDR_H(\psi)=\frac{E[\sum_{i\leq N} (1-\ell_i)\psi_i]}{N\pi_1},\] hence
	$\mTDR_H(\vphi)\geq \mTDR_H(\psi)$ as claimed.
	\end{proof}

\begin{lemma}\label[lemma]{lem:OptimalTestNearContinuity}
	In the setting of \cref{thm:BFNRoptimality}, define the map
	\[ g: x\mapsto \sup\braces{\mTDR(\psi) : \mFDR(\psi)\leq x},\]
	where the supremum is defined over multiple testing procedures $\psi$. Then for sequences $x_N,y_N$ such that $\abs{x_N-y_N}\to 0$, we have
	\[ \abs{g(x_N)- g(y_N)}\to 0 \quad \text{as }N\to \infty. \]
	[Note that $g$ depends implicitly on $N$, so that this does not simply say that $g$ is continuous.]
\end{lemma}
\begin{proof}
Prompted by	\cref{lem:OptimalityOfOracleClass}, we focus on tests $\psi$ of the form $\vphi_{\lambda,H}$, $\lambda\in[0,1]$ and define, for $N\geq 1$, 
\begin{align*}
	\lambda_N &= \sup\braces{\lambda : \mFDR_H(\vphi_{\lambda,H})\leq x_N}, \\
	\mu_N &= \sup\braces{\lambda: \mFDR_H(\vphi_{\lambda,H})\leq y_N}.
\end{align*}
As with the postFDR (recall \cref{eqn:HatLambdaCharacterisesPostFDRSign}) one has the dichotomies, implied by the fact that $\mFDR_H(\vphi_{\lambda,H})$ is non-decreasing and left-continuous in $\lambda$,
\begin{align*} \mFDR_H(\vphi_{\lambda,H})\leq x_N &\iff \lambda\leq \lambda_N, \\
	\mFDR_H(\vphi_{\lambda,H})\leq y_N &\iff \lambda\leq \mu_N,
\end{align*} and we set $\tilde{\lambda}_N=\min(\lambda_N+1/N,1)$ and $\tilde{\mu}_N=\min(\mu_N+1/N,1)$.
Then \cref{lem:OptimalityOfOracleClass} (and the remarks thereafter, for the cases $\tilde{\lambda}_N=1$, $\tilde{\mu}_N=1$) implies that 
\begin{align*} \mTDR_H(\vphi_{\lambda_N,H})\leq g(x_N) &\leq \mTDR_H(\vphi_{\tilde{\lambda}_N,H}) \\
	\mTDR(\vphi_{\mu_N,H})\leq g(y_N)&\leq \mTDR(\vphi_{\tilde{\mu}_N,H}).
\end{align*}

We now show that if $\lambda'>\lambda$, then for $N$ large enough
\begin{equation} \label[equation]{eqn:NearStrictness} \mFDR_H(\vphi_{\lambda',H})>\mFDR_H(\vphi_{\lambda,H});\end{equation} 
it will follow that necessarily $\abs{\lambda_N-\mu_N}\to 0$, since otherwise we cannot have $\abs{x_N-y_N} \to 0$ (this is trivial when $y_N$ is a constant sequence, and follows for general $y_N$ by compactness). 
Writing $a=\mFDR_H(\vphi_{\lambda,H})$ we note that as in the proof of \cref{lem:OptimalityOfOracleClass} we have $a<\lambda$, and for any test $\psi$,  
\[ \mFDR_H(\psi)\leq a \iff E \sum_{i\leq N} (\ell_i-a)\psi_i \leq 0.\] 
Then
\[ E \sum_{i \leq N} (\ell_i-a)\II\braces{\ell_i<\lambda'} = E\sum_{i\leq N}(\ell_i-a)\II\braces{\ell_i<\lambda} + E\sum_{i\leq N}(\ell_i-a)\II\braces{\lambda\leq \ell_i<\lambda'}.\]
The first term on the right equals zero, and we show that the second is strictly positive, for large $N$. 
Indeed, by \cref{lem:LHatApproximatedBylinfty} there exists a sequence $\xi_N\to 0$ such that $E\#\braces{i: \abs{\ell_i-\ell_i^\infty}>\xi_N}/N \to 0$ as $N\to \infty$; then the term in question is lower bounded by $(\lambda-a)$ multiplied by
\[E \#\braces{i : \lambda+\xi_N \leq \ell_i^\infty<\lambda'-\xi_N}- E\#\braces{i : \abs{\ell_i-\ell_i^\infty}>\xi_N}.\] \Cref{lem:liinftyFullSupport} tells us that under the assumptions of \cref{thm:BFNRoptimality} the distribution function of $\ell_i^\infty$ is strictly increasing, so that for $N$ large enough that $\lambda+\xi_N<\lambda'-\xi_N$ the first term on the right of the latest display is of order $N$ and the second is of smaller order, so that indeed the difference is positive, proving \cref{eqn:NearStrictness}.

Finally we prove that, as a consequence of the fact that $\abs{\lambda_N-\mu_N}\to0$, we have $\abs{\mTDR_H(\vphi_{\lambda_N,H})-\mTDR_H(\vphi_{\mu_N,H})}\to 0$. Since also $\abs{\tilde{\lambda}_N-\lambda_N}\to 0$, $\abs{\tilde{\mu}_N-\mu_N}\to 0$, the same proof will imply that each of $\mTDR_H(\vphi_{\lambda_N,H}),$ $\mTDR_H(\vphi_{\tilde{\lambda}_N,H})$, $\mTDR_H(\vphi_{\mu_N,H})$ and $\mTDR_H(\vphi_{\tilde{\mu}_N,H})$ differ by at most $o(1)$, allowing us to conclude. 

Assume for notational convenience that $\lambda_N\geq \mu_N$. The denominator in the expressions defining each of the mTDR's is $E\#\braces{i : \theta_i=1}=N\pi_1$, and we see that
	\[\mTDR_H(\vphi_{\lambda_N,H})= \mTDR_H(\vphi_{\mu_N,H})+\frac{E\#\braces{i: \theta_i=1,\mu_N\leq \ell_i<\lambda_N}}{N\pi_1}.\]
	
As used above, by \cref{lem:LHatApproximatedBylinfty} there exists a sequence $\xi_N\to 0$ such that $E\#\braces{i: \abs{\ell_i-\ell_i^\infty}>\xi_N}/N \to 0$ as $N\to \infty$. 
	\Cref{lem:liinftyFullSupport} tells us that the distribution function of $\ell_i^\infty$ is continuous -- and hence uniformly continuous -- and we see that
	\[  \begin{split} &N^{-1} E\#\braces{i: \theta_i=1,\mu_N\leq \ell_i<\lambda_N} \\ \leq &\Pi_H(\mu_N-\xi_N \leq \ell_1^\infty<\lambda_N+\xi_N) +N^{-1} E\#\braces{i : \abs{\ell_i-\ell_i^\infty}>\xi_N} \to 0, \end{split}\] as $N\to \infty$, proving the claim. \qedhere 
\end{proof}

\section{Auxiliary Results for the upper bounds of \cref{sec:SupNormEstimation}}\label[appendix]{sec:lemmasforLinftyEstimation}

\subsection{Well-definedness of the Estimators}\label[appendix]{sec:WellDefinednessOfEstimators}
\begin{lemma}
	\label[lemma]{lem:ExistenceOfSuitableHl} In the setting of \cref{thm:LinftyEstimation}, there exist $(h_l)_{l\in \NN}$ (not depending on $H$) uniformly supremum-norm bounded such that $O^{L_0}= (E[h_l(X_1) \mid \theta_1 = j]_{l\leq L_0,j\leq J})\in \RR^{L_0\times J}$ satisfies
	\[ \sigma_J(O^{L_0})\geq C,\] uniformly in $L_0\geq \underline{L}$, for some $C,\underline{L}$ depending on the parameters $f_j,~j\leq J$.
\end{lemma}
\begin{proof}
	For $L>L'$, $\sigma_J(O^L)>\sigma_J(O^{L'})$ 
	because $O^{L'}$ is a submatrix of $O^L$, see e.g.\ \cite[Chapter 1, Theorem 4.4]{SS90}.
	So it suffices to show that $\sigma_J(O^{\underline{L}})>0$ for some $\underline{L}$.

	Choose a countable family of sets $\Aa=\braces{A_1,\dots}$ generating the Borel $\sigma$--algebra on $\RR$, for example $\Aa=\braces{(-\infty,q): q\in\QQ}$, and let $h_l=\II_{A_l}$. Suppose for a contradiction that $\sigma_J(O^L)=0$ for all $L\in \NN$, or, put another way, that the $J$ vectors $(\ip{h_l,f_j}_{l\leq L})\in \RR^{L}$, $j\leq J$ are linearly dependent for all $L\in \NN$, so that there exist $a_1^{L},\dots,a_J^{L}\in [-1,1]$ for which $\sum_j \abs{a_j^{L}}=1$ and $\sum_{j} a_j^{L} \ip{h_l,f_j}= 0$ for all $l\leq L$. 
	By Bolzano--Weierstrass, there is a sequence $L_n\to \infty$ such that for each $j\leq J$, $a_j^{L_n}$ converges to some $a_j^\infty$, and note that necessarily $(a_j^\infty)_{j\leq J}$ is not the zero vector.
	For each $l \in \NN$, we have that
	\[ \ip{h_l,\sum_{j\leq J} a_j^\infty f_j}=\sum_{j\leq J} a_j^\infty \ip{h_l,f_j}=\lim_{n\to \infty}  \sum_{j\leq J} a_j^{L_n} \ip{h_l,f_j}=0.\] Since $\braces{h_l : l\in \NN}$ generates the Borel $\sigma$--algebra, it follows that $\sum_j a_j^\infty f_j$ corresponds to the zero measure hence, since it is a continuous function, is the zero function, contradicting that the functions $f_j,~j\leq J$ are linearly independent.
\end{proof}

\begin{lemma}\label[lemma]{lem:EventA}
	Under the assumptions of \cref{thm:LinftyEstimation}, define $\hat{P}$ and $(\hat{M}^x,\hat{B}^x,~x\in\RR)$ as in \cref{alg:Linfty} for $L$ such that $2^L\asymp (N/\log N)^{1/(1+2s)}$. Then
	\begin{lemenum}
		\item The map $x\mapsto \hat{M}^x$ is 
		 continuous. 
		For any $\kappa>0$, there exists $c=c(\kappa,\Hh)$ such that the event \[ \Aa=\braces{\norm{\hat{P}-P}\leq cL_0 r_N,~\sup_{x\in\RR}\norm{\hat{M}^x-M^x}\leq cL_0^2 r_N}\] (is measurable and) has probability at least $1-N^{-\kappa}$ for $N$ large.
		\item On $\Aa$, for $N$ large enough $\hat{P}$ has rank $J$
		, and the matrices 	\begin{equation} \label[equation]{eqn:def:tildeB} \tilde{B}^x=(\hat{V}^\transpose P\hat{V})^{-1} \hat{V}^\transpose M^x \hat{V}, \quad x\in\RR,
		\end{equation}
		are well defined.
		\item On $\Aa$, for some $C>0$ depending on both the constant $c$ of $\Aa$ and on $\Hh$, we have for $N$ large enough 
\begin{align} \sup_{x\in \RR} \max (\norm{\hat{B}^x},\norm{\tilde{B}^x})&\leq CL_0^{1/2},
	\\ 
\label[equation]{eqn:supTildeB-HatB}\sup_{x\in \RR} \norm{\tilde{B}^x-\hat{B}^x}&\leq CL_0^2 r_N.\end{align}
\end{lemenum}
\end{lemma}
\begin{proof}
	\Cref{lem:ControlOfP} and \cref{lem:ControlOfSupMx} together tell us that for suitable $c=c(\kappa,\Hh)$,
	\[ \Pi(\norm{\hat{P}-P}\leq cL_0r_N,~ \sup_{x\in\QQ} \norm{\hat{M}^x-M^x}\leq cL_0^2r_N)\geq 1- N^{-\kappa}.\] [In fact a union bound yields this with $2N^{-\kappa}$ in place of $N^{-\kappa}$, but the factor 2 can be removed by initially considering some $\kappa'>\kappa$.] We prove the claimed 
	continuity of the map $x\mapsto \hat{M}^x$; it will follow that  \[ \braces{\sup_{x\in \QQ} \norm{\hat{M}^x-M^x}\leq cL_0^2 r_N }=\braces{\sup_{x\in \RR} \norm{\hat{M}^x-M^x}\leq cL_0^2 r_N },\] which implies measurability and the probability bound for $\Aa$. This continuity results from the assumed Lipschitz continuity of $K$. 
	Indeed, if $\Lambda$ is the Lipschitz constant for $K$, observe that if $\abs{x-y}<\delta$ then for any $n$
	\[\abs{K_L(x,X_{n+1})-K_L(y,X_{n+1})}\leq \sup_{t\in \RR} \abs{K_L(x,t)-K_L(y,t)}\leq \sup_{\abs{u-v}<2^L\delta} 2^L\abs{K(u)-K(v)}\leq 2^{2L}\Lambda \delta,\] hence, for some $C=C(\Hh)$,
	\[ \norm{\hat{M}^x-\hat{M}^y} \leq \frac{L_0}{N} \max_l \norm{h_l}_\infty \max_{n\leq N}\abs{K_L(x,X_{n+1})-K_L(y,X_{n+1})} \leq C  \frac{L_02^{2L}}{N}\abs{x-y}. \] 
	
	Next, in view of the assumption on $O$ made in the \namecref{alg:Linfty}, \cref{lem:PfullRank} implies that $\sigma_J(P)$ is bounded away from zero for large $N$ and consequently by \cref{lem:InversionsPossible}, on $\Aa$ and for $N$ large we have that $\hat{P}$ is of rank $J$ and that $\hat{V}^\transpose P \hat{V}$ is invertible (recall that $L_0^5r_N\to 0$ by assumption, so that the condition of \cref{lem:LemsFromdCGL} -- that $\norm{\hat{P}-P}<\sigma_J(P)/3$ --  holds eventually). Then \cref{lem:AlgorithmIntuition} tells us that $\tilde{B}^x$ is well defined for each $x\in \RR$ and can be expressed as $(QO^\transpose \hat{V})^{-1} D^x QO^\transpose \hat{V}$. It follows, using \cref{lem:kappaQObound} and eq.\ \cref{eqn:lambdamaxBounded}, that on $\Aa$, for a constant $c=c(\Hh)$ and any $x\in\RR$ we have 
	\begin{equation*} 	\norm{\tilde{B}^x}\leq \kappa(QO^\transpose \hat{V}) \max_j \abs{ K_L[f_j](x)}\leq cL_0^{1/2} \end{equation*}
	for $N$ large.
	
	Finally, \cref{lem:TildeBCloseToHatB} tells us that on $\Aa$, for any $x\in \RR$ and for $N$ large enough that $cL_0r_N<\sigma_J(P)/3$, \[\norm{\tilde{B}^x-\hat{B}^x}\leq 3.2 \sqbrackets[\Big]{\frac{\norm{\hat{M}^x-M^x}}{\sigma_J(P)}+\frac{\norm{M^x}\norm{\hat{P}-P}}{\sigma_J(P)^2}}, \quad \forall x\in\RR.\] Noting that $\norm{M^x}\leq cL_0$ for some $c=c(\Hh)$ by \cref{lem:PfullRank}, we deduce \cref{eqn:supTildeB-HatB}. The bound for $\norm{\hat{B}^x}$ then follows from the bound for $\norm{\tilde{B}^x}$ by the triangle inequality.
\end{proof}

\begin{lemma}\label[lemma]{lem:ExistenceOfHatAHatU} 
			Recall $\sep(B)$ denotes the eigen-separation of a matrix $B$, in that if $B$ has eigenvalues $\lambda_1,\dots, \lambda_J$ then $\sep(B)=\min_{j\neq j'} \abs{\lambda_{j}-\lambda_{j'}}$.  
			On the event $\Aa$ of \cref{lem:EventA}, let $2^L\asymp (N/\log N)^{1/(1+2s)}$ and define $B^{a,u}\equiv \tilde{B}^{a,u}$ as in \cref{alg:Linfty} for $V=\hat{V}$:
	\[ \tilde{B}^{a,u}=\sum a_i \tilde{B}^{u_i},\qquad \tilde{B}^x=(\hat{V}^\transpose P\hat{V})^{-1} \hat{V}^\transpose M^x \hat{V}.\] Let $\Dd_N$ be an increasing sequence of finite sets consisting of dyadic rationals whose union $\cup_N \Dd_N$ is dense in $\RR$. Define
	\[\DD_N =\braces{(a,u)\in \Dd_N^{J(J-1)/2}\times \Dd_N^{J(J-1)/2} : \sum_i \abs{a_i}\leq 1}\]
	Then there exists a constant $c$ depending only on $f_1,\dots,f_J$ and positive when they are all distinct such that, on $\Aa$, \begin{equation*}\max \braces{\sep(\tilde{B}^{a,u}) : (a,u)\in\DD_N} \geq c,\end{equation*} 
	for all $N$ large.
\end{lemma}
\begin{remark*}
	Recall, as remarked after \cref{alg:Linfty}, that proving this result for $V=\hat{V}$ implies it holds for any $V$ such that $B^x=(V^\transpose P V)^{-1} (V^\transpose M^x V)$ is well-defined.
\end{remark*}
\begin{proof}
	In view of \cref{lem:AlgorithmIntuition}, $\tilde{B}^{a,u}$, being a linear combination of simultaneously diagonalisable matrices, is diagonalisable for any $a,u$, with eigenvalues 
	\[(\sum_i a_{i} K_L[f_j](u_{i}))_{j\leq J}.\] 
	Recall that $\norm{K_L[f_j]-f_j}_\infty\to 0$ as $L=L(N)\to \infty$ by \cref{eqn:KLfapproximatesf}. It follows by the triangle inequality that
	\[\max_{\DD_N} \abs{\sum_i a_{i} \brackets[big]{K_L[f_j](u_{i})-f_j(u_{i})}}\to 0,\] hence \begin{equation}\label[equation]{eqn:supauPositive} \max_{\DD_N}\sep(\tilde{B}^{a,u})=\max_{\DD_N} \min_{j\neq j'}\abs[\Big]{\sum_i a_{i} K_L[f_j-f_{j'}](u_{i})}> \tfrac{1}{2}  \max_{\DD_N} \min_{j\neq j'} \abs[\Big]{\sum_i a_{i}\brackets[\big]{f_j(u_{i})- f_{j'}(u_{i}) }},\end{equation} 
	for $N$ large, provided this latter quantity is strictly positive. 
	
	Next, let $U_N$ be a sequence of sets, increasing to $\RR$, such that $\sup_{u\in U_N} \min_{d\in \Dd_N} \abs{u-d}\to 0.$ 	Observe that, since $f\in C^s(\RR)$,
	\[(a,u)\mapsto \min_{j\neq j'} \abs{\sum_i a_i (f_j(u_i)-f_{j'}(u_i))}\] is uniformly continuous on $\RR^{J(J-1)/2}\times\RR^{J(J-1)/2}$, so that 
	\begin{equation}\label[equation]{eqn:DNdense} \tfrac{1}{2}  \max_{(a,u)\in\DD_N} \min_{j\neq j'} \abs[\Big]{\sum_i a_{i} \brackets[\big]{f_j(u_{i})- f_{j'}(u_{i})}}>\tfrac{1}{4}\sup_{a}\sup_{u\in U_N}\min_{j\neq j'} \abs[\Big]{\sum_i a_{i} \brackets[\big]{f_j(u_{i})- f_{j'}(u_{i})}}\end{equation} 
	for $N$ large, provided this latter quantity is strictly positive. The supremum on the right can be extended: while at first we must take the supremum over ($a$ such that $\sum \abs{a_i}\leq 1$ and) $u\in U_N$, the result remains true taking the supremum instead over all $u\in \RR^{J(J-1)/2}$, using that $f_j(u)\to 0$ as $u\to \infty$.
	We now prove that
	\[ \sup_{a,u} \min_{j\neq j'} \abs[\Big]{\sum_i a_i\brackets[\big]{f_j(u_i)-f_{j'}(u_i)}}>0.\]

	Choose for each pair $j\neq j'$ some $x\in \RR$ such that $f_j(x)\neq f_{j'}(x)$, and collect these $x$ into the vector $u$. For each $j\neq j'$, the set $\braces{ v\in \RR^{J(J-1)/2} : \ip{v, f_j(u_i)-f_{j'}(u_i)}=0}$ is a proper subspace of $\RR^{J(J-1)/2}$, so the union over these $J(J-1)/2$ spaces is not equal to $\RR^{J(J-1)/2}$ (for example it has Lebesgue measure zero) 
	and we may choose $a$ in the complement of the union. Scale invariance means that moreover we may assume $a$ satisfies $\sum_i \abs{a_i}=1$.
	Then $\abs{\sum_i a_i (f_j(u_i)-f_{j'}(u_i))}>0$ for each $j\neq j'$, as required.
	
	Finally, combining also with \cref{eqn:supauPositive,eqn:DNdense} we deduce that
	\[\max (\sep(\tilde{B}^{a,u}) : (a,u)\in\DD_N) > \tfrac{1}{4} \sup_{a,u}\min_{j\neq j'} \abs[\Big]{\sum_i a_{i} (f_j(u_{i})- f_{j'}(u_{i}))}>0,\] concluding the proof.
\end{proof}

\begin{lemma}\label[lemma]{lem:EmpiricalEigengap}
	In the setting of \cref{thm:LinftyEstimation}, let $\Aa$ be the event of \cref{lem:EventA}.
	Define $\hat{B}^x=\hat{B}^{x,L_0,L}$ and $\hat{B}^{a,u}$ as in \cref{alg:Linfty} for $2^L\asymp (N/\log N)^{1/(1+2s)}$. 
	There exists a constant $c=c(\Hh)>0$ such that on the event $\Aa$ we have
	\begin{equation} \label[equation]{eqn:SepHatB}
		\sep(\hat{B}^{\hat{a},\hat{u}})>c,
	\end{equation}
	for $N$ large and, defining $\tilde{B}^{a,u}=\sum a_i\tilde{B}^{u_i}$ for $\tilde{B}^x$ as in \cref{eqn:def:tildeB}, we also have	\begin{equation}\label[equation]{eqn:SepTildeB} \sep(\tilde{B}^{\hat{a},\hat{u}})>c.
	\end{equation}
	Note that \cref{eqn:SepHatB} implies in particular that $\hat{B}^{\hat{a},\hat{u}}$ has $J$ distinct eigenvalues and so is diagonalisable.
\end{lemma}

\begin{proof}
	By \cref{lem:EventA}, on $\Aa$ the matrices $\hat{B}^x,\tilde{B}^x$ are well-defined and bounded up to a constant by $L_0^{1/2}$, and satisfy for some $C=C(\Hh)$ \[ \sup_x \norm{\tilde{B}^x-\hat{B}^x}\leq  CL_0^2 r_N, \quad \sup_x\max(\norm{\tilde{B}^x},\norm{\hat{B}^x})\leq CL_0^{1/2}.\] By the triangle inequality, we deduce that \[\norm{\hat{B}^{a,u}}\leq \sum\abs{a_i}\norm{\hat{B}^{u_i}}\leq \sup_x \norm{\hat{B}^x}\leq CL_0^{1/2},\] and similarly $\norm{\tilde{B}^{a,u}}\leq CL_0^{1/2}$. Let $(a_N,u_N)\in \argmax_{\DD_N}(\sep(\tilde{B}^{a,u}))$ and recall by assumption that \[\sep(\tilde{B}^{a_N,u_N})>c \quad \text{ uniformly in $N$ large enough, for some $c>0$.}\]
	[As noted in the remark after \cref{lem:ExistenceOfHatAHatU}, choosing $V=\hat{V}$ in \cref{alg:Linfty}, and hence replacing $B^x$ defined therein with $\tilde{B}^x$, is valid on $\Aa$.]
	We apply the Ostrowski--Elsner theorem (\cref{thm:Ostrowski--Elsner}) to $A=\hat{B}^{a,u}$, $B=\tilde{B}^{a,u}$  to see for a constant $C=C(\Hh)$ that for any $a,u$ we have
	\[ \min_\tau \max_j \abs{\lambda_{\tau(j)}(\tilde{B}^{a,u})-\lambda_{j}(\hat{B}^{a,u})}\leq C L_0^{(J-1)/(2J)}(L_0^2r_N)^{1/J}, \] where $\lambda_j, j\leq J$ are maps taking matrices to their eigenvalues. This last expression tends to zero as $N\to \infty$ (since by assumption $L_0^{(J+3)/2}r_N\to 0$) and in particular it is smaller than $\sep(\tilde{B}^{a_N,u_N})/5$ for $N$ large.
	
	By the triangle inequality we deduce that on $\Aa$,
	\[ \sep(\hat{B}^{a_N,u_N}) \geq \sep(\tilde{B}^{a_N,u_N})-2\sup_{a,u} \min_\tau \max_j \abs{\lambda_{\tau(j)}(\tilde{B}^{a,u})-\lambda_{j}(\hat{B}^{a,u})} \geq  (3/5)\sep(\tilde{B}^{a_N,u_N}).\]
	It follows by definition of $\hat{a},\hat{u}$ that \[\sep(\hat{B}^{\hat{a},\hat{u}}) \geq \sep(\hat{B}^{a_N,u_N})\geq  (3/5)\sep(\tilde{B}^{a_N,u_N}),\] proving \cref{eqn:SepHatB}.
	Applying the triangle inequality again we conclude that
	\[ \sep(\tilde{B}^{\hat{a},\hat{u}})\geq (1/5)\sep(\tilde{B}^{a_N,u_N}),\] proving \cref{eqn:SepTildeB}.
\end{proof}

\subsection{Concentration of Empirical Estimators}\label[appendix]{sec:ConcentrationOfEmpiricalEstimators}

We note the following concentration results for Markov chains, adapted as in \cite[Proposition 13]{dCGL16} from results of \cite{Paulin2015}, which will allow us to control the errors of the empirical estimators $\hat{P}$ and $\hat{M}^x$. The pseudo-spectral gap of a chain is defined in \cite{Paulin2015}, wherein it is noted that its reciprocal is equivalent to the mixing time. The bracketing numbers $N_{[]}(\Tt,\norm{\cdot}_{L^2(P)},\eps)$ are defined as the smallest number of pairs of functions $(\underline{f},\bar{f})$ such that every $g\in \Tt$ is bracketed by one of the pairs, where $(\underline{f},\bar{f})$ brackets $g$ if $\underline{f}\leq g\leq \bar{f}$ pointwise.
\begin{lemma}\label[lemma]{lem:ConcentrationResultsFromPaulin}
Let $Y$ be a stationary Markov chain taking values in $\Yy$ with pseudo-spectral gap $\gamma_{\textnormal{ps}}>0$, with law denoted $P$. Let $\Tt$ be some countable class of real valued and measurable functions on $\Yy$. Assume there exist $\sigma,b>0$ such that for all $t\in \Tt$, $\norm{t}_{L^2(P)}\leq \sigma$ and $\norm{t}_\infty \leq b$. Suppose that the $L^2(P)$ bracketing entropy 
\[H_{[]}(\Tt,\norm{\cdot}_{L^2(P)},\eps):=\log N_{[]}(\Tt,\norm{\cdot}_{L^2(P)},\eps),\] is upper bounded by some $\bar{H}(\eps)$, achievable using brackets of $L^\infty$--diameter at most $b$. 
Then for fixed $t\in \Tt$ we have 
\begin{equation}\label[equation]{eqn:PaulinConcentration} P(\abs{\sum (h(Y_i)-E h(Y_1))}\geq x)\leq 2 \exp\brackets[\Big]{- \frac{x^2 \gamma_{\textnormal{ps}}}{8 (N+1/\gamma_{\textnormal{ps}})\sigma^2+20bx}},\end{equation} and there exists $C>0$ depending only on a lower bound for $\gamma_{\textnormal{ps}}$ such that
\begin{equation}\label[equation]{eqn:PaulinSupremumConcentration}P\brackets[\Big]{\sup_{t\in \Tt} \sum_{n=1}^N (t(Y_n)- E t)\geq C\sqbrackets{A+\sigma\sqrt{Nx}+bx}}\leq \exp(-x),\end{equation}
where
\[ A=\sqrt{N}\int_0^\sigma \sqrt{\bar{H}(u)\wedge N}\du + (b+\sigma)\bar{H}(\sigma).\]
\end{lemma}
\begin{proof} 
For the first claim, see \cite[Theorem 3.4]{Paulin2015} (but note there is an updated version of the paper on arXiv). For the second, observe that the proof of the same theorem gives the following bound for the Laplace transform of 
$S= \sum (t(Y_n)-Et)/b$:
\begin{equation}\label[equation]{eqn:LaplaceTransform} E \exp(\lambda S) \leq \exp\brackets[\Big]{\frac{2(N+1/\gamma_{\textnormal{ps}})(\sigma^2/b^2)}{\gamma_{\textnormal{ps}}}\lambda^2 \brackets[\Big]{1-\frac{10\lambda}{\gamma_{\textnormal{ps}}}}^{-1}}.\end{equation}
One now appeals to \cite[Theorem 6.8]{Massart07} and the consequent Corollary 6.8. While the theorem is stated for independent random variables, the proof uses this condition only when applying Lemma 6.6 of the same reference, a version of which holds also in the current setting thanks to \cref{eqn:LaplaceTransform}.
\end{proof}

\begin{lemma}\label[lemma]{lem:ControlOfP}
In the setting of \cref{thm:LinftyEstimation} and defining $P,\hat{P}$ as in \cref{alg:Linfty}, for any $\kappa>0$ there exists $C=C(\kappa,\Hh)$ such that
\[ \Pi_H \brackets[\Big]{ \norm{\hat{P}-P} > CL_0 (N/\log N)^{-1/2}} \leq N^{-\kappa}\]
\end{lemma}
\begin{proof} 
Noting that $Y_n=(X_n,X_{n+1},X_{n+2},\theta_n,\theta_{n+1},\theta_{n+2})$ defines a stationary Markov chain, we apply \cref{eqn:PaulinConcentration} to deduce that
\[\begin{split} &\Pi\brackets[\Big]{\abs[\Big]{\frac{1}{N} \sum_{i=1}^{N} h_{ij}(Y_n)-E[h_{ij}]}>C \brackets[\Big]{\frac{\log N}{N}}^{1/2} } \\ \leq & 2 \exp\brackets[\Big]{ - \frac{C^2 \gamma_{\textnormal{ps}} N \log N}{8 (N+1/\gamma_{\textnormal{ps}})\Var_\pi (h_{ij})+20C (N\log N)^{1/2}\norm{h_{ij}}_\infty}},\end{split}\]
where $h_{ij}(Y_n)=h_i(Y_{n,1})h_j(Y_{n,3})$ and where $\gamma_{\textnormal{ps}}$ is the pseudo-spectral gap of the chain $Y_n$. 
We note that $\Var_\pi (h_{ij})\leq \norm{h_{ij}}_\infty^2\leq  \norm{h_{i}}_\infty^2\norm{h_j}_\infty^2$ is bounded by assumption. The pseudo spectral gap is also bounded: by \cite[Proposition 3.4]{Paulin2015} its reciprocal is controlled up to a constant by the mixing time of the Markov chain $Y_n$, which is equal to the mixing time of the chain $(\theta_n,\theta_{n+1},\theta_{n+2})_n$. This latter quantity is bounded since the assumption that $Q$ is irreducible and aperiodic on a finite state space implies that $\theta$ mixes exponentially, at a rate governed (again, in view of \cite[Proposition 3.4]{Paulin2015}) by the pseudo-spectral gap of $Q$ itself and $\min_j \pi_j$.
  
We deduce that for a constant $c=c(\Hh)$ we have
\[\Pi\brackets[\Big]{\abs[\Big]{\frac{1}{N} \sum h_{ij}(Y_n)-E[h_{ij}]}>C \brackets[\Big]{\frac{\log N}{N}}^{1/2} } \leq 2\exp(-C^2 c\log (N)).\]
For any $\kappa>0$, choosing $C=C(\kappa,c)$ large enough, this last probability is smaller than $N^{-\kappa}$ as claimed.
\end{proof}

\begin{lemma}\label[lemma]{lem:ControlOfSupMx}
	In the setting of \cref{thm:LinftyEstimation}, define $M^x=M^{x,L_0,L}$, $\hat{M}^x=\hat{M}^{x,L_0,L}$ as in \cref{alg:Linfty}, and recall that we choose $L$ such that $2^L\asymp (N/\log N)^{1/(1+2s)}$ and assumed that $L_0^5 r_N\to 0$. For any $\kappa>0$ there exists $C=C(\kappa,\Hh)$ such that
\[\Pi_H \brackets[\Big]{\sup_{x\in \QQ} \norm{\hat{M}^x-M^x} \geq C L_0^2 (N/\log N)^{-s/(1+2s)}} \leq N^{-\kappa}.\] 
\end{lemma}
\begin{proof} 
As in \cref{lem:ControlOfP} we note that the pseudo-spectral gap of the chain \[Y_n=(X_n,X_{n+1},X_{n+2},\theta_n,\theta_{n+1},\theta_{n+2})\] is bounded away from zero provided the same is true of $\min_j \pi_j$ and the pseudo-spectral gap of $Q$ itself, which holds by \cref{ass:Qpi'}.
We apply \cref{lem:ConcentrationResultsFromPaulin} to the family $\Tt=\braces{\pm h_i\otimes K_L(x,\cdot) \otimes h_j : i,j\leq L_0, x\in \QQ}$. Recall we assume that $\max(\norm{h_l}_\infty : l\leq L_0)$ is bounded independently of $L_0$. \Cref{lem:BracketingEntropyOfKernelTranslates} implies the bracketing entropy bound
\[ H_{[]}(\Tt,\norm{\cdot}_{L^2(\Pi_H)},\eps)\leq \bar{H}(\eps)=  C\log (L_0^{1/4} 2^L\eps^{-1}),\quad \eps\leq \sigma,\] where we may take
\begin{align*} b&=2^{L+2} \max_i(\norm{h_i}_\infty^2)\norm{K}_\infty = C2^L,\\
\sigma^2&\leq 2^L \max_i (\norm{h_i}_\infty^4) \norm{f_\pi}_\infty \int K(z)^2 \dz= C 2^{L/2};
\end{align*}
to bound $\sigma^2$ we have substituted $z=2^L(x-y)$ into $\int K_L(x,y)^2 f_\pi(y)\dy \leq \norm{f_\pi}_\infty \int 2^{2L} K(2^L(x-y))^2\dy$.

An application of Jensen's inequality yields the standard bound 
\begin{equation}\label[equation]{eqn:IntSqrtLogBound} \int_0^x \sqrt{\log(1/u)}\du \leq x \sqrt{1+\log(1/x)}\leq x\brackets[\Big]{1+\sqrt{\log(1/x)}}.\end{equation} Performing suitable substitutions we deduce that
\[ \int_0^\sigma \sqrt{\log (L_0^{1/4}2^L/u)}\du = L_0^{1/4} 2^L \int_0^{\sigma /(2^{L}L_0^{1/4})} \sqrt{\log (1/v)}{\dif v} \leq \sigma \brackets[\Big]{1+\sqrt{\log (L_0^{1/4}2^L/\sigma)}}\leq C \sqrt{L2^L}, \] for some constant $C$, since by assumption $L_0^5r_N\to 0$, which implies that $\log(L_0)\leq \log N\asymp L$.
Noting that $(b+\sigma)\bar{H}(\sigma)\leq  CL2^L$ for some $C$, we deduce that 
\[ \Pi_H \brackets[\Big]{\sup_{t\in \Tt} \sum_{n=1}^N  \brackets[\big]{t(Y_n) - E_{H} t} \geq C\sqbrackets{ \sqrt{N2^L}(\sqrt{L}+\sqrt{\kappa \log N})+2^L(L+\kappa \log N)}} \leq \exp(-\kappa \log N).\]
Since $2^L\asymp (N/\log N)^{1/(1+2s)}$ we find, bounding the operator norm by the $L_0^2$ times the maximum of the entries, that as claimed,
\begin{equation}\label[equation]{eqn:supMinQ} \Pi_H \brackets[\Big]{ \sup_{x\in \QQ} \norm{\hat{M}^x-M^x} \geq C L_0^2 (N/\log N)^{-s/(1+2s)}} \leq N^{-\kappa}.\qedhere \end{equation} 
\end{proof}

\begin{lemma}\label[lemma]{lem:BracketingEntropyOfKernelTranslates}
Let $\Tt=\braces{h_i\otimes K_L(t,\cdot) \otimes h_j  : i,j\leq L_0, t\in \RR}$. Then we have the following bound for the bracketing numbers:
\begin{equation}
N_{[]}(\Tt,\norm{\cdot}_{L^2(\Pi_H)},\eps) \leq C L_0^2 \max( 2^{8L}\eps^{-8},1).
\end{equation}
for some constant $C>0$. This bound is achieved with brackets whose $L^\infty$--diameter is at most $2^{L+2}\norm{K}_\infty \max_i\norm{h_i}_\infty^2$.
\end{lemma}
\begin{proof} 
The kernel $K$ is assumed to be bounded, continuous, and supported in $[-1,1]$, see before \cref{eqn:def:KL}; note that necessarily $K(1)=K(-1)=0$.
Let $\Uu=\braces{\II_{(-\infty,u)} : u\in \RR}$, let $\Vv=\braces{\II_{(a,b]}-\II_{(c,d]} : a,b,c,d \in \RR}$. We show that
\begin{equation}\label[equation]{eqn:bracketingbounds} L_0^{-2} N_{[]}(\Tt,\norm{\cdot}_{L^2(\Pi_H)},\eps \norm{K_L}_\infty)\leq  N_{[]}(\Vv,\norm{\cdot}_{L^2(\Pi_H)},\eps) \leq N_{[]}(\Uu,\norm{\cdot}_{L^2(\Pi_H)},\eps/4)^4\end{equation} The first inequality follows from the fact that, given brackets $\sqbrackets{\underline{v}_k,\overline{v}_k}, k\leq N_\Vv$ of $L^2(\Pi_H)$--diameter $\eps$ for $\Vv$, we can define \[\underline{t}_{ikj} =\norm{K_L}_\infty h_i\otimes \underline{v}_k \otimes h_j, \quad \overline{t}_{ikj}=\norm{K_L}_\infty h_i\otimes \overline{v}_k \otimes h_j\] to obtain brackets $\sqbrackets{\underline{t}_{ikj},\overline{t}_{ikj}},i,j\leq L_0,k\leq N_\Vv$ for $\Tt$ whose $L^2(\Pi_H)$--diameter is $\norm{K_L}_\infty \eps$. For the second inequality, observe that any $v\in \Vv$ can be written in the form $(u_1-u_2)-(u_3-u_4)$ for $u_1,u_2,u_3,u_4\in \Uu$. Then, given brackets $\sqbrackets{\underline{u}_k,\overline{u}_k},k\leq N_\Uu$ for $\Uu$, it follows that $\sqbrackets{\underline{v}_{ijkl},\overline{v}_{ijkl}},i,j,k,l\leq N_\Uu$ form brackets for $\Vv$, where \[\underline{v}_{ijkl}=(\underline{u}_i-\overline{u}_j)-(\overline{u}_k-\underline{u}_l), \quad \overline{v}_{ijkl}=(\overline{u}_i-\underline{u}_j)-(\underline{u}_k-\overline{u}_l),\] and the $L^2(\Pi_H)$--diameter of such a bracket, if $\norm{\overline{u}_k-\underline{u}_k}_{L^2(\Pi_H)}=(E[\overline{u}_k-\underline{u}_k]^2)^{1/2}\leq \eps/4$ for each $k$, is by Cauchy--Schwarz at most \[ \brackets[\big]{E_{\Pi_H}[(v_{i_1,i_2,i_3,i_4}-v_{i_1,i_2,i_3,i_4})^2]}^{1/2}=\brackets[\big]{E_{\Pi_H} \sqbrackets[\big]{  \sum_{j} \overline{u}_{i_j}-\underline{u}_{i_j}}^2}^{1/2}\leq \brackets[\big]{4\sum_{j} E\sqbrackets[\big]{\overline{u}_{i_j}-\underline{u}_{i_j}}^2}^{1/2}\leq \eps. \]
It remains to bound $N_{[]}(\Uu,\norm{\cdot}_{L^2(\Pi_{H})},\eps)$. One argues as in the proof of the Glivenko--Cantelli theorem: let $R= \ceil{\eps^{-2}}$, set $x_0=-\infty$, $x_R=\infty$ and choose $x_k$ such that $\Pi_H(X_1\in [x_{k-1},x_k))=R^{-1}\leq \eps^2$. [This is possible because the distribution of $X_1$ has a density so is non-atomic, but the proof would require only minor adjustments to accomodate distributions with atoms.] Define 
\[ \underline{u}_k = \II_{(-\infty,x_{k-1})}, \quad \overline{u}_k=\II_{(-\infty,x_k)}, \quad 1\leq k\leq R,\] and note any $u\in \Uu$ is contained in one of the brackets $\sqbrackets{\underline{u}_k,\overline{u}_k}$. The $L^2(\Pi_H)$--diameter of such a bracket is at most
\[ \brackets[\big]{\Pi_{H}\braces[\big]{X_1\in [x_{k-1},x_k)}}^{1/2}= R^{-1/2}\leq \eps.\] It follows that
\[ N_{[]}(\Uu,\norm{\cdot}_{L^2(\Pi_H)},\eps)\leq \ceil{\eps^{-2}}.\]
The bracketing bound \cref{eqn:bracketingbounds} follows for a suitable constant $C$ by substituting into \cref{eqn:bracketingbounds} upon noting that $\ceil{\eps^{-2}}\leq 2\max(\eps^{-2},1)$ and that $\norm{K_L}_\infty=2^L \norm{K}_\infty$.

Finally, we note that $\overline{u}_i(x)-\underline{u}_i(x)\in[0,1]$ for $x\in \RR$ hence $\overline{v}_{ijkl}-\underline{v}_{ijkl}\leq 4$ (in fact one can improve the constant 4 to 1 by noting that in the expression $v=(u_1-u_2)-(u_3-u_4)$ we may assume $u_2\leq u_1$, $u_3\leq u_4$, and carefully considering the consequences). The brackets for $\Tt$ have $L^\infty$--diameter at most $2^{L+2}\norm{K}\max_i\norm{h_i}_\infty^2$ as a consequence.
\end{proof}

\subsection{Matrix Approximation Theory Arguments}\label[appendix]{sec:MatrixApproximationTheory}

\begin{lemma}\label[lemma]{lem:RhatEstimatesR}
	Define $\Aa$ as in \cref{lem:EventA}. In the setting of \cref{thm:LinftyEstimation}, define $\hat{R}$ as in \cref{alg:Linfty} for $2^L\asymp (N/\log N)^{1/(1+2s)}$, and define $\tilde{R}$ to have columns equal to the normalised columns of $QO^\transpose \hat{V}$.  Then, on $\Aa$, $\hat{R}$ is well-defined and
	\[ \norm{\hat{R}-\tilde{R}_\tau}\leq \norm{\hat{R}-\tilde{R}_\tau}_F\leq CL_0^{7/2}  r_N, \] for some $C=C(\Hh)$ and some permutation $\tau$, where $\tilde{R}_\tau$ is obtained by permuting the columns of $\tilde{R}$ according to $\tau$.
\end{lemma}
\begin{proof}
	\Cref{lem:EmpiricalEigengap} tells us on $\Aa$ that $\hat{B}^{\hat{a},\hat{u}}$ is diagonalisable, so that $\hat{R}$ is well defined, and moreover that \[\min\brackets[\big]{\sep(\hat{B}^{\hat{a},\hat{u}}),\sep(\tilde{B}^{\hat{a},\hat{u}})}>c,\] for some constant $c=c(\Hh)>0$.
	Now we apply \cite[Lemma C.3]{AHK12}, which says, as a consequence of the Bauer--Fike theorem, that if \[\eps=\kappa(\tilde{R})\sep(\tilde{B}^{\hat{a},\hat{u}})^{-1}\norm{\hat{B}^{\hat{a},\hat{u}}-\tilde{B}^{\hat{a},\hat{u}}}\] is smaller than $1/2$, then
	\[\norm{\hat{R}-\tilde{R}}\leq \norm{\hat{R}-\tilde{R}}_F\leq 4J^{1/2}(J-1)\norm{\tilde{R}^{-1}}\eps.\] By construction $\sum \abs{\hat{a}_i}\leq 1$, hence by the triangle inequality and \cref{lem:EventA}, on $\Aa$ we have 
	\[\norm{\hat{B}^{\hat{a},\hat{u}}-\tilde{B}^{\hat{a},\hat{u}}}\leq \sum_{i} \abs{\hat{a}_i} \norm{\hat{B}^{\hat{u}_i}-\tilde{B}^{\hat{u}_i}}\leq \sup_x \norm{\hat{B}^x-\tilde{B}^x}\leq CL_0^2 r_N,\] for some $C=C(\Hh)$. By \cref{lem:kappaQObound}, we have $\kappa(\tilde{R})\leq CL_0$ and $\norm{\tilde{R}^{-1}}\leq CL_0^{1/2}$. We deduce that $\eps\to 0$ on $\Aa$, hence is smaller than 1/2 for large $N$, and the result follows.
\end{proof}

One could directly use the Ostrowski--Elsner theorem (\cref{thm:Ostrowski--Elsner}) to obtain a version of \cref{thm:LinftyEstimation} with a suboptimal estimation rate. We here go through the slightly circuitous route of using \cref{thm:Ostrowski--Elsner} to prove an eigen-separation condition (i.e.\ \cref{lem:EmpiricalEigengap}) and deducing \cref{lem:RhatEstimatesR} 
because we may then apply the following \namecref{lem:C.4AHK12}, adapted from \cite[Lemma C.4]{AHK12},  to obtain a near-minimax rate instead.

\begin{lemma} \label[lemma]{lem:C.4AHK12} 
	Suppose $(A_t : t\in \Tt)$ are $J\times J$ matrices simultaneously diagonalised by a matrix $R$ with unit norm columns:
	\[R^{-1}A_t R=\diag(\lambda_{t,1},\dots,\lambda_{t,J}),~t\in\Tt.\]
	Let $\hat{R}$ be a matrix such that for some permutation $\tau$ of $\braces{1,\dots, J}$ we have \[\norm{\hat{R}-R_\tau}:= \eps_R\leq (1/2)\norm{R^{-1}}^{-1},\] where $R_\tau$ has is obtained by permuting the columns of $R$ according to $\tau$. Assume \[\lambda_{\max}:=\sup_t \max _{j}\abs{\lambda_{t,j}}<\infty.\] For matrices $(\hat{A}_t : t\in \Tt)$, write \[\eps_A:=\sup_t \norm{A_t-\hat{A}_t},\]
	and define \[\hat{\lambda}_{t,j}=e_j^\transpose \hat{R}^{-1}\hat{A}_t \hat{R} e_j.\] Then 
	\[ \sup_t \max_j \abs{\hat{\lambda}_{t,j}-\lambda_{t,\tau(j)}} \leq 4\kappa(R)[\eps_A+\lambda_{\max}\norm{R^{-1}}\eps_R].\]
\end{lemma}

\begin{proof}	
	Let $\hat{\zeta}_j^\transpose$ be the $j$th row of $\hat{R}^{-1}$, let $\hat{\xi}_j$ be the $j$th column of $\hat{R}$, and define $\zeta_j,\xi_j$ correspondingly with respect to the matrix $R_\tau$ obtained by permuting the columns of $R$ according to $\tau$. Then $\lambda_{t,\tau(j)}=\zeta_j^\transpose A_t\xi_j$, $\hat{\lambda}_{t,j}=\hat{\zeta}_j^\transpose \hat{A}_t\hat{\xi}_j$, and we have
	\begin{equation*}
		\begin{split}
			\abs{\hat{\lambda}_{t,j}-\lambda_{t,\tau(j)}}&=\abs{\hat{\zeta}_j^\transpose\hat{A}_t\hat{\xi}_j-\zeta_j^\transpose A_t\xi_j }\\
			&=\abs{\hat{\zeta}_j^\transpose \hat{A}_t (\hat{\xi}_j-\xi_j) + \hat{\zeta}_j^\transpose (\hat{A}_t-A_t)\xi_j + (\hat{\zeta}_j^\transpose-\zeta_j^\transpose)A_t\xi_j} \\
			&\leq\norm{\hat{\zeta}_j^\transpose}\norm{\hat{A}_t}\norm{\hat{\xi}_j-\xi_j}+\norm{\hat{\zeta}_j^\transpose}\norm{\xi_j}\eps_A+\norm{A_t \xi_j}\norm{\hat{\zeta}_j-\zeta_j}
		\end{split}
	\end{equation*}
	Using \cref{lem:BasicApproximationTheoryResults}, we have that
	\[\norm{\hat{R}^{-1}-R_\tau^{-1}}\leq \norm{R^{-1}}^2 \eps_R /(1-\norm{R^{-1}}\eps_R),\] and we further note the following:
	\begin{itemize}
		\item $\norm{\zeta_j^\transpose}=\norm{e_{\tau(j)}^\transpose R^{-1}}\leq \norm{R^{-1}}$, and $\norm{\hat{\zeta}_j^\transpose-\zeta_j^\transpose}\leq \norm{\hat{R}^{-1}-R_\tau^{-1}}\leq \norm{R^{-1}}^2\eps_R/(1-\norm{R^{-1}}\eps_R)$, so that also $\norm{\hat{\zeta}_j^\transpose}\leq \norm{\zeta_{\tau(j)}^\transpose}+\norm{\hat{\zeta}_j-\zeta_{\tau(j)}}\leq \norm{R^{-1}}/(1-\norm{R^{-1}}\eps_R)$ .
		\item $\norm{\xi_j}\leq \norm{R}$, and $\norm{\hat{\xi}_j-\xi_j}\leq \norm{\hat{R}-R_\tau}= \eps_R$.
		\item  $\norm{A_i}=\norm{R\diag(\lambda_{i,\cdot})R^{-1}}\leq \kappa(R)\lambdamax$, and $\norm{\hat{A}_t}\leq \norm{A_t}+\eps_A\leq \kappa(R)\lambda_{\max}+\eps_A$.
		\item $\norm{A_t\xi_j}=\abs{\lambda_{t,\tau(j)}}\norm{\xi_j}\leq \lambda_{\max}\norm{R}$. 
	\end{itemize}	
	Then, continuing the inequalities from the display, we have
	\begin{equation*}
		\begin{split}
			\abs{\hat{\lambda}_{t,j}-\lambda_{t,\tau(j)}}& \leq \frac{\norm{R^{-1}}}{1-\norm{R^{-1}}\eps_R} \sqbrackets[\Big]{ (\kappa(R)\lambdamax +\eps_A) \eps_R + \norm{R}\eps_A} + \lambda_{\max}\norm{R}\norm{R^{-1}}^2 \frac{\eps_R}{1-\norm{R^{-1}}\eps_R} \\
			&\leq \frac{\kappa(R)+\norm{R^{-1}}\eps_R}{1-\norm{R^{-1}}\eps_R}\eps_A+2\lambdamax\kappa(R) \frac{\norm{R^{-1}}\eps_R}{1-\norm{R^{-1}}\eps_R} \\
			&\leq (1+2\kappa(R))\eps_A+4\lambda_{\max}\norm{R^{-1}}\kappa(R)\eps_R,
		\end{split}
	\end{equation*}
	where for the last line we have used that $\norm{R^{-1}}\eps_R\leq 1/2$ by assumption.
	Taking the supremum over $t\in \Tt$ concludes the result since necessarily $1+2\kappa(R)\leq 3\kappa(R)<4\kappa(R)$.
\end{proof}

\begin{lemma}\label[lemma]{lem:PfullRank}
Define $O=O^{L_0},P=P^{L_0},(M^x=M^{x,L,L_0}: x\in \RR)$ as in \cref{lem:AlgorithmIntuition} for functions $(h_l)_{l\leq L_0}$ satisfying a sup-norm bound uniformly in $L_0$ and assume that $\sigma_J(O)\geq c>0$ uniformly in $L_0\geq \underline{L}$ for some $\underline{L}=\underline{L}(\Hh)$ (for example, by choosing $(h_l: l\leq L_0)$ as in \cref{lem:ExistenceOfSuitableHl}). Then \[\kappa(O)\leq C L_0^{1/2}, \quad  \sigma_J(P)\geq c', \quad \text{and}\quad \norm{M^x}\leq C'L_0,\] for some constants $c',C,C'>0$, uniformly in $L_0\geq \underline{L}$ and all $L$.
\end{lemma}
\begin{proof}
	Given the assumed bound on $\sigma_J(O)$, to control $\kappa(O)$ it remains to bound $\norm{O}$, since one has the standard expression $\kappa(O):=\norm{O}\norm{O^{-1}}\equiv \norm{O}/\sigma_J(O).$ 
	Then it suffices to note, using Cauchy--Schwarz and the fact that $\abs{\ip{f_j,h_l}}=\abs{\int h_l(x)f_j(x)\dx}\leq \norm{h_l}_\infty$, that
	\begin{equation}\label[equation]{eqn:NormOBound} \norm{O}^2=\sup_{\norm{v}=1} \sum_j (\sum_l v_l \ip{f_j,h_l})^2\leq  \max_l \norm{h_l}_\infty^2 JL_0.\end{equation}
	
	Next, \cref{ass:Qpi'} implies $\sigma_J(Q)>0$ and $\sigma_J(\diag(\pi))=\min_{j}\pi_j>0$. Using submultiplicativity of $\sigma_J$ (see \cref{lem:BasicApproximationTheoryResults}) and the expression $P=O\diag(\pi)Q^2O^\transpose$ (from \cref{lem:AlgorithmIntuition}), we have \[\sigma_J(P)=\sigma_J(O\diag(\pi)Q^2O^\transpose)\geq \sigma_J(O)\sigma_J(\diag(\pi))\sigma_J(Q)^2\sigma_J(O^\transpose)\geq c'(H)>0.\]
For $M^x$, the expression $M^x=O\diag(\pi)QD^xQO^\transpose$  from \cref{lem:AlgorithmIntuition} similarly yields
\[ \norm{M^x}\leq \norm{O}^2\norm{Q}^2 \max_j \abs{K_L[f_j](x)}.\] Recalling that $\norm{K_L[f_j]}_\infty$ is bounded (see \cref{eqn:lambdamaxBounded}) we deduce the result.
\end{proof}

The following collects several useful results from \cite{dCGlC17} and \cite{AHK12}.
\begin{lemma}\label[lemma]{lem:LemsFromdCGL} Assume $\sigma_J(O)\geq c>0$ uniformly in $L_0\geq \underline{L}$, so that by \cref{lem:PfullRank} we also have $\sigma_J(P)>0$ and $\kappa(O)\leq CL_0^{1/2}$ for some $C$. On the event $\Bb=\braces{\norm{\hat{P}-P}<\sigma_J(P)/3}$, for $L_0\geq \underline{L}$ and $N$ large enough we have the following.
	\begin{lemenum} 
		\item \label[lemma]{lem:InversionsPossible} $\sigma_J(\hat{P})>c/2$. Writing $\hat{V}$ and $V$ for matrices of orthonormal right singular vectors of $\hat{P}$ and $P$ respectively we have $\sigma_J(\hat{V}^\transpose V)^2\geq 3/4$, and consequently $\hat{V}^\transpose P \hat{V}$ is invertible. 
		\item \label[lemma]{lem:kappaQObound} $\kappa (QO^\transpose \hat{V})\leq C L_0^{1/2}$, $\norm{\tilde{R}^{-1}}\leq C'L_0^{1/2}$ and $\kappa(\tilde{R})\leq C''  L_0$, where $\tilde{R}$ is the matrix whose columns are those of $QO^\transpose \hat{V}$ but rescaled to have unit norm. 
	\item \label[lemma]{lem:TildeBCloseToHatB} For any $x\in \RR$, \[\norm{\tilde{B}^x-\hat{B}^x}\leq 3.2 \sqbrackets[\Big]{\frac{\norm{\hat{M}^x-M^x}}{\sigma_J(P)}+\frac{\norm{M^x}\norm{\hat{P}-P}}{\sigma_J(P)^2}}.\]
\end{lemenum}
\end{lemma}
\begin{proof}
	We throughout use various basic properties of $\sigma_J,\kappa$, which are summarised in \cref{lem:BasicApproximationTheoryResults} below.
	\begin{enumerate}[a.]
			\item  By \cref{lem:PfullRank}, $\sigma_J(P)>0$. The result then follows from standard approximation theory. In particular  \cite[Lemma C.1, part 2]{AHK12} tells us that $\sigma_J(\hat{P})>\sigma_J(P)/3>0$. That $\sigma_J(V^\transpose \hat{V})^2\geq 3/4$ on $\Bb$ is given by \cite[Lemma C.1, part 3]{AHK12} and submultiplicativity of $\sigma_J$ yields \[\sigma_J(\hat{V}^\transpose P \hat{V})=\sigma_J(\hat{V}^\transpose(V V^\transpose) P (V V^\transpose) \hat{V})\geq \sigma_J(V^\transpose \hat{V})^2 \sigma_J(V^\transpose PV)\geq (3/4)\sigma_J(P)>0,\] which implies invertibility of $\hat{V}^\transpose P \hat{V}$.
	\item  Observe that
	\[ \kappa (QO^\transpose \hat{V})= \frac{\norm{QO^\transpose\hat{V}}}{\sigma_J(QO^\transpose \hat{V})}\leq\frac{ \norm{QO^\transpose}}{\sigma_J(QO^\transpose V)\sigma_J(V^\transpose \hat{V})}.\]
		We have $\sigma_J(QO^\transpose V)=\sigma_J(QO^\transpose)$ and we deduce that $\kappa(QO^\transpose \hat{V})\leq ({4/3})^{1/2}\kappa(QO^\transpose)\leq 2 \kappa(Q)\kappa(O)$ by part a.
		\Cref{ass:Qpi'} implies $\kappa(Q)<\infty$.
		For $R$, see \cite[Lemma C.5]{AHK12}, which tells us that $\norm{\tilde{R}^{-1}}\leq \kappa(QO^\transpose\hat{V})$ and $\kappa(\tilde{R})\leq\kappa(QO^\transpose\hat{V})^2$.  
		\item One decomposes \[ \norm{\tilde{B}^x-\hat{B}^x}\leq \norm{(\hat{V}^\transpose \hat{P} \hat{V})^{-1}}\norm{\hat{V}^\transpose (M^x-\hat{M}^x)\hat{V}}+\norm{\hat{V}^\transpose M^x \hat{V}}\norm{(\hat{V}^\transpose P \hat{V})^{-1}-(\hat{V}^\transpose P \hat{V})^{-1}},\] then uses \cref{lem:BasicApproximationTheoryResults} with $\hat{A}=\hat{V}^\transpose \hat{B}^x\hat{V}$, $A=\hat{V}^\transpose \tilde{B}^x\hat{V}$, noting that in part a we showed $\norm{(\hat{V}^\transpose P \hat{V})^{-1}}\equiv \sigma_J(\hat{V}^\transpose P\hat{V})^{-1}\leq (4/3)\sigma_J(P)^{-1}.$  See the proof of \cite[Lemma F.4, on p28]{dCGlC17}, which adapts to the current setting.
	\end{enumerate}
\end{proof}

\begin{theorem}[{Ostrowski--Elsner, e.g.\ \cite[Chapter IV, Theorem 1.4]{SS90}}]
	\label[theorem]{thm:Ostrowski--Elsner}
		For a matrix $U\in\RR^{J\times J}$, write $(\lambda_i(U): i\leq J)$ for the eigenvalues of $U$. Then for matrices $A,B\in \RR^{J\times J}$ we have
	\begin{equation}\label[equation]{eqn:OstrowskiElsner}
		\min_\tau \max_j \abs{\lambda_{\tau(j)} (A)-\lambda_j(B)}\leq (2J-1) (\norm{A}+\norm{B})^{(J-1)/J} \norm{A-B}^{1/J},\end{equation} where the minimum is over permutations $\tau$.
	\end{theorem}
\begin{lemma}
	\label[lemma]{lem:BasicApproximationTheoryResults}
		Let $A$ and $\hat{A}$ be matrices such that $A$ is invertible and $\norm{A-\hat{A}}<\norm{A^{-1}}^{-1}$. Then $\hat{A}$ is invertible and 
	\[  \norm{\hat{A}^{-1}-A^{-1}}\leq \frac{\norm{A^{-1}}^2\norm{A-\hat{A}}}{1-\norm{A^{-1}}\norm{A-\hat{A}}}.\]
	We also have the following:
	$\kappa(A)=\kappa(A^\transpose)$; $\sigma_J(A)=\sigma_J(A^\transpose)$; $\sigma_J(A)=\sigma_J(AW^\transpose)$ for any matrix $W$ whose columns are orthonormal and whose domain is $\RR^J$; $\sigma_J(AB)\geq \sigma_J(A)\sigma_J(B)$, and $\kappa(AB)\leq \kappa(A)\kappa(B)$ for matrices $A,B$.
\end{lemma}
\begin{proof}  
For the first see \cite[Chapter III, Theorem 2.5]{SS90} The other results can be found in Chapter I.4 of the same reference.
\end{proof}

\subsection{Sketch Proof of \cref{thm:DiscreteEstimation}}\label[appendix]{sec:ProofOfDiscreteEstimation}
	The arguments used to prove \cref{thm:LinftyEstimation} work also in this discrete setting, given the following observations and slight adaptations. To ease notation we assume that $f_j(x)=0$ for all $x\leq 0$ and $j\leq J$. We make the following definitions, which correspond to taking $h_l=\II_l,$ i.e.\ $h_l(x)=\II\braces{x=l}$, and replacing $K_L(x,y)$ by $\II\braces{x=y}$:
	\begin{align*}
		M^x = M^{x,L_0}&=\Pi_H(X_1=l,X_2=x,X_3=m)_{l,m\leq L_0}, \quad x\in \NN \\
		P=P^{L_0} &=\Pi_H(X_1=l,X_3=m)_{l,m\leq L_0},\\
		O=O^{L_0} &= \Pi_H(X_1=l \mid \theta_1=j)_{l\leq L_0,j\leq J},\\
		D^x &= (\diag \Pi_H(X_2=x \mid \theta_2=j)_{j\leq J})\equiv \diag ((O_{xj})_j).
	\end{align*}
	The proof of \cref{lem:AlgorithmIntuition} is unchanged with these adjusted definitions, and we adapt the definitions in \cref{alg:Linfty} correspondingly:
	\begin{align*}
		\hat{M}^x &=\brackets[\Big]{\frac{1}{N} \sum_{n\leq N} \II_l(X_n)\II_x(X_{n+1})\II_m(X_{n+2})}_{l,m\leq L_0},\\
		\hat{P} &=\brackets[\Big]{\frac{1}{N} \sum_{n\leq N} \II_l(X_n)\II_m(X_{n+2})}_{l,m\leq L_0},\\
		\hat{B}^x &= \brackets{\hat{V}^\transpose\hat{P}\hat{V}}^{-1} \hat{V}^\transpose \hat{M}^x\hat{V},
	\end{align*}
	for $\hat{V}$ comprising right singular vectors of $\hat{P}$. 
	
	Observe that the proofs of \cref{lem:ExistenceOfSuitableHl,lem:PfullRank} work in the current setting for the current choice of the $h_l$ [indeed, thanks to the disjoint support of $h_l,h_m$ for $l\neq m$ one can improve the bound in eq.~\cref{eqn:NormOBound} to $\norm{O^{L_0}}\leq J$], and similarly a version of \cref{lem:ExistenceOfHatAHatU} holds by choosing $\DD_N=\Aa_N\times \Uu_N$ for sequences of finite sets $\Aa_N\subset \RR,$ $\Uu_N\subset\NN$ such that $\cup_N \Uu_N=\NN$ and $\cup_N\Aa_N$ is dense in $\braces{a\in\RR : \sum_i \abs{a_i}\leq 1}$.
	
	Next note that a version of the Glivenko--Cantelli theorem gives control over $\sup_{x\in\NN} \norm{\hat{M}^x-M^x}$ for our new definitions of $\hat{M}^x,M^x$; we give here a slightly indirect proof of this fact by reusing the machinery of \cref{lem:ControlOfSupMx}. Indeed, inspecting the proof of \cref{lem:BracketingEntropyOfKernelTranslates}, one deduces that 
	\[ H_{[]}(\Tt, \norm{\cdot}_{L^2(\Pi_H)},\eps)\leq 16 L_0^2 \max(\eps^{-1},1)^4\] for $\Tt=\braces{\II_i\otimes \II_x\otimes \II_j : x\in \NN,~ i,j\leq L_0}$. It follows, in view of a standard bound (see \cref{eqn:IntSqrtLogBound}) \[ \int_0^x \sqrt{\log (1/u)}\du \leq x(1+\sqrt{\log(1/x)}),\]  and recalling as in \cref{lem:ControlOfP} that the chain  \[Y_n=(X_n,X_{n+1},X_{n+2},\theta_n,\theta_{n+1},\theta_{n+2})\] has pseudo-spectral gap bounded away from zero by \cref{ass:Qpi'}, that \cref{lem:ConcentrationResultsFromPaulin}, applied with $b=\sigma=1$, yields
	\[ \Pi\brackets[\big]{\sup_{x\in\NN} \abs{\hat{M}^x_{ij}-M^x_{ij}}>C(N^{-1/2}+N^{-1/2}\sqrt{u}+N^{-1}u)}\leq \exp(-u).\] We note that $\norm{\hat{M}^x-M^x}\leq L_0 \max_{ij} \abs{\hat{M}^x_{ij}-M^x_{ij}}$. Combining with \cref{lem:ControlOfP}, for any $c_N\to \infty$, we may choose suitable $L_0\to \infty$ and $u\to \infty$ to deduce
	\[ \Pi_H (\norm{\hat{P}-P}\leq c_N N^{-1/2},~\sup_{x\in\NN} \norm{\hat{M}^x-M^x}\leq c_N N^{-1/2})\to 1.\]
	The rest of the proof exactly mirrors that of \cref{thm:LinftyEstimation}.
	
	\section{Proof of the Lower Bound}\label[appendix]{sec:ProofOfLowerBound}
For the lower bound for simplicity we consider the (in view of the multiple testing application) most relevant case $J=2$. Let $\si_{2}$ denote the set of all permutations of $\braces{0,1}$. Define,  for $s, R>0$ and for the H\"older space $C^s$ defined in \cref{ass:smoothness},
\[ \cC^s(R)=\braces[\Big]{f\in C^s : f\geq 0,~\int_\RR f =1,~ \norm{f}_{C^s}\leq R}.\] 

\paragraphi{Parameters} The unknown parameters are $H=(Q,\pi,\bef)$, where $\bef=(f_0,f_1)$ denotes the vector of  emission densities.  Denoting by $P_{f_i}$ the distribution of density $f_i$ on $\RR$, $i=0,1$, the distribution of the observations $X=(X_1,\ldots,X_N)$ is
\[ \Pi_H=\Pi_H^{(N)} = \sum_{\bev\in \braces{0,1}^N} w_\bev\, \bigotimes_{j=1}^N P_{\bef_{v_j}}, \]
where $w_\bev$ denotes the probability under the Markov chain to observe the successive sequence of states $(v_1,\ldots,v_N)\in\braces{0,1}^N$, that is 
$w_{(v_1,\ldots,v_N)}=\pi_{v_1}Q_{v_1,v_2}\cdots Q_{v_{N-1},v_N}$.\\

\paragraphi{Class  $\Hh_{sep}$ of well--separated parameters}	
Let $\cF_{sep}$ be a class of pairs $\bef=(f_0,f_1)$ that are well--separated in the following sense, for a (small)  $d>0$ to be chosen:
\begin{equation} \label{fsep}
	\cF_{sep} = \left\{ \bef=(f_0,f_1) \in \cC^s(R):\  |(f_1-f_0)(0)|\ge d,\ 
	|P_{f_1}([-1,1])-P_{f_0}([-1,1])|\ge d \right\}.
\end{equation}
We define, for given $Q,\pi$,
\begin{equation}\label{defte}
	\Hh_{sep}=\Hh_{sep}(Q,\pi,R,d,s)=\left\{H=(Q,\pi,\bef):\ 	\bef\in\cF_{sep}\right\}.\\
\end{equation}

{\em Minimax risk.} 
For $\bef=(f_0,f_1)$ and $\beg=(g_0,g_1)$ two pairs of real functions, denote
\begin{equation}\label{rhod}
	\rho(\bef,\beg) = \min_{\vphi\in\sigma_2} \left( \|g_{\vphi(0)}-f_0\|_\infty + \|g_{\vphi(1)}-f_1\|_\infty \right).
\end{equation}
The loss $\rho$ is a pseudo--metric, verifying the axioms of a distance except that one can have $\rho(\bef,\beg)=0$ for $\bef\neq \beg$. We note that one could also consider the equivalent loss obtained by replacing the sum in \cref{rhod} with a maximum.

Let us consider the minimax risk 
\begin{equation} \label{minir}
	R_n=R_n(\Hh_{sep}) = \inf_{\bet=(T_0,T_1)} \sup_{H \in \Hh_{sep}}
	E_H \left[\rho\left(\bet,\bef\right)\right].
\end{equation}
Since $E[\min(X,Y)]\le \min(EX,EY)$, one notes that
\begin{equation} \label{miniris}
	R_n \le \inf_{\bet=(T_1,T_2)} \sup_{H \in H_{sep}}
	\left[  \min_{\vphi\in \sigma_2} \left(E_H
	\|T_{\vphi(1)} -f_1\|_\infty + E_H \|T_{\vphi(2)} -f_2\|_\infty \right) \right].
\end{equation}
In view of \cref{sec:Uniformity} (and constructing $\hat{f}_0,\hat{f}_1$ using $L_0=2$, $h_1=1$, $h_2=[-1,1]$ in \cref{alg:Linfty}), \Cref{thm:LinftyEstimation} provides a procedure for which the last quantity is bounded from above by (any rate slower than) $r_N=(N/\log N)^{-s/(2s+1)}$.  The next result provides the corresponding minimax lower bound. 
Note that the lower bound in \cref{thm:mmsimple} is pointwise in $Q$ and $\pi$, and thus continues to hold if $\pi,Q$ are allowed to vary in some set. 
\begin{proposition} \label[proposition]{thm:mmsimple}
	Consider $J=2$ classes, and fix both $\pi=(\pi_0,\pi_1)\in[0,1]^2$ and $Q$ a $2\times 2$ transition matrix. Given $s,R,d>0$, let $\Hh_{sep}$ be as in \eqref{defte} and let $R_n=R_n(\Hh_{sep})$ be as in \eqref{minir}. Then there exists $C=C(s,R)>0$ such that, for $N$ large enough, 
	\[
	R_n(\Hh_{sep}) \ge C \left(\frac{\log{N}}{N} \right)^{\frac{s}{2s+1}}.
	\]
\end{proposition}

\begin{proof} 
	We reduce the estimation problem to a classification problem in a standard way. 
	Suppose the two sets of densities $\{f_{0}^{(m)},\ 0\le m\le M\}$ and $\{f_1^{(m)},\ 0\le m\le M\}$ are such that for some $0<s_1,s_2<C_0$, 
	\begin{align} \label{c:sep12}
		\min\{ \|f_1^{(i)}-f_0^{(j)}\|_\infty,\ 0\le i,j\le M \} &\ge C_0,\\
		\label{c:sep12b}	\min\braces{ \norm{f_0^{(i)}-f_0^{(j)}} : 0\leq i,j\leq M, i\neq j}&\geq 2s_0,\\
		\label{c:sep12c}	\min\braces{ \norm{f_1^{(i)}-f_1^{(j)}} : 0\leq i,j\leq M, i\neq j}&\geq 2s_1.
	\end{align} 
	It follows that the family of functions $\bef^{(m)}=(f_0^{(m)},f_1^{(m)})$ is $2(s_0+s_1)$--separated in terms of $\rho$, since for $m\neq m'$,
	\begin{align*}
		\rho(\bef^{(m)},\bef^{(m')}) & \ge \min\left(\|f_0^{(m)}-f_{0}^{(m')}\|_\infty + 
		\|f_1^{(m)}-f_{1}^{(m')}\|_\infty,  
		\|f_1^{(m)}-f_{0}^{(m')}\|_\infty + 
		\|f_0^{(m)}-f_{1}^{(m')}\|_\infty \right)\\
		& \ge \min(2(s_0+s_1),2C_0)=2(s_0+s_1)=:2S.
	\end{align*}
	
	For a given estimator $\bet$ of $\bef\in\{\bef^{(0)},\ldots,\bef^{(M)}\}$, let $j^*(\bet)$ be the index $j$ such that  $\bef^{(j)}$ is the closest to $\bet$ in the $\rho$ pseudo-distance. Since the family $(\bef^{(m)},\, m\in\{0,\ldots,M\})$ is $2S$--separated, we have  $\rho(\bet,\bef^{(m)}) \ge S\1\{j^*(\bet)\neq m\}$. Writing $H_m=(Q,\pi,\bef^{(m)})$, we have
	\begin{align}
		\sup_{H \in \Hh_{sep}}
		E_H \left[\rho\left(\bet,\bef\right)\right] & \ge 
		\max_{0\le m\le M} 
		E_{H_m} \left[\rho\left(\bet,\bef^{(m)}\right)\right] \notag \\
		& \ge S \max_{0\le m\le M} \Pi_{H_m}[j^*(\bet)\neq m] \ge S p_{e,M}, \label{minlb}
	\end{align} 
	where $p_{e,M} = \inf_{\psi} \max_{0\le m\le M} \Pi_{H_m}[\psi\neq m],$ with the infimum being over all classifiers $\psi$. Taking the infimum with respect to $\bet$ in \cref{minlb}, one obtains $R_n(\Hh_{sep})\ge S p_{e,M}$.
	
	
	Lemma \ref{lem:pem} shows that in order to bound $p_{e,M}$ from below it suffices to bounds $\KL(\Pi_{H_m},\Pi_{H_0})$ from above, where $\KL(P,Q)$ denotes the Kullback-Leibler divergence between distributions $P$ and $Q$ with densities $p,q$,
	\begin{equation}\label{eqn:def:Kullback}
		\KL(P,Q)= E_P\sqbrackets[\Big]{\log\brackets[\Big]{\frac{p}{q}}}.
	\end{equation} By convexity of the map $(x,y)\to x\log(x/y)$, writing $\bev=(v_j)\in\braces{0,1}^N$, one obtains
	\[ \KL(\Pi_{H_m},\Pi_{H_0})\le \sum_{\bev\in\braces{0,1}^N} w_\bev 
	\KL\left(\bigotimes_{j=1}^N P_{f^{(m)}_{v_j}},\bigotimes_{j=1}^N P_{f^{(0)}_{v_j}}\right).\]
	For a given $\bev\in\braces{0,1}^N$, let $n_i(\bev)$, $i=0,1$, denote the number of elements of $\bev$ equal to $i$. The tensorisation property of the KL divergence implies
	\[ \KL\left(\bigotimes_{j=1}^N P_{f^{(m)}_{v_j}},\bigotimes_{j=1}^N P_{f^{(0)}_{v_j}}\right) = n_1(\bev)\KL(P_{f_1^{(m)}},P_{f_1^{(0)}}) + n_2(\bev)\KL(P_{f_2^{(m)}},P_{f_2^{(0)}}),\]
	where $n_1(\bev), n_2(\bev)$ are both at most $N$. 
	
	Let us now choose functions $f_0^{(m)}, f_1^{(m)}$, satisfying \cref{c:sep12,c:sep12b,c:sep12c} for which we have good control over $\KL(f_j^{(m)},f_j^{(0)})$, $j=0,1$ and $1\leq m\leq M$. For $\phi$ the standard normal density and $g_{m,A}$ defined as in \cref{lem:pert}, set,
	\begin{alignat*}{2}
		&f_0^{(m)}(x) = g_{m,A_0}(x), \quad &&m\geq 1,\qquad  f_0^{(0)}(x) =r\phi(r x),\qquad   \\
		&f_1^{(m)}(x)  = g_{m,A_1}(x-2/r),~ &&m\geq 1, \qquad f_1^{(0)}(x)  =r\phi(r(x-2/r)),
	\end{alignat*}
	where we choose 	\[ A_0 = c_0\left(\frac{\log{N}}{N} \right)^{\frac{s}{2s+1}}, \quad A_1 = c_1\left(\frac{\log{N}}{N} \right)^{\frac{s}{2s+1}},\quad
	M= d \ceil[\Big]{\brackets[\Big]{\frac{N}{\log{N}} }^{\frac{1}{2s+1}}}, \]
	with $r,c_0,c_1$ small, but fixed, positive constants. Note firstly that for $r,c_0,c_1$ small enough (and $N$ large enough) each pair $(f_0^{(m)},f_1^{(m)})$ is in $\cF_{sep}$ for some $d>0,R>0$. Indeed, examining the definition of $g_{m,A}$ from \cref{lem:pert}, we see for all $0\leq m \leq M$ that we have
	\[ \abs{f_1^{(m)}(0)-f_0^{(m)}(0)}=r\abs{\phi(2)-\phi(0)};\]
	that $P_{f_0^{(m)}}[-1,1]\ge r\phi(r/2)$; and, using $\int_x^\infty \phi(u)du=:\bar\Phi(x)\le\phi(x)/x$ for $x>0$ and $r<2$, that
	\[ P_{f_1^{(m)}}[-1,1] = \int_{-r-2}^{r-2}\phi\le \bar\Phi(2-r)\le \frac{\phi(2-r)}{2-r}.\]
	We further note by \cref{lem:pert} that for suitable $c_0,c_1,r,d,R$, we have both that \cref{c:sep12b,c:sep12c} hold for $s_0=A_0/2$, $s_1=A_1/2$, and also that
	\[ \KL\left(P_{f_0^{(m)}},P_{f_0^{(0)}}\right)\le C\frac{A_0^2}{M} \le \frac{C}{d} \frac{c_0^2 \log{N}}{N},\qquad
	\KL\left(P_{f_1^{(m)}},P_{f_1^{(m)}}\right)\le C\frac{A_1^2}{M} \le \frac{C}{d}\frac{c_1^2 \log{N}}{N}.
	\] 
	Putting the previous bounds together leads to
	\begin{align*} 
		\KL(\Pi_{H_m},\Pi_{H_0}) & \le 
		N\cdot \KL\left(P_{f_0^{(m)}},P_{f_0^{(0)}}\right)+N\cdot \KL\left(P_{f_1^{(m)}},P_{f_1^{(0)}}\right)
		\\
		& \le \frac{C\log{N}}{d}\left[ c_0^2 + c_1^2\right].
	\end{align*}
	In particular, one can bound from above
	\[ \frac1M\sum_{m=1}^M \KL(\Pi_{H_m},\Pi_{H_0}) 
	\le (C/d)(c_0^2+c_1^2) \log N \le (\log M)/10,
	\]
	provided $c_0, c_1$ are small enough constants, and we deduce by \cref{lem:pem} that $p_{e,M}:= \inf_\psi \max_{0\leq m\leq M} \Pi_{H_m}[\psi\neq m]$ is greater than a positive constant. Finally, recalling  \cref{minlb},we 
	\[ R_n(\Hh_{sep})\geq Sp_{e,M},\]
	with $S=2(s_0+s_1)=A_0+A_1$. The \namecref{prop:MinimaxLowerBound} follows by definition of $A_0, A_1$.
\end{proof}


\begin{lemma} \label[lemma]{lem:pert}
	Let $\psi$ be a $\cC^\infty$ function with support in $(-1/2,1/2)$ such that $\|\psi\|_\infty=1$ and $\int_\RR \psi=0$. Let $\phi(\cdot)$ denote the standard normal density and for $m\in\braces{1,\dots,M}$ and some $A, r>0$ and integer $M\ge 2$, set $g_0(x)=r\phi(rx)$ and
	\[ g_{m,A}(x) = r\phi(r x) + A \psi(Mx-m+1/2).\]
	Then for $s,R>0$, the functions $g_{m,A}$ are densities belonging to $\cC^s(R)$ provided $AM^s\le R/2$ and $r, A$ are small enough, 
	\[ \| g_{m,A}-g_{p,A}\|_\infty =A, \qquad (\text{for all }m\neq p),\]
	and, for $P_{g}$ the distribution with density $g$ on $\RR$, and some $C=C(r)>0$, any $m\in\braces{1,\dots,M}$,
	\[ \KL(P_{g_{m,A}},P_{g_{0,A}}) \le CA^2/M. \]
\end{lemma}
\begin{proof} 
	For the statement on supremum norms, it suffices to note that $x\to \psi(Mx-m)$ have disjoint support. For the KL--bounds, one expands the logarithm at the order $2$ in a neighborhood of $0$.
\end{proof}

\begin{lemma} \label[lemma]{lem:pem}
	For a family of points $(H_m)_{0\le m\le M}$  in $H$ with $M\ge 2$, let 
	\begin{equation} \label{def:pem}
		p_{e,M} = \inf_{\psi} \max_{0\le m\le M} \Pi_{H_m}[\psi\neq m], 
	\end{equation} 
	where the infimum is over all possible measurable $\psi$ taking values in $\braces{1,\dots,M}$. Suppose, for $\al<1/8$, 
	\[ \frac1M\sum_{m=1}^M \KL(\Pi_{H_m},\Pi_{H_0}) \le \al\log{M}.\]
	Then
	\[p_{e,M} \ge \frac{\sqrt{M}}{1+\sqrt{M}}\left(1-2\al-\sqrt{\frac{2\al}{\log{M}}}\right). \]
\end{lemma}
\begin{proof} 
	This follows from combining Proposition 2.3 and (the proof of) Theorem 2.5 in \cite{Tsybakov09}. 
\end{proof}

\section{Notation} \label[appendix]{sec:Notation}
We give notation assuming, as in \cref{sec:SupNormEstimation}, that there are a (known) number $J$ of hidden states $\braces{1,\dots, J}$ (recall that  $J=2$ for \cref{sec:EmpiricalBayesProcedure} and the proofs of results therein, with hidden states labelled 0 and 1, and the notation is adapted accordingly).\\
\paragraphi{HMM parameters}
\begin{description}
	\item [$X=(X_n)_{n\leq N}$] (or $(X_n)_{n\leq N+2}$ for convenience, or $(X_n)_{n\in \NN}$ for some of the proofs and lemmas) the data, drawn from the HMM \cref{eqn:def:model}.
	\item [$\theta=(\theta_n)_{n\leq N}$] the vector of hidden states, taking values in $\braces{1,\dots,J}^N$.
	\item [$Q,\pi$] the transition matrix of $\theta$ and its stationary (and initial) distribution.
	\item [$\mu$] a dominating measure on the space $\Xx=\RR$ (equipped with the usual Borel $\sigma$--algebra) in which $X_1$ takes values. Throughout we take $\mu$ to equal Lebesgue measure on $\RR$ or counting measure on $\ZZ\subset \RR$.
	\item [$f_1,\dots,f_J$] the emission densities, i.e.\ $f_j$ is the density of $X_1$ conditional on $\theta_1=j$.
	\item [$f_\pi$] the density of $X_1$; this is only used in the two-state case so $f_\pi=\pi_0f_0+\pi_1f_1$.
	\item [$H$]$=(Q,\pi,f_1,\dots,f_J)$, $\hat{H}=(\hat{Q},\hat{\pi},\hat{f}_1,\dots,\hat{f}_J)$.
	\item [$\Pi_H,E_H$] the law of $X$ for parameter $H$ and the associated expecation operator.
	\item [$\Hh,\Ii$:] see \cref{sec:Uniformity}. [Also note that $C=C(\Hh)$ is allowed to depend on the kernel $K$, the functions $(h_l)_{l\in \NN}$ and the sets $\DD_N$ since these can be chosen universally.]
	\item [$\nu,x^*$] constants as in \cref{ass:fpi}. 
	\item[$\delta$] a lower bound for $\min_{i,j}Q_{ij}$.
\end{description}
\paragraphi{Multiple testing}
\begin{description}
	\item [$\FDP,\FDR,\postFDR,\mFDR,\mTDR,$] see \cref{eqn:def:BFDR,eqn:def:mFDR,eqn:def:FDP,eqn:def:FDR,eqn:def:mTDR,eqn:def:postFDR1} (also \cref{eqn:postFDR} for an alternative characterisation of $\postFDR$).
	\item [$\ell_i$]$\equiv \ell_i(X)\equiv \ell_{i,H}(X)=\Pi_H(\theta_i=0 \mid X)$; $\hat{\ell}_i=\ell_{i,\hat{H}}$; $\ell_i'=\Pi_H(\theta_i=0\mid X_{i-A},\dots,X_{i+A})$ for some $A$; $\ell_i^\infty=\Pi_H(\theta_i=0 \mid (X_n)_{n\in \ZZ})$.
	\item[$\Phi_i^\infty$]$=\Pi_H(\theta_i=0 \mid (X_n : n\in\ZZ, n\leq i))$.
	\item [$\vphi_{\lambda,H}$]$=(\II\braces{\ell_{i,H}<\lambda})_{i\leq N}$
	\item [$\hat{\lambda}$]$=\sup\braces{\lambda : \postFDR_{\hat{H}}(\vphi_{\lambda,\hat{H}})\leq t}$.  $\lambda^*$ the solution to $E[\ell_i^\infty \mid \ell_i^\infty < \lambda^*]=\min(t,\pi_0)$.
	\item [$\hat{\vphi}$]$\equiv\hat{\vphi}^{(t)}=\vphi_{\hat{\lambda},\hat{H}}$ when there are no ties in $\ell$--values, and is defined by \cref{def:HatPhi} when there may be ties.
	\item [$\hat{S}_0$]$=\braces{i : \hat{\vphi}_i=1}$, $\hat{K}=\abs{\hat{S}_0}$.
	\end{description}
\paragraphi{Estimation}
\begin{description}
	\item[$h_1,\dots,h_{L_0}$], where $L_0$ is either constant or diverges slowly to infinity; bounded functions such that ``witness'' the linear independence of $f_1,\dots,f_J$ (see \cref{alg:Linfty} and \cref{lem:ExistenceOfSuitableHl}).
	\item[$K,K_L,$...] a convolution kernel, see \cref{eqn:def:KL}
	\item[$M^x$]$\equiv M^{x,L_0,L}=(E_H[h_i(X_1)K_L(x,X_2)h_j(X_3)]_{i,j\leq L_0})\in\RR^{L_0 \times L_0}$. 
	\item[$P$]$\equiv P^{L_0}=(E_H[h_i(X_1)h_j(X_3)]_{i,j\leq L_0})\in\RR^{L_0\times L_0}$
	\item[$O$]$=O^{L_0}=(E_H[h_i(X_1) \mid \theta_1=a]_{i\leq L_0, a\leq J})\in \RR^{L_0 \times J}$
	\item[$D$]$=D^x=\diag((K_L[f_j](x))_{j\leq J})$, i.e.\ the diagonal matrix whose diagonal entries are $D_{jj}=K_L[f_j](x)$.
	\item[$V$]$=V^{L_0}\in \RR^{L_0\times J}$ a matrix such that $V^\transpose PV$ is invertible. Specifically, we either take $V$ to equal a matrix of orthonormal right singular vectors of $P$ (so that $\sigma_J(V^\transpose P V)=\sigma_J(P)$) or, on the event of \cref{lem:EventA}, to equal $\hat{V}$ (defined in \cref{alg:Linfty}). 	\item[$B^x$]$=B^{x,L_0}=[V^\transpose PV]^{-1}   V^\transpose M^xV\equiv [QO^\transpose V]^{-1} D^x QO^\transpose V$.
	\item[$\hat{M}^x,\hat{P},\hat{O},\hat{V}$] empirical versions of $M^x,P,O,V,B^x$ (see \cref{alg:Linfty}, p\pageref{alg:Linfty}).
	\item[$\hat{B}^x$]$=[\hat{V}^\transpose \hat{P}\hat{V}]^{-1}  \hat{V}^\transpose \hat{M}^x \hat{V}$, $\tilde{B}^x=[\hat{V}^\transpose P\hat{V}]^{-1} \hat{V}^\transpose M^x \hat{V}$, $\hat{B}^{a,u}=\sum a_i \hat{B}^{u_i}$ and $\tilde{B}^{a,u}=\sum a_i \tilde{B}^{u_i}$ for $a,u\in \RR^{J(J-1)/2}$ such that $\sum\abs{a_i}\leq 1$.
	\item[$\sep(B)$]$=\min_{i\neq j}\abs{\lambda_i-\lambda_j}$ the ``eigen-separation'' of a matrix $B\in \RR^{J\times J}$, with eigenvalues $\lambda_1,\dots,\lambda_J$.
	\item[$\hat{a},\hat{u},\DD_N$] See \cref{alg:Linfty}, p\pageref{alg:Linfty}.
	\item[$\hat{R}$] a matrix of normalised columns diagonalising $\hat{B}^{\hat{a},\hat{u}}$, $\tilde{R}$ a matrix whose columns are those of $QO^\transpose \hat{V}$ but scaled to have unit Euclidean norm (which therefore diagonalises $\tilde{B}^{a,u}$ for any $a,u$).
	\item [$\Aa$]$= \braces{ \norm{\hat{P}-P}\leq cL_0 r_N,~ \norm{\hat{M}^x-M^x}\leq cL_0^2 r_N ~\forall x\in \RR}$ the event of \cref{lem:EventA}.
		\item[$C^s$] the usual H\"older space (see \cref{ass:smoothness} equipped with the usual norm $\norm{\cdot}_{C^s}$.
\end{description}
\paragraphi{Minimax lower bound}
\begin{description}
	\item[$\cC^s(R)$] the subspace of $C^s$ consisting of densities with H\"older norm bounded by $R$.
	\item[$\sigma_2$] the set of all permutations on $\braces{0,1}$.
	\item[$\rho(\bef,\beg)$]$= \min_{\vphi\in\sigma_2} \brackets[\big]{ \|g_{\vphi(0)}-f_0\|_\infty + \|g_{\vphi(1)}-f_1\|_\infty },$ for $\bef=(f_0,f_1)$, $\beg=(g_0,g_1)$.
	\item[$\cF_{sep}$]$= \braces[\big]{ \bef=(f_0,f_1) \in \cC^s(R):\  |(f_1-f_0)(0)|\ge d,\ 
	|P_{f_1}([-1,1])-P_{f_0}([-1,1])|\ge d }$.
\item[$\Hh_{sep}$]$=\left\{H=(Q,\pi,\bef):\ 
	\bef\in\cF_{sep}\right\},$ for some arbitrary (fixed) $Q,\pi$.
\end{description}
\paragraphi{Miscellaneous}
\begin{description}
	\item[$\norm{\cdot},\norm{\cdot}_F,\norm{\cdot}_\infty$] the Euclidean norm on vectors or the corresponding operator norm on matrices, the Frobenius norm on matrices, and the $L^\infty$ (supremum) norm on functions taking values in $\RR$.
	\item[$\sigma_j(A)$] the $j$th largest singular value of a matrix $A$.
	\item[$\kappa(A)$]= $\sigma_1(A)/\sigma_J(A)=\norm{A}\norm{A^{-1}}$ for a matrix with smaller dimension $J$, the condition number of the matrix $A$.
	\item[$o(1),o_p(1)$] The usual little-oh notation: $a_N=o(1)$ if $a_N\to 0$ as $N\to infty$, $a_N=o_p(1)$ if $a_N\to 0$ in probability as $N\to \infty$.
	\item[$C^s(\RR)$] the usual space of locally H\"older smooth functions, equipped with the usual H\"older norm $\norm{\cdot}_{C^s(\RR)}$ (see \cref{ass:smoothness}). Note that since we consider density functions, we could equivalently use the space of globally H\"older smooth functions.
	\item [$r_N$]$=(N/\log N)^{-s/(1+2s)}$. $\eps_N$ some rate of consistency of estimators in \cref{eqn:ConsistencyAssumption}.
	\item [$N_{[]},H_{[]}$:] The bracketing numbers/entropy, wherein $N_{[]}(\Tt,\norm{\cdot}_{L^2(P)},\eps)$ is the smallest number of pairs of functions $(\underline{f},\bar{f})$ such that every $g\in \Tt$ is bracketed by one of the pairs, where $(\underline{f},\bar{f})$ brackets $g$ if $\underline{f}\leq g\leq \bar{f}$ pointwise, and $H_{[]}(\Tt,\norm{\cdot}_{L^2(P)},\eps):=\log N_{[]}(\Tt,\norm{\cdot}_{L^2(P)},\eps).$ 
\end{description}

\bibliography{bibliography}
\end{document}